\documentclass[12pt,oneside,reqno]{amsart}
\usepackage{mathrsfs}
\usepackage{graphics}
\usepackage{cite}
\usepackage{amssymb}
\usepackage{enumerate}
\newcommand{\dif}{\mathrm{d}}

\newcommand{\be}{\begin{eqnarray}}
\newcommand{\ee}{\end{eqnarray}}
\newcommand{\ce}{\begin{eqnarray*}}
\newcommand{\de}{\end{eqnarray*}}
\newtheorem{theorem}{Theorem}[section]
\newtheorem{lemma}[theorem]{Lemma}
\newtheorem{remark}[theorem]{Remark}
\newtheorem{definition}[theorem]{Definition}
\newtheorem{proposition}[theorem]{Proposition}
\newtheorem{Example}[theorem]{Example}
\newtheorem{corollary}[theorem]{Corollary}
\newtheorem{condition}[theorem]{Condition}
\def\e{\varepsilon}
\def\s{\sigma}
\def\t{\theta}
\def\a{\alpha}

\def\b{\beta}
\def\d{\delta}
\def\p{\partial}
\def\g{\gamma}
\def\l{\lambda}

\def\[{{\Big[}}
\def\]{{\Big]}}
\def\<{{\langle}}
\def\>{{\rangle}}
\def\({{\Big(}}
\def\){{\Big)}}

\def\no{\nonumber}
\def\bt{\begin{theorem}}
\def\et{\end{theorem}}
\def\bl{\begin{lemma}}
\def\el{\end{lemma}}
\def\br{\begin{remark}}
\def\er{\end{remark}}
\def\bx{\begin{Example}}
\def\ex{\end{Example}}
\def\bd{\begin{definition}}
\def\ed{\end{definition}}
\def\bp{\begin{proposition}}
\def\ep{\end{proposition}}
\def\bc{\begin{corollary}}
\def\ec{\end{corollary}}
\def\bco{\begin{condition}}
\def\eco{\end{condition}}

\def\cB{{\mathcal B}}

\def\cD{{\mathcal D}}

\def\cN{{\mathcal N}}
\def\cO{{\mathcal O}}

\def\mE{{\mathbb E}}

\def\mH{{\mathbb H}}

\def\mN{{\mathbb N}}

\def\mP{{\mathbb P}}

\def\mR{{\mathbb R}}
\def\mS{{\mathbb S}}

\def\sA{{\mathscr A}}
\def\sB{{\mathscr B}}

\def\sF{{\mathscr F}}

\def\sV{{\mathscr V}}

\def\geq{\geqslant}
\def\leq{\leqslant}
\def\epsilon{\varepsilon}

\begin{document}

\allowdisplaybreaks

\title{Asymptotic behaviors of multiscale multivalued stochastic systems with small noises}

\author{Huijie Qiao}

\dedicatory{School of Mathematics,
Southeast University,\\
Nanjing, Jiangsu 211189, P.R.China\\
hjqiaogean@seu.edu.cn}

\thanks{{\it AMS Subject Classification(2020):} 60H10, 60F10, 60F15}

\thanks{{\it Keywords: Multiscale multivalued stochastic systems with small noises; averaging principles; large deviation principles}}

\thanks{This work was partly supported by NSF of China (No.12071071).}

\subjclass{}

\date{}

\begin{abstract}
In this paper, we consider asymptotic behaviors of multiscale multivalued stochastic systems with small noises. First of all, for general, fully coupled systems for multivalued stochastic differential equations of slow and fast motions with small noises in the slow components, we prove an averaging principle in the strong convergence sense. Moreover, a convergence rate is given in a special case. Next, for these systems, we establish the large deviation principle by the weak convergence approach. Then for a special case the rate function is explicitly characterized. Finally, we explain our results by an example.
\end{abstract}

\maketitle \rm

\section{Introduction}

Given a complete probability space $(\Omega,\sF,\mP)$. Let $(W^1_{t}), (W^2_{t})$ be mutually independent $d_1$- and $d_2$-dimensional standard Brownian motions defined on $(\Omega,\sF,\mP)$, respectively. Let $\{\sF_{t}\}_{t\in[0,\infty)}$ be the natural filtration generated by $(W^1_{t})$ and $(W^2_{t})$. Consider the following slow-fast system on $\mR^{n} \times \mR^{m}$:
\be\left\{\begin{array}{l}
\dif X_{t}^{\e,\g}\in -A_1(X_{t}^{\e,\g})\dif t+b_{1}(t,X_{t}^{\e,\g},Y_{t}^{\e,\g})\dif t+\sqrt{\e}\s_{1}(t,X_{t}^{\e,\g},Y_{t}^{\e,\g})\dif W^1_{t},\\
X_{0}^{\e,\g}=x_0\in\overline{\cD(A_1)},\quad  t\geq 0,\\
\dif Y_{t}^{\e,\g}\in -A_2(Y_{t}^{\e,\g})\dif t+\frac{1}{\g}b_{2}(t,X_{t}^{\e,\g},Y_{t}^{\e,\g})\dif t+\frac{1}{\sqrt{\g}}\s_{2}(t,X_{t}^{\e,\g},Y_{t}^{\e,\g})\dif W^2_{t},\\
Y_{0}^{\e,\g}=y_0\in\overline{\cD(A_2)},\quad  t\geq 0,
\end{array}
\right.
\label{Eq1i}
\ee
where $A_1, A_2$ are two multivalued maximal monotone operators with $\text{Int}(\cD(A_1))\ne\emptyset,\\ \text{Int}(\cD(A_2))\ne\emptyset$ (See Subsection \ref{mmo}), and these mappings $b_{1}: [0,\infty)\times\mR^{n}\times\mR^{m}\rightarrow\mR^{n}$, $\s_{1}: [0,\infty)\times\mR^{n}\times\mR^{m}\rightarrow\mR^{n\times d_1}$, $b_{2}: [0,\infty)\times\mR^{n}\times\mR^{m}\rightarrow\mR^{m}$, $\s_{2}: [0,\infty)\times\mR^{n}\times\mR^{m}\rightarrow\mR^{m\times d_2}$ are all Borel measurable. Here $0<\e<1$ is a small parameter and $0<\g<1$ represents the other small parameter which characterizes the ratio of timescales between processes $X_{\cdot}^{\e,\g}$ and $Y_{\cdot}^{\e,\g}$. 

Here we call the system (\ref{Eq1i}) a multiscale multivalued stochastic system with small noises, where $X_{\cdot}^{\e,\g}$, $Y_{\cdot}^{\e,\g}$ are slow and fast components, respectively. This type of systems arises from physics, engineering and economics and is usually used to solve optimization and equilibrium problems involving random data; see \cite{bmy, bt1, bt2, cepa1, cepa2, rw} and the references therein. For these systems, the averaging principle is an effective tool in capturing the main behaviors of the slow components while avoiding the complexities caused by the details of the fast components of them. If $A_1=A_2=0$ and $\e=1$, the system (\ref{Eq1i}) is just a common multiscale system or slow-fast system. Moreover, there have been many related results (cf. e.g. \cite{br, ee,elve, fw1, fw2,  lrsx, kh, yk, ck, av1}). However, for the case of $A_1\neq 0$ or $A_2\neq 0$, there are few related results about the averaging principle for the system (\ref{Eq1i}). Hence, the first goal of this paper is to show an averaging principle for the system (\ref{Eq1i}). Concretely speaking, under some suitable assumptions we prove that if the parameters $\e, \g$ satisfy 
\ce
\lim_{\e\rightarrow 0}\frac{\g}{\e}=\left\{\begin{array}{l} 0,\\
\varrho\in (0,\infty), 
\end{array}
\right.
\de
then
\be
\lim\limits_{\e\rightarrow 0}\mE\(\sup_{0\leq t\leq T}|X_{t}^{\e,\g}-\bar{X}_{t}|^{2}\)=0,
\label{averprin0}
\ee
where $\bar{X}$ solves the following equation
\ce\left\{\begin{array}{l}
\dif \bar{X}_{t}\in -A_1(\bar{X}_{t})\dif t+\bar{b}_{1}(t,\bar{X}_{t})\dif t,\\
\bar{X}_{0}=x_0\in\overline{\cD(A_1)},\quad  0\leq t\leq T,
\end{array}
\right.
\de
$\bar{b}_{1}(t,x):=\int_{\mR^{m}}b_{1}(t,x,y)\nu^{t,x}(\dif y)$ and $\nu^{t,x}$ is the unique invariant probability measure of Eq.(\ref{Eq2}) (See Subsection \ref{fun}). Moreover, a convergence rate is given in a special case. In order to prove two averaging principles, we will face three difficulties. Firstly, since the operator $A_1$ is multivalued, nonlinear and not smooth, it is difficult to  directly estimate $|X_{t+l}^{\e,\g}-X_{t}^{\e,\g}|$, which is crucial in the time discretization method. Here we use the properties of maximal monotone operators to deal with it. Secondly, some finite variation processes appear and we will use their total variation estimates. Usually it is difficult to straight calculate these total variations. Hence, we make use of some inequalities to treat these total variations. Thirdly, dependence on the time variable and maximal monotone operators makes some estimations more complicated and difficult than that in common cases. Therefore, we need some new ideas and approaches so as to reach our goals.

Recently, Chen and Wu \cite{cw} considered the following multiscale stochastic system:
\be\left\{\begin{array}{l}
\dif X_{t}^{\g}\in -\p \varphi(X_{t}^{\g})\dif t+\tilde b_{1}(X_{t}^{\g},Y_{t}^{\g})\dif t+\tilde \s_{1}(X_{t}^{\g})\dif W^1_{t},\\
X_{0}^{\g}=x_0\in\overline{\cD(\p \varphi)},\quad  0\leq t\leq T,\\
\dif Y_{t}^{\g}\in -\p \psi(Y_{t}^{\g})\dif t+\frac{1}{\g}\tilde b_{2}(X_{t}^{\g},Y_{t}^{\g})\dif t+\frac{1}{\sqrt{\g}}\tilde\s_{2}(X_{t}^{\g},Y_{t}^{\g})\dif W^2_{t},\\
Y_{0}^{\g}=y_0\in\overline{\cD(\p \psi)},\quad  0\leq t\leq T,
\end{array}
\right.
\label{Eq1p}
\ee
where $\varphi, \psi$ are lower semicontinuous convex functions on $\mR^n, \mR^m$, respectively, $\p \varphi, \p \psi$ are corresponding subdifferential operators  and $\tilde b_{1}:\mR^{n}\times\mR^{m}\rightarrow\mR^{n}, \tilde\s_{1}:\mR^{n}\rightarrow\mR^{n\times d_1}, \tilde b_{2}: \mR^{n}\times\mR^{m}\rightarrow\mR^{m}, \tilde\s_{2}:\mR^{n}\times\mR^{m}\rightarrow\mR^{m\times d_2}$ are Borel measurable. There they also obtained the similar result. Note that by {\bf Example \ref{exmmo1}}, $\p \varphi, \p \psi$ are maximal monotone operators and the system (\ref{Eq1p}) is a multiscale multivalued stochastic system. Thus, comparing Theorem \ref{xbarxp} (See Subsection \ref{fun}) with Theorem 4.3 in \cite{cw}, we find that all our coefficients depend on the time variable, and our $\s_1$ can depend on $y$. That is, the small noise in the slow component can improve the result, which is one of our motivations. 

Next, it would be useful to further study the probability of the deviation $X_{\cdot}^{\e,\g}-\bar{X}_{\cdot}$. The study of these rare events in the multiscale framework is a difficult problem due to the presence of the underlying fast component. The first  step is to develop the associated large deviation theory. 

The large deviation principle (LDP in short), which describes the asymptotic behaviors of the trajectory $X_{t}^{\e,\g}-\bar{X}_{t}$ for $0\leq t\leq T$ as $\e\rightarrow 0$ and $\g\rightarrow 0$, is an important topic in the probability theory and has widespread applications in many fields such as thermodynamics, statistics, information theory and engineering. Moreover, the research on the LDP for the usual multiscale stochastic systems has fruitful results in the literature (cf. \cite{ds, ghl, rl, kp, ku, aap, ks1, av2, av3}). As to the case in this paper, to the best of our knowledge there are few related results in the literature. Therefore, the second goal of this paper is to show the LDP for the following multiscale multivalued stochastic system:
\be\left\{\begin{array}{l}
\dif X_{t}^{\e,\g}\in -A_1(X_{t}^{\e,\g})\dif t+b_{1}(X_{t}^{\e,\g},Y_{t}^{\e,\g})\dif t+\sqrt{\e}\s_{1}(X_{t}^{\e,\g})\dif W^1_{t},\\
X_{0}^{\e,\g}=x_0\in\overline{\cD(A_1)},\quad  0\leq t\leq T,\\
\dif Y_{t}^{\e,\g}\in -A_2(Y_{t}^{\e,\g})\dif t+\frac{1}{\g}b_{2}(X_{t}^{\e,\g},Y_{t}^{\e,\g})\dif t+\frac{1}{\sqrt{\g}}\s_{2}(X_{t}^{\e,\g},Y_{t}^{\e,\g})\dif W^2_{t},\\
Y_{0}^{\e,\g}=y_0\in\overline{\cD(A_2)},\quad  0\leq t\leq T.
\end{array}
\right.
\label{Eq1ldp}
\ee
Concretely speaking, if the parameters $\e, \g$ satisfy 
$$
\lim_{\e\rightarrow 0}\frac{\g}{\e}=0,
$$
we prove that the family $\{X^{\epsilon,\g},\epsilon\in(0,1)\}$ satisfies the LDP in $\mS:=C([0,T],\overline{\mathcal{D}(A_1)})$. Comparing Theorem \ref{ldpmmsde} with \cite[Theorem 3.3 (Regime 1)]{ks1}, we find that our conditions are more general. Furthermore, we explicitly characterize the rate function in a special case, which in general can't be done (\cite{ds, flqz, ghl}). 

In \cite{ku}, Kushner considered the following multiscale stochastic system:
\be\left\{\begin{array}{l}
\dif \tilde X_{t}^{\e,\g}=\tilde b_{1}(\tilde X_{t}^{\e,\g},\tilde Y_{t}^{\e,\g})\dif t+\sqrt{\e}\tilde\s_{1}(\tilde X_{t}^{\e,\g})\dif W^1_{t}+\dif \tilde Z_{t}^{\e,\g},\\
\tilde X_{0}^{\e,\g}=\tilde x_0,\quad  0\leq t\leq T,\\
\dif \tilde Y_{t}^{\e,\g}=\frac{1}{\g}\tilde b_{2}(\tilde X_{t}^{\e,\g},\tilde Y_{t}^{\e,\g})\dif t+\frac{1}{\sqrt{\g}}\hat\s_{2}(\tilde X_{t}^{\e,\g})\dif W^2_{t},\\
\tilde Y_{0}^{\e,\g}=\tilde y_0,\quad  0\leq t\leq T,
\end{array}
\right.
\label{Eq1ldpc}
\ee
where $\tilde Z_{\cdot}^{\e,\g}$ is a boundary reflection process depending on $\tilde X_{\cdot}^{\e,\g}$ and $\hat\s_{2} :\mR^{m}\rightarrow\mR^{m\times d_2}$ is Borel measurable. There he mentioned that when 
$$
\lim_{\e\rightarrow 0}\frac{\g}{\e}=0,
$$
and $b_{1}, b_2$ are linear, the family $\{\tilde X^{\epsilon,\g},\epsilon\in(0,1)\}$ satisfies the LDP. Note that by {\bf Example \ref{exmmo2}}, if we take $A_1=\p I_{\cO}$ and $A_2=0$, where $\cO$ is a closed convex subset of $\mathbb{R}^n$, the system (\ref{Eq1ldp}) can become the system (\ref{Eq1ldpc}). Therefore our result (Theorem \ref{ldpmmsde}) can cover this result.

Here we mention that based on the order that $\e, \g$ converge to $0$, there are three different regimes, i.e.
\ce
\lim_{\e\rightarrow 0}\frac{\g}{\e}=\left\{\begin{array}{l} 0, \qquad\qquad\qquad~~ ~\mbox{Regime}~ 1;\\
\varrho\in (0,\infty), ~\mbox{Regime}~  2;\\
\infty, ~\qquad\quad\quad~\mbox{Regime}~ 3.
\end{array}
\right.
\de
For Regime $2$ and Regime $3$, the LDP of current setting becomes quite complicated and our current approach does not seem to work directly anymore. In \cite{ks1}, Spiliopoulos investigated the LDP for all three regimes by means of occupation measures and the characterization of optimal controls. However, for the system (\ref{Eq1ldp}), due to containing the multivalued operators $A_1, A_2$, it is difficult to express the needed Poisson equations, which is important to prove the tightness for a family of occupation measures. Therefore, in the forthcoming work we would consider the cases of Regime $2$ and Regime $3$.

The rest of this paper is organized as follows. In Section \ref{pre}, we introduce some notations, maximal monotone operators and multivalued stochastic differential equations (SDEs for short). In Section \ref{mare}, we formulate the main results. Before proving main results, we make preparations in Section \ref{prep}. Main results are showed in Section \ref{proofirs} and \ref{prooseco}, respectively. Finally, we explain our results through an example in Section \ref{exam}.

The following convention will be used throughout the paper: $C$ with or without indices will denote different positive constants whose values may change from one place to another.

\section{Notations and assumptions}\label{pre}

In the section, we will recall some notations and some concepts of maximal monotone operators and multivalued SDEs.

\subsection{Notations}\label{nn}
In this subsection, we introduce some notations used in the sequel.

Let $|\cdot|, \|\cdot\|$ be the norms of a vector and a matrix, respectively. Let $\langle\cdot,\cdot\rangle$ be the inner product of vectors on $\mR^n$. $U^{*}$ denotes the transpose of the matrix $U$.

Let $\cB_b(\mR^{n})$ be the set of all bounded Borel measurable functions on $\mR^n$. Let $C(\mR^n)$ be the set of all  functions which are continuous on $\mR^n$. $C^{2}(\mR^n)$ represents the collection of all functions in $C(\mR^n)$ with continuous derivatives of order up to 2. $C_b^2(\mR^n)$ stands for the subspace of $C^2(\mR^n)$, consisting of functions whose derivatives up to order 2 and themselves are bounded.

\subsection{Maximal monotone operators}\label{mmo}

In the subsection, we introduce maximal monotone operators.

For a multivalued operator $A: \mR^n\mapsto 2^{\mR^n}$, where $2^{\mR^n}$ stands for all the subsets of $\mR^n$, set
\ce
&&\cD(A):= \left\{x\in \mR^n: A(x) \ne \emptyset\right\},\\
&&Gr(A):= \left\{(x,y)\in \mR^{2n}:x \in \cD(A), ~ y\in A(x)\right\}.
\de
We say that $A$ is monotone if $\langle x_1 - x_2, y_1 - y_2 \rangle \geq 0$ for any $(x_1,y_1), (x_2,y_2) \in Gr(A)$, and $A$ is maximal monotone if 
$$
(x_1,y_1) \in Gr(A) \iff \langle x_1-x_2, y_1 -y_2 \rangle \geq 0, \quad \forall (x_2,y_2) \in Gr(A).
$$

For readers to understand maximal monotone operators very well, we give two examples.

\bx\label{exmmo1}
For a lower semicontinuous convex function $\varphi:\mR^n\mapsto(-\infty, +\infty]$, we assume ${\rm Int}(Dom(\varphi))\neq \emptyset$, where $Dom(\varphi)\equiv\{x\in\mR^n; \varphi(x)<\infty\}$ and $\operatorname{Int}(Dom(\varphi))$ is the interior of $Dom(\varphi)$. Define the subdifferential operator of the function $\varphi$:
$$
\partial\varphi(x):=\{y\in\mR^n: \<y,z-x\>+\varphi(x)\leq \varphi(z), \forall z\in\mR^n\}.
$$
Then $\partial\varphi$ is a maximal monotone operator. 
\ex

\bx\label{exmmo2}
For a closed convex subset $\mathcal{O}$ of $\mathbb{R}^n$, we suppose $\operatorname{Int}(\mathcal{O})\neq\emptyset$. Define the indicator function of $\mathcal{O}$ as follows:
$$
I_{\mathcal{O}}(x):= \begin{cases}0, & \text { if } x \in \mathcal{O}, \\ 
+\infty, & \text { if } x \notin \mathcal{O}.\end{cases}
$$
The subdifferential operator of $I_{\mathcal{O}}$ is given by
$$
\begin{aligned}
\partial I_{\mathcal{O}}(x) & :=\left\{y \in \mathbb{R}^n:\langle y, z-x\rangle \leq 0, \forall z \in \mathcal{O}\right\} \\
& = \begin{cases}\emptyset, & \text { if } x \notin \mathcal{O}, \\
\{0\}, & \text { if } x \in \operatorname{Int}(\mathcal{O}), \\
\Lambda_x, & \text { if } x \in \partial \mathcal{O},\end{cases}
\end{aligned}
$$
where $\Lambda_x$ is the exterior normal cone at $x$. By simple deduction, we know that $\partial I_{\mathcal{O}}$ is a maximal monotone operator.
\ex
In the following, we recall some properties of a maximal monotone operator $A$ (cf.\cite{cepa1,q1}):
\begin{enumerate}[(i)]
\item
${\rm Int}(\cD(A))$ and $\overline{\mathrm{\cD}(A)}$ are convex subsets of $\mR^n$ with ${\rm Int}\left( \overline{\mathrm{\cD}(A)} \right) = {\rm Int}\( \mathrm{\cD}(A) \) 
$, where ${\rm Int}(\cD(A))$ denotes the interior of the set $\cD(A)$. 
\item For every $x\in\mR^n$, $A(x)$ is a closed and convex subset of $\mR^n$.
\end{enumerate}

Take any $T>0$ and fix it. Let $\sV_{0}$ be the set of all continuous functions $K: [0,T]\mapsto\mR^n$ with finite variations and $K_{0} = 0$. For $K\in\sV_0$ and $s\in [0,T]$, we shall use $|K|_{0}^{s}$ to denote the variation of $K$ on $[0,s]$ and write $|K|_{TV}:=|K|_{0}^{T}$. Set
\ce
&&\sA:=\Big\{(X,K): X\in C([0,T],\overline{\cD(A)}), K \in \sV_0, \\
&&\qquad\qquad\quad~\mbox{and}~\langle X_{t}-x, \dif K_{t}-y\dif t\rangle \geq 0 ~\mbox{for any}~ (x,y)\in Gr(A)\Big\}.
\de
And about $\sA$ we have the following two results (cf.\cite[Proposition 4.1 and 4.4]{cepa2} or \cite[Proposition 3.3 and 3.4]{ZXCH}).

\bl\label{equi}
For $X\in C([0,T],\overline{\cD(A)})$ and $K\in \sV_{0}$, the following statements are equivalent:
\begin{enumerate}[(i)]
	\item $(X,K)\in \sA$.
	\item For any $(x,y)\in C([0,T],\mR^n)$ with $(x_t, y_t)\in Gr(A)$, it holds that 
	$$
	\left\langle X_t-x_t, \dif K_t-y_t\dif t\right\rangle \geq0.
	$$
	\item For any $(X^{'},K^{'})\in \sA$, it holds that 
	$$
	\left\langle X_t-X_t^{'},\dif K_t-\dif K_t^{'}\right\rangle \geq0.
	$$
\end{enumerate}
\el

\bl\label{inteineq}
Assume that $\text{Int}(\cD(A))\ne\emptyset$. For any $a\in \text{Int}(\cD(A))$, there exist $M_1 >0$, and $M_{2},M_{3}\geq0$ such that  for any $(X,K)\in \sA$ and $0\leq s<t\leq T$,
$$
\int_s^t{\left< X_r-a, \dif K_r \right>}\geq M_1\left| K \right|_{s}^{t}-M_2\int_s^t{\left| X_r-a\right|}\dif r-M_3\left( t-s \right) .
$$
\el

\subsection{Multivalued SDEs}

In the subsection, we introduce multivalued SDEs. 

Consider the following multivalued SDE on $\mR^n$:
\be
\dif X_t\in -A(X_t)\dif t+b(t,X_t)\dif t+\s(t,X_t)\dif W^1_t, \quad t\geq 0, \label{eq1}
\ee
where $A$ is a maximal monotone operator with $\text{Int}(\cD(A))\ne\emptyset$, the coefficients $b: [0,\infty)\times\mR^n\mapsto{\mR^n}, \,\,\sigma: [0,\infty)\times\mR^n\mapsto{\mR^n}\times{\mR^{d_1}}$ are Borel measurable and $W^1_{\cdot}$ is a $d_1$-dimensional Brownian motion on the filtered probability space $(\Omega, \mathscr{F}, \{\mathscr{F}_t\}_{t\in[0,\infty)}, \mP)$.

\bd\label{strosolu}
We say that Eq.$(\ref{eq1})$ admits a strong solution with the initial value $X_0\in\overline{\cD(A)}$ if there exists a pair of adapted processes $(X,K)$ on $(\Omega, \mathscr{F}, \{\mathscr{F}_t\}_{t\in[0,\infty)}, \mP)$ such that

$(i)$ $X_t\in{\mathscr{F}_t^{W^1}}$, where $\{\mathscr{F}_t^{W^1}\}_{t\in[0,\infty)}$ stands for the $\sigma$-field filtration generated by $W^1$,

$(ii)$ $(X_{\cdot}(\omega),K_{\cdot}(\omega))\in \sA$ $\mP$ a.s.,

$(iii)$ it holds that for any $T>0$
\ce
\mP\left\{\int_0^T(|b(s,X_s)|+\|\sigma(s,X_s)\|^2)\dif s<+\infty\right\}=1,
\de
and
\ce
X_t=X_0-K_{t}+\int_0^tb(s,X_s)\dif s+\int_0^t\sigma(s,X_s)\dif W^1_s, \quad 0\leq{t}\leq{T}, \quad \mP~a.s..
\de
\ed

\subsection{A general criterion of LDPs}

In this subsection, we present a general criterion to establish the LDP. 

Let $\mS$ be a Polish space. For each $\e>0$, let $X^{\e}$ be a $\mS$-valued random variable given on $(\Omega, \mathscr{F}, \{\mathscr{F}_t\}_{t\in[0,\infty)}, \mP)$.

\bd\label{compleve}
The function $I$ on $\mathbb{S}$ is called a rate function if for each $M<\infty$, $\{\varsigma\in \mathbb{S}:I(\varsigma)\leq M\}$ is a compact subset of $\mathbb{S}$.
\ed

\bd
$(i)$ We say that the family of $\{X^{\epsilon}, \e\in(0,1)\}$ satisfies the Laplace principle in $(\mS,\sB(\mS))$ with the rate function $I$,  if for any real bounded continuous function $g$ on $\mathbb{S}$,
\ce
\lim\limits_{\varepsilon\rightarrow 0}\varepsilon \log \mE\left\{\exp\left[-\frac{g(X^{\epsilon})}{\epsilon}\right]\right\}=-\inf\limits_{\varsigma\in \mathbb{S}}\(g(\varsigma)+I(\varsigma)\).
\de

$(ii)$ We say that the family of $\{X^{\epsilon}, \e\in(0,1)\}$ satisfies the LDP in $(\mS,\sB(\mS))$ with the rate function $I$, if for any closed subset $B_1\in \sB(\mS)$,
$$
\limsup\limits_{\e\rightarrow 0}\e\log\mP(X^{\epsilon}\in B_1)\leq -\inf\limits_{\varsigma\in B_1}I(\varsigma),
$$
and for any open subset $B_2\in \sB(\mS)$,
\ce
\liminf_{\e\rightarrow 0}\e\log\mP(X^{\epsilon}\in B_2)\geq -\inf\limits_{\varsigma\in B_2}I(\varsigma).
\de
\ed

 Note that the Laplace principle is equivalent to the LDP. Thus, in order to obtain that the family of $\{X^{\epsilon}, \e\in(0,1)\}$ satisfies the LDP in $(\mS,\sB(\mS))$ with the rate function $I$, we only need to prove that the family of $\{X^{\epsilon}, \e\in(0,1)\}$ satisfies the Laplace principle in $(\mS,\sB(\mS))$ with the rate function $I$. Let us state the conditions under which the Laplace principle holds. 
 
 Take any $T>0$. Set $\mathbb{H}:=L^{2}([0,T]; \mR^{d_1+d_2})$ and $\|h\|_{\mathbb{H}}:=(\int_{0}^{T}|h(t)|^{2}\dif t)^{\frac{1}{2}}$ for $h\in\mH$. Let $\mathcal{A}$ be the collection of predictable processes $u(\omega, \cdot)$ belonging to $\mathbb{H}$ a.s. $\omega$. For each $N\in\mN$ we define two following spaces:
\ce
\mathbf{D}_{2}^{N}:=\left\{h\in\mathbb{H}: \|h\|_{\mathbb{H}}^{2}\leq N \right\}, \quad \mathbf{A}_{2}^{N}:=\left\{u\in\mathcal{A}: u(\omega, \cdot)\in\mathbf{D}_{2}^{N}, \mP \ a.s.~\omega \right\}.
\de
We equip $\mathbf{D}_{2}^{N}$ with the weak convergence topology in $\mH$. So, $\mathbf{D}_{2}^{N}$ is metrizable as a compact Polish space. In the sequel, $\mathbf{D}_{2}^{N}$ will be always endowed with this topology.

\bco\label{cond}
Let $\Psi^{\epsilon} : C([0,T];\mathbb{R}^{d_1+d_2})\mapsto\mS$ be  a family of measurable mappings. There exists a
measurable mapping $\Psi^{0} : C([0,T];\mathbb{R}^{d_1+d_2})\mapsto\mS$ such that

$(i)$ for $N\in\mathbb{N}$ and $\{h_{\e}, \e>0\}\subset\mathbf{D}_2^{N}$, $h\in \mathbf{D}_2^{N}$, if $h_{\e}\rightarrow h$ as $\e\rightarrow 0$, then
\ce
\Psi^{0}\left(\int_{0}^{\cdot}h_{\e}(s)\dif s\right)\longrightarrow \Psi^{0}\left(\int_{0}^{\cdot}h(s)\dif s\right).
\de

$(ii)$ for $N\in\mathbb{N}$ and $\{u_{\epsilon},\epsilon>0\}\subset \mathbf{A}_{2}^{N}$, $u\in\mathbf{A}_{2}^{N}$, if $u_{\epsilon}$ converges in distribution to $u$ as $\e\rightarrow 0$, then
$$
\Psi^{\e}\left(\sqrt{\e}W_\cdot+\int_{0}^{\cdot}u_{\e}(s)\dif s\right) \overset{d}{\longrightarrow} \Psi^{0}\left(\int_{0}^{\cdot}u(s)\dif s\right),
$$
where $W_\cdot$ is a $d_1+d_2$-dimensional Brownian motion.
\eco

 Given $\varsigma\in\mS$, let ${\bf D}_{\varsigma}=\{h\in\mH: \varsigma=\Psi^{0}(\int_{0}^{\cdot}h(s)\dif s)\}$. Let $I:\mathbb{S}\mapsto [0,\infty]$ be defined by
$$
I(\varsigma)=\frac{1}{2}\inf\limits_{h\in{\bf D}_{\varsigma}}\|h\|_{\mathbb{H}}^{2}.
$$
The following result is due to \cite[Theorem 4.2]{BDM2}.

\bt\label{ldpbase}
Set $X^{\epsilon}:=\Psi^{\epsilon}(\sqrt{\e}W)$. Assume that Condition \ref{cond} holds. Then $\{X^{\epsilon}\}$ satisfies the Laplace principle with the rate function $I$ given above.
\et

\section{Formulation for main results}\label{mare}

In this section, we provide main results of this paper. 

\subsection{Averaging principles}\label{fun}

In this subsection, we state averaging principle results.

Recall the slow-fast system (\ref{Eq1i}), i.e. for any $T>0$
\be\left\{\begin{array}{l}
\dif X_{t}^{\e,\g}\in -A_1(X_{t}^{\e,\g})\dif t+b_{1}(t,X_{t}^{\e,\g},Y_{t}^{\e,\g})\dif t+\sqrt{\e}\s_{1}(t,X_{t}^{\e,\g},Y_{t}^{\e,\g})\dif W^1_{t},\\
X_{0}^{\e,\g}=x_0\in\overline{\cD(A_1)},\quad  0\leq t\leq T,\\
\dif Y_{t}^{\e,\g}\in -A_2(Y_{t}^{\e,\g})\dif t+\frac{1}{\g}b_{2}(t,X_{t}^{\e,\g},Y_{t}^{\e,\g})\dif t+\frac{1}{\sqrt{\g}}\s_{2}(t,X_{t}^{\e,\g},Y_{t}^{\e,\g})\dif W^2_{t},\\
Y_{0}^{\e,\g}=y_0\in\overline{\cD(A_2)},\quad  0\leq t\leq T.
\end{array}
\right.
\label{Eq1}
\ee

Assume that
\begin{enumerate}[$(\mathbf{H}_{A_1})$]
\item $0\in{\rm Int}(\cD(A_1))$.
\end{enumerate}
\begin{enumerate}[$(\mathbf{H}^1_{b_{1}, \s_{1}})$]
\item
There exists a constant $L_{b_{1}, \s_{1}}>0$ such that for $t_i\in[0,T], x_{i}\in\mR^n$, $y_{i}\in\mR^m$, $i=1, 2$,
\ce
&&|b_{1}(t_1, x_{1},y_{1})-b_{1}(t_2, x_{2},y_{2})|^{2}+\|\s_{1}(t_1,x_{1},y_{1})-\s_{1}(t_2,x_{2},y_{2})\|^{2}\\
&\leq& L_{b_{1},\s_{1}}\(|t_1-t_2|^2+|x_{1}-x_{2}|^{2}+|y_{1}-y_{2}|^{2}\).
\de
\end{enumerate}
\begin{enumerate}[$(\mathbf{H}_{A_2})$]
\item $0\in{\rm Int}(\cD(A_2))$.
\end{enumerate}
\begin{enumerate}[$(\mathbf{H}^1_{b_{2}, \s_{2}})$]
\item
There exists a constant $L_{b_{2}, \s_{2}}>0$ such that for $t_i\in[0,T], x_{i}\in\mR^n$, $y_{i}\in\mR^m$, $i=1, 2$,
\ce
&&|b_{2}(t_1,x_{1},y_{1})-b_{2}(t_2,x_{2},y_{2})|^{2}+\|\s_{2}(t_1,x_{1},y_{1})-\s_{2}(t_2,x_{2},y_{2})\|^{2}\\
&\leq& L_{b_{2}, \s_{2}}\(|t_1-t_2|^2+|x_{1}-x_{2}|^{2}+|y_{1}-y_{2}|^{2}\).
\de
\end{enumerate}

\br
$(i)$ $(\mathbf{H}^1_{b_{1}, \s_{1}})$ yields that there exists a constant $\bar{L}_{b_{1}, \s_{1}}>0$ such that for $t\in[0,T]$, $x\in\mR^n$, $y\in\mR^m$,
\be
|b_{1}(t,x,y)|^{2}+\|\s_{1}(t,x,y)\|^{2}\leq \bar{L}_{b_{1}, \s_{1}}(1+|x|^{2}+|y|^{2}).
\label{b1line}
\ee

$(ii)$ $(\mathbf{H}^1_{b_{2}, \s_{2}})$ implies that there exists a constant $\bar{L}_{b_{2}, \s_{2}}>0$ such that for $t\in[0,T]$, $x\in\mR^n$, $y\in\mR^m$, 
\be
|b_{2}(t,x,y)|^{2}+\|\s_{2}(t,x,y)\|^{2}
\leq \bar{L}_{b_{2}, \s_{2}}(1+|x|^{2}+|y|^{2}).
\label{b2nu}
\ee
\er

 Under $(\mathbf{H}_{A_{1}})$, $(\mathbf{H}_{A_{2}})$, $(\mathbf{H}^1_{b_{1}, \s_{1}})$ and $(\mathbf{H}^1_{b_{2}, \s_{2}})$, by the similar deduction to  that of \cite[Theorem 6.1]{cw},  
 we obtain the following result.

\bt\label{well}
Assume that $(\mathbf{H}_{A_{1}})$, $(\mathbf{H}_{A_{2}})$, $(\mathbf{H}^1_{b_{1}, \s_{1}})$, $(\mathbf{H}^1_{b_{2}, \s_{2}})$ hold. Then the system (\ref{Eq1}) has a unique strong solution $(X_{\cdot}^{\e,\g},K_{\cdot}^{1,\e,\g},Y_{\cdot}^{\e,\g},K_{\cdot}^{2,\e,\g})$.
\et

In order to obtain an averaging principle, we assume more:
\begin{enumerate}[$(\mathbf{H}^2_{b_{2}, \s_{2}})$]
\item
There exists a constant $\b>0$ satisfying $\b>2L_{b_2,\sigma_2}$ such that for $t\in[0,T], x\in\mR^n$, $y_{i}\in\mR^m$, $i=1, 2$,
\ce
2\<y_{1}-y_{2},b_{2}(t,x,y_{1})-b_{2}(t,x,y_{2})\>
+\|\s_{2}(t,x,y_{1})-\s_{2}(t,x,y_{2})\|^{2}\leq -\b|y_{1}-y_{2}|^{2}.
\de
\end{enumerate}

\br
By $(\mathbf{H}^1_{b_{2}, \s_{2}})$ and $(\mathbf{H}^2_{b_{2}, \s_{2}})$, it holds that for $t\in[0,T], x\in\mR^n$, $y\in\mR^m$,
\be
2\<y,b_{2}(t,x,y)\>+\|\s_{2}(t,x,y)\|^{2}\leq -\a|y|^{2}+C_{L_{b_{2}, \s_{2}}}(1+|x|^{2}),
\label{bemu}
\ee
where $\a:=\b-2L_{b_2,\sigma_2}>0$ and $C_{L_{b_{2}, \s_{2}}}>0$ are two constants.
\er

Take any $t\in[0,T], x\in \mR^{n}$ and fix it. Consider the following multivalued SDE:
\be\left\{\begin{array}{l}
\dif Y_{s}^{t,x}\in -A_2(Y_{s}^{t,x})\dif s+b_{2}(t,x,Y_{s}^{t,x})\dif s+\s_{2}(t,x,Y_{s}^{t,x})\dif W^2_{s},\\
Y_{0}^{t,x}=y_0\in\overline{\cD(A_2)}, \quad 0 \leq s \leq T.
\end{array}
\right.
\label{Eq2}
\ee
Under $(\mathbf{H}_{A_{2}})$ and $(\mathbf{H}^1_{b_{2}, \s_{2}})$, we know that the above equation has a unique strong solution $(Y_{\cdot}^{t,x,y_0},K_{\cdot}^{2,t,x,y_0})$ (\cite{rwzx}). Moreover, under $(\mathbf{H}^2_{b_{2}, \s_{2}})$, by the same deduction as that of \cite[Theorem 3.2]{q1}, one could obtain that there exists a unique invariant probability $\nu^{t,x}$ for Eq.(\ref{Eq2}). So, we construct an averaging equation as follows:
\be\left\{\begin{array}{l}
\dif \bar{X}_{t}\in -A_1(\bar{X}_{t})\dif t+\bar{b}_{1}(t,\bar{X}_{t})\dif t,\\
\bar{X}_{0}=x_0\in\overline{\cD(A_1)},\quad  0\leq t\leq T,
\end{array}
\right.
\label{Eq3}
\ee
where $\bar{b}_{1}(t,x):=\int_{\mR^{m}}b_{1}(t,x,y)\nu^{t,x}(\dif y)$. Note that Eq.(\ref{Eq3}) is a multivalued ordinary differential equation. By Lemma \ref{averc}, we know that Eq.(\ref{Eq3}) has a unique solution $(\bar{X}_{\cdot},\bar{K}_{\cdot})$.

Now, it is the position to state the first main result in this paper.

\bt \label{xbarxp}
Suppose that $(\mathbf{H}_{A_{1}})$, $(\mathbf{H}_{A_{2}})$, $(\mathbf{H}^1_{b_{1}, \s_{1}})$, $(\mathbf{H}^{1}_{b_{2}, \s_{2}})$ and $(\mathbf{H}^{2}_{b_{2}, \s_{2}})$ hold. If
\ce
\lim_{\e\rightarrow 0}\frac{\g}{\e}=\left\{\begin{array}{l} 0,\\
\varrho\in (0,\infty), 
\end{array}
\right.
\de
it holds that
\ce
\lim_{\e\rightarrow 0}\mE\(\sup_{0\leq t\leq T}|X_{t}^{\e,\g}-\bar{X}_{t}|^{2}\)=0,
\de
\et

For the sake of presentation, the proof of Theorem \ref{xbarxp} is deferred to Section \ref{proofirs}.

Next, we take $A_1=\p I_{\cO}$, where $\cO$ is a closed and convex domain in $\mR^n$ with ${\rm Int}(\cO)\neq\emptyset$ and obtain the following result.

\bc\label{convrate}
Suppose that $(\mathbf{H}_{A_{1}})$, $(\mathbf{H}_{A_{2}})$, $(\mathbf{H}^1_{b_{1}, \s_{1}})$, $(\mathbf{H}^{1}_{b_{2}, \s_{2}})$ and $(\mathbf{H}^{2}_{b_{2}, \s_{2}})$ hold. If
\ce
\lim_{\e\rightarrow 0}\frac{\g}{\e}=\left\{\begin{array}{l} 0,\\
\varrho\in (0,\infty), 
\end{array}
\right.
\de
it holds that for $0<\iota<1$
\ce
\mE\(\sup_{0\leq t\leq T}|X_{t}^{\e,\g}-\bar{X}_{t}|^{2}\)\leq C(\g^{\iota/2}+\g^\iota+\g^{2\iota}+\g^{\frac{1}{2}(1-\iota)}+\e).
\de
\ec

The proof of Corollary \ref{convrate} is also placed in Section \ref{proofirs}.

\br
If $\e=0, \iota=\frac{1}{2}$, by the above corollary we obtain that the convergence rate is $\frac{1}{8}$. Note that in \cite[Theorem 1]{wrd}, Wang, Robert and Duan proved that the convergence rate is $\frac{1}{2}$ for slow-fast stochastic reaction-diffusion equations. Therefore, we conclude that the maximal monotone operator $A_1$ reduces the convergence order.
\er

\subsection{Large deviation principles}

In this subsection, we require that $b_1, \s_1, b_2, \s_2$ is independent of $t$ and $\s_1$ is independent of $y$ and state the LDP results. 

First of all, we remember the system (\ref{Eq1ldp}), i.e.
\be\left\{\begin{array}{l}
\dif X_{t}^{\e,\g}\in -A_1(X_{t}^{\e,\g})\dif t+b_{1}(X_{t}^{\e,\g},Y_{t}^{\e,\g})\dif t+\sqrt{\e}\s_{1}(X_{t}^{\e,\g})\dif W^1_{t},\\
X_{0}^{\e,\g}=x_0\in\overline{\cD(A_1)},\quad  0\leq t\leq T,\\
\dif Y_{t}^{\e,\g}\in -A_2(Y_{t}^{\e,\g})\dif t+\frac{1}{\g}b_{2}(X_{t}^{\e,\g},Y_{t}^{\e,\g})\dif t+\frac{1}{\sqrt{\g}}\s_{2}(X_{t}^{\e,\g},Y_{t}^{\e,\g})\dif W^2_{t},\\
Y_{0}^{\e,\g}=y_0\in\overline{\cD(A_2)},\quad  0\leq t\leq T.
\end{array}
\right.
\label{Eq1ldp2}
\ee
Then under $(\mathbf{H}_{A_{1}})$, $(\mathbf{H}_{A_{2}})$, $(\mathbf{H}^1_{b_{1}, \s_{1}})$ and $(\mathbf{H}^{1}_{b_{2}, \s_{2}})$, $(X_{\cdot}^{\e,\g},K_{\cdot}^{1,\e,\g},Y_{\cdot}^{\e,\g},K_{\cdot}^{2,\e,\g})$ still stands for the unique strong solution of the above system. 

So as to establish the LDP result for the system (\ref{Eq1ldp2}), we also assume:
\begin{enumerate}[$(\mathbf{H}^3_{\s_{2}})$]
\item $\s_2(x,y)$ is bounded.
\end{enumerate}

Now, it is the position to describe the LDP result for the system (\ref{Eq1ldp2}), which is the second main result in this paper.

\bt\label{ldpmmsde}
Assume that $(\mathbf{H}_{A_{1}})$, $(\mathbf{H}_{A_{2}})$, $(\mathbf{H}^1_{b_{1}, \s_{1}})$, $(\mathbf{H}^{1}_{b_{2}, \s_{2}})$, $(\mathbf{H}^{2}_{b_{2}, \s_{2}})$ and $(\mathbf{H}^3_{\s_{2}})$ hold. If 
\ce
\lim_{\e\rightarrow 0}\frac{\g}{\e}=0,
\de
the family $\{X^{\epsilon,\g},\epsilon\in(0,1)\}$ satisfies the LDP in $\mS:=C([0,T],\overline{\mathcal{D}(A_1)})$ with the rate function given by
$$
I(\varsigma)=\frac{1}{2} \inf\limits_{h\in {\bf D}_{\varsigma}: \varsigma=\bar{X}^{h}}\|h\|_{\mH}^2,
$$
where $(\bar{X}^{h}, \bar{K}^{h})$ solves the following equation
\be\left\{\begin{array}{l}
\dif \bar{X}^h_{t}\in -A_1(\bar{X}^h_{t})\dif t+\bar{b}_{1}(\bar{X}^h_{t})\dif t+\s_{1}(\bar{X}^h_{t})\pi_1h(t)\dif t,\\
\bar{X}^h_{0}=x_0\in\overline{\cD(A_1)},
\end{array}
\right.
\label{ldppsioequ}
\ee
and $\pi_1: \mR^{d_1+d_2}\mapsto \mR^{d_1}$ is a projection operator.
\et

Next, we take $A_1=\p I_{D}$, where for a function $\kappa\in C_b^2(\mR^n)$, $D=\{x\in\mR^n, \kappa(x)\geq 0\}$ is a closed convex domain with $\p D=\{x\in\mR^n, \kappa(x)=0\}$ and $|\triangledown \kappa(x)|=1$ for $x\in\p D$. Note that at any boundary point $x\in\p D$, $\triangledown \kappa(x)$ is a unit normal vector to the boundary, pointing towards the interior of $D$. Let $\phi\in C([0,T], D), \psi\in C([0,T], \mR^n), V\in\sV_0$ such that
\ce
\phi(t)=\psi(t)+V(t), \quad V_t=\int_0^t\triangledown \kappa(\phi_r)\dif |V|_0^r, \quad |V|_0^t=\int_0^tI_{\{\phi_r\in\p D\}}\dif |V|_0^r.
\de
For $\phi, \psi$ as above, we define a mapping $\Gamma: C([0,T], \mR^n)\rightarrow C([0,T], D)$ by
$$
\phi=\Gamma(\psi).
$$

Now we assume:
\begin{enumerate}[$(\mathbf{H}^2_{\s_{1}})$]
\item $(\s_1\s_1^*)(x)$ is uniformly nondegenerate.
\end{enumerate}

In the following corollary, we characterize the rate function.
\bc\label{charate}
Assume that $(\mathbf{H}_{A_{1}})$, $(\mathbf{H}_{A_{2}})$, $(\mathbf{H}^1_{b_{1}, \s_{1}})$, $(\mathbf{H}^{1}_{b_{2}, \s_{2}})$, $(\mathbf{H}^{2}_{b_{2}, \s_{2}})$, $(\mathbf{H}^3_{\s_{2}})$ and $(\mathbf{H}^2_{\s_{1}})$ hold. If 
\ce
\lim_{\e\rightarrow 0}\frac{\g}{\e}=0,
\de
the family $\{X^{\epsilon,\g},\epsilon\in(0,1)\}$ satisfies the LDP in $\mS:=C([0,T],\overline{\mathcal{D}(\p I_{D})})$ with the rate function given by
\ce
I(\phi)=\left\{\begin{array}{l}\frac{1}{2}\inf\limits_{\psi:\phi=\Gamma(\psi)}\int_0^T\(\dot{\psi}(t)-\bar{b}_1(\phi(t))\)^*(\s_1\s_1^*)^{-1}(\phi(t))\(\dot{\psi}(t)-\bar{b}_1(\phi(t))\)\dif t,\\
 ~\mbox{if}~ \psi\in AC([0,T],\mR^n) ~\mbox{and}~ \psi(0)=x_0,\\
\infty, \quad otherwise,
\end{array}
\right.
\de
where $AC([0,T],\mR^n)$ is the set of absolutely continuous functions from $[0,T]$ to $\mR^n$.
\ec

\br
We mention that if $D=\mR^n$, that is, the mapping $\Gamma$ is identical, the above rate function becomes
\ce
I(\phi)=\left\{\begin{array}{l}\frac{1}{2}\int_0^T\(\dot{\phi}(t)-\bar{b}_1(\phi(t))\)^*(\s_1\s_1^*)^{-1}(\phi(t))\(\dot{\phi}(t)-\bar{b}_1(\phi(t))\)\dif t,\\
 ~\mbox{if}~ \phi\in AC([0,T],\mR^n) ~\mbox{and}~ \phi(0)=x_0,\\
\infty, \quad otherwise,
\end{array}
\right.
\de
which is the same to that in \cite[Theorem 3.4]{ks1}.
\er

The proofs of Theorem \ref{ldpmmsde} and Corollary \ref{charate} are placed in Section \ref{prooseco}.

\section{Preparations}\label{prep}

In this section, we make preparations so as to prove the main results.

\subsection{Some estimates for the frozen equation (\ref{Eq2})}\label{froequ}

\bl
Under $(\mathbf{H}_{A_{2}})$, $(\mathbf{H}^1_{b_{2}, \s_{2}})$, $(\mathbf{H}^2_{b_{2}, \s_{2}})$, it holds that for any $s,t, t_1,t_2\in[0,T], x,x_1,\\x_2\in\mR^n, y_1,y_2\in\mR^m$
\be
&&\mE|Y_{s}^{t,x,y_0}|^{2}\leq|y_0|^{2}e^{-\frac{\a}{2} s}+C(1+|x|^{2}), \label{memu2}\\
&&\mE|Y_{s}^{t_1,x_1,y_1}-Y_{s}^{t_2,x_2,y_2}|^{2}\leq |y_{1}-y_{2}|^{2}e^{-\a s}+C(|t_1-t_2|^2+|x_1-x_2|^2).\label{memu20}
\ee
\el
\begin{proof}
Applying the It\^{o} formula to $|Y_{s}^{t,x,y_0}|^{2}e^{\l_1 s}$ for any $\l_1>0$ and taking the expectation, we get that
\ce
\mE|Y_{s}^{t,x, y_0}|^{2}e^{\l_1 s}
&=&|y_0|^{2}+\l_1\mE\int_{0}^{s}e^{\l_1 u}|Y_{u}^{t,x, y_0}|^{2}\dif u-2\mE\int_{0}^{s}e^{\l_1 u}\<Y_{u}^{t,x, y_0}, \dif K_u^{2,t,x,y_0}\>\\
&&+2\mE\int_{0}^{s}e^{\l_1 u}\<Y_{u}^{x, y_0}, b_{2}(t,x,Y_{u}^{t,x, y_0})\>\dif u+\mE\int_{0}^{s}e^{\l_1 u}\|\s_{2}(t,x, Y_{u}^{t,x, y_0})\|^{2}\dif u\\
&\leq&|y_0|^{2}+\l_1\mE\int_{0}^{s}e^{\l_1 u}|Y_{u}^{t,x, y_0}|^{2}\dif u+2\mE\int_{0}^{s}e^{\l_1 u}|v||Y_{u}^{t,x, y_0}|\dif u\\
&&+\mE\int_{0}^{s}e^{\l_1 u}\(-\a|Y_{u}^{t,x,y_0}|^{2}+C(1+|x|^{2})\)\dif u\\
&\leq&|y_0|^{2}+\l_1\mE\int_{0}^{s}e^{\l_1 u}|Y_{u}^{t,x, y_0}|^{2}\dif u+\frac{\a}{2}\mE\int_{0}^{s}e^{\l_1 u}|Y_{u}^{t,x, y_0}|^{2}\dif u+C\int_{0}^{s}e^{\l_1 u}|v|^2\dif u\\
&&+\mE\int_{0}^{s}e^{\l_1 u}\(-\a|Y_{u}^{t,x,y_0}|^{2}+C(1+|x|^{2})\)\dif u\\
&\leq&|y_0|^{2}+(\l_1+\frac{\a}{2}-\a)\mE\int_{0}^{s}e^{\l_1 u}|Y_{u}^{t,x, y_0}|^{2}\dif u+C(1+|x|^{2})\int_{0}^{s}e^{\l_1 u}\dif u,
\de
where $v\in A_2(0)$ and Lemma \ref{equi} is used. By the above deduction and letting $\l_1=\frac{\a}{2}$, it holds that
\ce
\mE|Y_{s}^{t,x,y_0}|^{2}\leq |y_0|^{2}e^{-\frac{\a}{2} s}+C(1+|x|^{2}).
\de

Next, we deal with the second estimate. Let $(Y_{\cdot}^{t_1,x_1,y_{1}},K_{\cdot}^{2,t_1,x_1,y_{1}})$ and $(Y_{\cdot}^{t_2,x_2,y_{2}},K_{\cdot}^{2,t_2,x_2,y_{2}})$ be the solutions to Eq.(\ref{Eq2}) with $t=t_1, x=x_1, Y^{t_1,x_1,y_1}_0=y_{1}$ and $t=t_2, x=x_2, Y^{t_2,x_2,y_2}_0=y_{2}$, respectively. Applying the It\^{o} formula to $|Y_{s}^{t_1,x_1,y_1}-Y_{s}^{t_2,x_2,y_2}|^{2}e^{\l_2 s}$ for any $\l_2>0$ and taking expectation on two sides, by Lemma \ref{equi} we obtain that
\ce
&&\mE|Y_{s}^{t_1,x_1,y_1}-Y_{s}^{t_2,x_2,y_2}|^{2}e^{\l_2 s}\\
&=&|y_{1}-y_{2}|^{2}+\l_2\mE\int_{0}^{s}e^{\l_2 u}|Y_{u}^{t_1,x_1,y_1}-Y_{u}^{t_2,x_2,y_2}|^{2}\dif u\\
&&-2\mE\int_{0}^{s}e^{\l_2 u}\< Y_{u}^{t_1,x_1,y_1}-Y_{u}^{t_2,x_2,y_2},\dif (K_{u}^{2,t_1,x_1,y_1}-K_{u}^{2,t_2,x_2,y_2})\>\\
&&+2\mE\int_{0}^{s}e^{\l_2 u}\<  Y_{u}^{t_1,x_1,y_1}-Y_{u}^{t_2,x_2,y_2}, b_{2}(t_1,x_1,Y_{u}^{t_1,x_1,y_1})-b_{2}(t_2,x_2,Y_{u}^{t_2,x_2,y_2})\>
\dif u\\
&&+\mE\int_{0}^{s}e^{\l_2 u}\|\s_{2}(t_1,x_1,Y_{u}^{t_1,x_1,y_1})-\s_{2}(t_2,x_2,Y_{u}^{t_2,x_2,y_2})\|^{2}\dif u\\
&\leq&|y_{1}-y_{2}|^{2}+\l_2\mE\int_{0}^{s}e^{\l_2 u}|Y_{u}^{t_1,x_1,y_1}-Y_{u}^{t_2,x_2,y_2}|^{2}\dif u\\
&&+2\mE\int_{0}^{s}e^{\l_2 u}\< Y_{u}^{t_1,x_1,y_1}-Y_{u}^{t_2,x_2,y_2}, b_{2}(t_1,x_1,Y_{u}^{t_1,x_1,y_1})-b_{2}(t_1,x_1,Y_{u}^{t_2,x_2,y_2})\>\dif u\\
&&+\mE\int_{0}^{s}e^{\l_2 u}\|\s_{2}(t_1,x_1,Y_{u}^{t_1,x_1,y_1})-\s_{2}(t_1,x_1,Y_{u}^{t_2,x_2,y_2})\|^{2}\dif u\\
&&+2\mE\int_{0}^{s}e^{\l_2 u}\< Y_{u}^{t_1,x_1,y_1}-Y_{u}^{t_2,x_2,y_2}, b_{2}(t_1,x_1,Y_{u}^{t_2,x_2,y_2})-b_{2}(t_2,x_2,Y_{u}^{t_2,x_2,y_2})\>\dif u\\
&&+\mE\int_{0}^{s}e^{\l_2 u}\|\s_{2}(t_1,x_1,Y_{u}^{t_1,x_1,y_1})-\s_{2}(t_1,x_1,Y_{u}^{t_2,x_2,y_2})\|^{2}\dif u\\
&&+2\mE\int_{0}^{s}e^{\l_2 u}\|\s_{2}(t_1,x_1,Y_{u}^{t_2,x_2,y_2})-\s_{2}(t_2,x_2,Y_{u}^{t_2,x_2,y_2})\|^{2}\dif u.
\de
Then $(\mathbf{H}^2_{b_{2}, \s_{2}})$ yields that
\ce
&&2\mE\int_{0}^{s}e^{\l_2 u}\< Y_{u}^{t_1,x_1,y_1}-Y_{u}^{t_2,x_2,y_2}, b_{2}(t_1,x_1,Y_{u}^{t_1,x_1,y_1})-b_{2}(t_1,x_1,Y_{u}^{t_2,x_2,y_2})\>\dif u\\
&&+\mE\int_{0}^{s}e^{\l_2 u}\|\s_{2}(t_1,x_1,Y_{u}^{t_1,x_1,y_1})-\s_{2}(t_1,x_1,Y_{u}^{t_2,x_2,y_2})\|^{2}\dif u\\
&\leq&-\b\mE\int_{0}^{s}e^{\l_2 u}|Y_{u}^{t_1,x_1,y_1}-Y_{u}^{t_2,x_2,y_2}|^{2}\dif u.
\de
And by $(\mathbf{H}^1_{b_{2}, \s_{2}})$ and the Young inequality, we infer that
\ce
&&2\mE\int_{0}^{s}e^{\l_2 u}\< Y_{u}^{t_1,x_1,y_1}-Y_{u}^{t_2,x_2,y_2}, b_{2}(t_1,x_1,Y_{u}^{t_2,x_2,y_2})-b_{2}(t_2,x_2,Y_{u}^{t_2,x_2,y_2})\>\dif u\\
&\leq&L_{b_2,\s_2}\mE\int_{0}^{s}e^{\l_2 u}|Y_{u}^{t_1,x_1,y_1}-Y_{u}^{t_2,x_2,y_2}|^{2}\dif u\\
&&+\frac{1}{L_{b_2,\s_2}}\mE\int_{0}^{s}e^{\l_2 u}|b_{2}(t_1,x_1,Y_{u}^{t_2,x_2,y_2})-b_{2}(t_2,x_2,Y_{u}^{t_2,x_2,y_2})|^2\dif u\\
&\leq&L_{b_2,\s_2}\mE\int_{0}^{s}e^{\l_2 u}|Y_{u}^{t_1,x_1,y_1}-Y_{u}^{t_2,x_2,y_2}|^{2}\dif u+\mE\int_{0}^{s}e^{\l_2 u}(|t_1-t_2|^2+|x_1-x_2|^2)\dif u,
\de
and
\ce
&&\mE\int_{0}^{s}e^{\l_2 u}\|\s_{2}(t_1,x_1,Y_{u}^{t_1,x_1,y_1})-\s_{2}(t_1,x_1,Y_{u}^{t_2,x_2,y_2})\|^{2}\dif u\\
&&+2\mE\int_{0}^{s}e^{\l_2 u}\|\s_{2}(t_1,x_1,Y_{u}^{t_2,x_2,y_2})-\s_{2}(t_2,x_2,Y_{u}^{t_2,x_2,y_2})\|^{2}\dif u\\
&\leq&L_{b_2,\s_2}\mE\int_{0}^{s}e^{\l_2 u}|Y_{u}^{t_1,x_1,y_1}-Y_{u}^{t_2,x_2,y_2}|^{2}\dif u+2L_{b_2,\s_2}\mE\int_{0}^{s}e^{\l_2 u}(|t_1-t_2|^2+|x_1-x_2|^2)\dif u.
\de
Collecting the above deduction and taking $\l_2=\a$, we conclude that
\ce
\mE|Y_{s}^{t_1,x_1,y_1}-Y_{s}^{t_2,x_2,y_2}|^{2}e^{\l_2 s}&\leq&|y_{1}-y_{2}|^{2}+(\a-\b+2L_{b_2,\s_2})\mE\int_{0}^{s}e^{\l_2 u}|Y_{u}^{t_1,x_1,y_1}-Y_{u}^{t_2,x_2,y_2}|^{2}\dif u\\
&&+C(|t_1-t_2|^2+|x_1-x_2|^2)\frac{e^{\l_2 s}-1}{\l_2},
\de
which together with $\a=\b-2L_{b_2,\s_2}$ implies that
$$
\mE|Y_{s}^{t_1,x_1,y_1}-Y_{s}^{t_2,x_2,y_2}|^{2}\leq |y_{1}-y_{2}|^{2}e^{-\a s}+C(|t_1-t_2|^2+|x_1-x_2|^2).
$$
The proof is complete.
\end{proof}

By the above lemma and the definition of $\nu^{t,x}$, it holds that
\ce
\int_{\overline{\cD(A_2)}}|y|^{2}\nu^{t,x}(\dif y)
&=&\int_{\overline{\cD(A_2)}}\mE|Y_{s}^{t,x,y}|^{2}\nu^{t,x}(\dif y)\leq \int_{\overline{\cD(A_2)}}\(|y|^{2}e^{-\frac{\a}{2} s}+C(1+|x|^{2})\)\nu^{t,x}(\dif y)\no\\
&=& e^{-\frac{\a}{2} s}\int_{\overline{\cD(A_2)}}|y|^{2}\nu^{t,x}(\dif y)+C(1+|x|^{2}),
\de
and furthermore
\be
\int_{\overline{\cD(A_2)}}|y|^{2}\nu^{t,x}(\dif y)\leq C(1+|x|^{2}).
\label{inu2}
\ee

\bl\label{emb1}
 Suppose that $(\mathbf{H}^1_{b_{1}, \s_{1}})$, $(\mathbf{H}_{A_{2}})$, $(\mathbf{H}^{1}_{b_{2}, \s_{2}})$, $(\mathbf{H}^2_{b_{2}, \s_{2}})$ hold. Then, for any $s\in[0, T]$ there exists a constant $C>0$ such that
\be
|\mE b_{1}(t,x,Y_{s}^{t,x,y_0})-\bar{b}_{1}(t,x)|^{2}\leq Ce^{-\a s}(1+|x|^{2}+|y_0|^{2}).
\label{meu2}
\ee
\el
\begin{proof}
Based on simple calculations, one can obtain that
\ce
&&|\mE b_{1}(t,x,Y_{s}^{t,x,y_0})-\bar{b}_{1}(t,x)|^{2}
=\Big|\mE b_{1}(t,x,Y_{s}^{t,x,y_0})-\int_{\overline{\cD(A_2)}}b_{1}(t,x,y)\nu^{t,x}(\dif y)\Big|^{2}\\
&=&\Big|\mE b_{1}(t,x,Y_{s}^{t,x,y_0})-\int_{\overline{\cD(A_2)}}\mE b_{1}(t,x,Y_{s}^{t,x,y})\nu^{t,x}(\dif y)\Big|^{2}\\
&\leq&\int_{\overline{\cD(A_2)}}\mE |b_{1}(t,x,Y_{s}^{t,x,y_0})- b_{1}(t,x,Y_{s}^{t,x,y})|^{2}\nu^{t,x}(\dif y)\\
&\leq& L_{b_{1},\s_{1}}\int_{\overline{\cD(A_2)}}\mE |Y_{s}^{t,x,y_0}-Y_{s}^{t,x,y}|^{2}\nu^{t,x}(\dif y)\\
&\leq& L_{b_{1},\s_{1}}\int_{\overline{\cD(A_2)}}|y_0-y|^{2}e^{-\a s}\nu^{t,x}(\dif y)\\
&\leq&2L_{b_{1},\s_{1}}e^{-\a s}\left(|y_0|^{2}+\int_{\overline{\cD(A_2)}}|y|^{2}\nu^{t,x}(\dif y)\right),
\de
where the second and third inequalities are based on $(\mathbf{H}^1_{b_{1}, \s_{1}})$ and (\ref{memu20}), respectively. Finally, (\ref{inu2}) implies the required estimate.
\end{proof}

Combining the above estimates, we obtain the following lemma which assures the well-posedness of the averaging equation.

\bl
Under $(\mathbf{H}_{A_{2}})$, $(\mathbf{H}^1_{b_{1}, \s_{1}})$, $(\mathbf{H}^{1}_{b_{2}, \s_{2}})$, $(\mathbf{H}^2_{b_{2}, \s_{2}})$, it holds that for any $x_1, x_2\in\mR^n$
\be
|\bar{b}_1(t_1,x_1)-\bar{b}_1(t_2,x_2)|^2\leq C(|t_1-t_2|^2+|x_1-x_2|^2).
\label{barb1lip}
\ee
\el
\begin{proof}
From (\ref{meu2}) and (\ref{memu20}), it follows that for any $x_1, x_2\in\mR^n$
\ce
|\bar{b}_1(t_1,x_1)-\bar{b}_1(t_2,x_2)|^2&\leq& 3|\bar{b}_1(t_1,x_1)-\mE b_{1}(t_1,x_1,Y_{s}^{t_1,x_1,y_0})|^2\\
&&+3|\mE b_{1}(t_1,x_1,Y_{s}^{t_1,x_1,y_0})-\mE b_{1}(t_2,x_2,Y_{s}^{t_2,x_2,y_0})|^2\\
&&+3|\mE b_{1}(t_2,x_2,Y_{s}^{t_2,x_2,y_0})-\bar{b}_1(t_2,x_2)|^2\\
&\leq& 3Ce^{-\a s}(2+|x_1|^{2}+|x_2|^{2}+|y_0|^{2})+3C(|t_1-t_2|^2+|x_1-x_2|^2).
\de
Letting $s\rightarrow \infty$, we obtain that
\ce
|\bar{b}_1(t_1,x_1)-\bar{b}_1(t_2,x_2)|\leq 3C(|t_1-t_2|^2+|x_1-x_2|^2),
\de
which completes the proof.
\end{proof}

\section{Proofs of Theorem \ref{xbarxp} and Corollary \ref{convrate}}\label{proofirs}

In this section, we prove Theorem \ref{xbarxp} and Corollary \ref{convrate}. 

\subsection{Proof of Theorem \ref{xbarxp}}

In this subsection, we show Theorem \ref{xbarxp}. The proof consists of two parts. In the first part (Subsubsection \ref{xtzthatxhatz}), we segment the time interval $[0, T]$ by the size $\d$, where $\d$ is a fixed positive number depending on $\g$, and introduce two auxiliary processes: 
\be\left\{\begin{array}{l}
\dif\hat{Y}_{t}^{\e,\g}\in -A_2(\hat{Y}_{t}^{\e,\g})\dif t+\frac{1}{\g}b_{2}(t(\d),X_{t(\d)}^{\e,\g},\hat{Y}_{t}^{\e,\g})\dif t
+\frac{1}{\sqrt{\g}}\s_{2}(t(\d),X_{t(\d)}^{\e,\g},\hat{Y}_{t}^{\e,\g})\dif W^2_{t}, \\
\hat{Y}_{0}^{\e,\g}=y_0\in\overline{\cD(A_2)},
\end{array}
\right.
\label{hatz}\\
\left\{\begin{array}{l}
\dif\hat{X}_t^{\e,\g}\in -A_1(\hat{X}_t^{\e,\g})\dif t+b_1(t(\d),X_{t(\d)}^{\e,\g},\hat{Y}_{t}^{\e,\g})\dif t+\sqrt{\e}\s_{1}(t,X_{t}^{\e,\g},Y_{t}^{\e,\g})\dif W^1_{t}, \\
\hat{X}_0^{\e,\g}=x_0\in\overline{\cD(A_1)},
\end{array}
\right.
\label{hatx}
\ee
where $t(\d)=[\frac{t}{\d}]\d$, and $[\frac{t}{\d}]$ denotes the integer part of $\frac{t}{\d}$. $(\hat{Y}^{\e,\g}, \hat{K}^{2,\e,\g}), (\hat{X}^{\e,\g},\hat{K}^{1,\e,\g})$ denote strong solutions of Eq.(\ref{hatz}) and (\ref{hatx}), respectively. Then we estimate $X^{\e,\g}, K^{1,\e,\g}, Y^{\e,\g}, \hat{X}^{\e,\g}, \\\hat{K}^{1,\e,\g}, \hat{Y}^{\e,\g}$. In the second part (Subsection \ref{aveequ}), we present some estimates for the averaging equation (\ref{Eq3}).

\subsubsection{Some necessary estimates}\label{xtzthatxhatz}

\bl \label{xtztc}
 Under the assumptions of Theorem \ref{xbarxp}, for $\e\in (0,1)$ there exists a constant $C>0$ independent of $\e, \g$ such that 
\be
&&\mE\left(\sup\limits_{t\in[0,T]}|X_{t}^{\e,\g}|^{2}\right)\leq C(1+|x_0|^{2}+|y_0|^{2}), \label{xeb}\\
&&\sup\limits_{t\in[0,T]}\mE|Y_{t}^{\e,\g}|^{2}\leq C(1+|x_0|^{2}+|y_0|^{2}), \label{zeb}\\
&&\mE|K^{1,\e,\g}|_0^T\leq C(1+|x_0|^{2}+|y_0|^{2}). \label{keb}
\ee
\el
\begin{proof}
First of all, we estimate $X_{\cdot}^{\e,\g}$. Note that $X_{\cdot}^{\e,\g}$ satisfies the following equation:
\ce
X_{t}^{\e,\g}=x_0-K_t^{1,\e,\g}+\int_0^t b_{1}(s,X_{s}^{\e,\g},Y_{s}^{\e,\g})\dif s+\sqrt{\e}\int_0^t\s_{1}(s,X_{s}^{\e,\g},Y_{s}^{\e,\g})\dif W^1_{s}.
\de
By applying the It\^o formula to $|X_{t}^{\e,\g}|^2$, it holds that 
\be
|X_{t}^{\e,\g}|^2&=&|x_0|^2-2\int_0^t\<X_{s}^{\e,\g}, \dif K_s^{1,\e,\g}\>+2\int_0^t\<X_{s}^{\e,\g}, b_{1}(s,X_{s}^{\e,\g},Y_{s}^{\e,\g})\>\dif s\no\\
&&+2\sqrt{\e}\int_0^t\<X_{s}^{\e,\g},\s_{1}(s,X_{s}^{\e,\g},Y_{s}^{\e,\g})\dif W^1_{s}\>+\e\int_0^t\|\s_{1}(s,X_{s}^{\e,\g},Y_{s}^{\e,\g})\|^2\dif s, 
\label{xegito}
\ee
and furthermore for any $v'\in A_1(0)$,
\ce 
|X_{t}^{\e,\g}|^2&\leq&|x_0|^2+2|v'|\int_0^t|X_{s}^{\e,\g}|\dif s+\int_0^t|X_{s}^{\e,\g}|^2\dif s+\int_0^t|b_{1}(s,X_{s}^{\e,\g},Y_{s}^{\e,\g})|^2\dif s\no\\
&&+2\left|\int_0^t\<X_{s}^{\e,\g},\s_{1}(s,X_{s}^{\e,\g},Y_{s}^{\e,\g})\dif W^1_{s}\>\right|+\int_0^t\|\s_{1}(s,X_{s}^{\e,\g},Y_{s}^{\e,\g})\|^2\dif s,\no
\de
where Lemma \ref{equi} is used in the above inequality. Then the BDG inequality and $(\ref{b1line})$ imply that
\ce
&&\mE\left(\sup\limits_{s\in[0,t]}|X_{s}^{\e,\g}|^{2}\right)\\
&\leq&(|x_0|^2+2|v'|T)+(2|v'|+1)\mE\int_0^t|X_{r}^{\e,\g}|^2\dif r+\bar{L}_{b_{1}, \s_{1}}\mE\int_0^t(1+|X_{r}^{\e,\g}|^2+|Y_{r}^{\e,\g}|^2)\dif r\no\\
&&+2\mE\sup\limits_{s\in[0,t]}\left|\int_0^s\<X_{r}^{\e,\g},\s_{1}(r,X_{r}^{\e,\g},Y_{r}^{\e,\g})\dif W^1_{r}\>\right|\no\\ 
&\leq&C(|x_0|^{2}+1)+C\mE\int_{0}^{t}|X_{r}^{\e,\g}|^2\dif r+C\mE\int_{0}^{t}|Y_{r}^{\e,\g}|^{2}\dif r\no\\
&&+C\mE\left(\int_0^t|X_{r}^{\e,\g}|^2\|\s_{1}(r,X_{r}^{\e,\g},Y_{r}^{\e,\g})\|^2\dif r\right)^{1/2}\no\\
&\leq&C(|x_0|^{2}+1)+C\mE\int_{0}^{t}|X_{r}^{\e,\g}|^2\dif r+C\mE\int_{0}^{t}|Y_{r}^{\e,\g}|^{2}\dif r\no\\
&&+\frac{1}{2}\mE\left(\sup\limits_{s\in[0,t]}|X_{s}^{\e,\g}|^{2}\right)+C\mE\int_0^t\|\s_{1}(r,X_{r}^{\e,\g},Y_{r}^{\e,\g})\|^2\dif r,
\de
and furthermore
\be
\mE\left(\sup\limits_{s\in[0,t]}|X_{s}^{\e,\g}|^{2}\right)\leq C(|x_0|^{2}+1)+C\int_{0}^{t}\mE|X_{r}^{\e,\g}|^{2}\dif r+C\int_{0}^{t}\mE|Y_{r}^{\e,\g}|^{2}\dif r.
\label{exqc}
\ee

For $Y_{t}^{\e,\g}$, fix $v\in A_2(0)$. Applying the It\^{o} formula to $|Y_{t}^{\e,\g}|^{2}e^{\l t}$ for any $\l>0$ and taking the expectation, one could obtain that
\ce
&&\mE|Y_{t}^{\e,\g}|^{2}e^{\l t}\\
&=& |y_0|^{2}+\l\mE\int_0^t|Y_{s}^{\e,\g}|^{2}e^{\l s}\dif s-2\mE\int_0^te^{\l s}\<Y_{s}^{\e,\g},\dif K_s^{2,\e,\g}\>\\
&&+\frac{2}{\g}\mE\int_{0}^{t}e^{\l s}\<Y_{s}^{\e,\g}, b_{2}(s,X_{s}^{\e,\g},Y_{s}^{\e,\g})\>\dif s+\frac{1}{\g}\mE\int_{0}^{t}e^{\l s}\|\s_{2}(s,X_{s}^{\e,\g},Y_{s}^{\e,\g})\|^2\dif s\\
&\leq& |y_0|^{2}+\l\mE\int_0^t|Y_{s}^{\e,\g}|^{2}e^{\l s}\dif s+2\mE\int_0^te^{\l s}|v||Y_{s}^{\e,\g}|\dif s\\
&&+\frac{1}{\g}\mE\int_{0}^{t}e^{\l s}\(-\a|Y_{s}^{\e,\g}|^{2}+C(1+|X_{s}^{\e,\g}|^{2})\)\dif s\\
&\leq& |y_0|^{2}+\l\mE\int_0^t|Y_{s}^{\e,\g}|^{2}e^{\l s}\dif s+\frac{\a}{2\g}\mE\int_0^t|Y_{s}^{\e,\g}|^{2}e^{\l s}\dif s+C\mE\int_0^te^{\l s}|v|^2\dif s\\
&&+\frac{1}{\g}\mE\int_{0}^{t}e^{\l s}\(-\a|Y_{s}^{\e,\g}|^{2}+C(1+|X_{s}^{\e,\g}|^{2})\)\dif s\\
&\leq& |y_0|^{2}+\(\l+\frac{\a}{2\g}-\frac{\a}{\g}\)\mE\int_0^t|Y_{s}^{\e,\g}|^{2}e^{\l s}\dif s+C\int_0^te^{\l s}|v|^2\dif s+\frac{C}{\g}\mE\int_{0}^{t}e^{\l s}(1+|X_{s}^{\e,\g}|^{2})\dif s,
\de
where the first inequality is based on Lemma \ref{equi} and (\ref{bemu}). From this and taking $\l=\frac{\a}{2\g}$, it follows that
\be
\mE|Y_{t}^{\e,\g}|^{2}\leq C(|y_0|^{2}+1)+C\mE\left(\sup\limits_{s\in[0,t]}|X_{s}^{\e,\g}|^{2}\right).
\label{zees}
\ee

Inserting (\ref{zees}) in (\ref{exqc}), by the Gronwall inequality one can get (\ref{xeb}) and (\ref{zeb}).

Finally, for $K_\cdot^{1,\e,\g}$, by (\ref{xegito}) and Lemma \ref{inteineq}, it holds that
\ce
|X_{T}^{\e,\g}|^2&=&|x_0|^2-2\int_0^T\<X_{s}^{\e,\g}, \dif K_s^{1,\e,\g}\>+2\int_0^T\<X_{s}^{\e,\g}, b_{1}(s,X_{s}^{\e,\g},Y_{s}^{\e,\g})\>\dif s\no\\
&&+2\sqrt{\e}\int_0^T\<X_{s}^{\e,\g},\s_{1}(s,X_{s}^{\e,\g},Y_{s}^{\e,\g})\dif W^1_{s}\>+\e\int_0^T\|\s_{1}(s,X_{s}^{\e,\g},Y_{s}^{\e,\g})\|^2\dif s\no\\
&\leq&|x_0|^2-2M_1\left| K^{1,\e,\g} \right|_{0}^{T}+2M _2\int_0^T{\left| X^{\e,\g}_s\right|}\dif s+2M_3T+\int_0^T|X_{s}^{\e,\g}|^2\dif s\\
&&+\int_0^T|b_{1}(s,X_{s}^{\e,\g},Y_{s}^{\e,\g})|^2\dif s+2\sqrt{\e}\int_0^T\<X_{s}^{\e,\g},\s_{1}(s,X_{s}^{\e,\g},Y_{s}^{\e,\g})\dif W^1_{s}\>\\
&&+\e\int_0^T\|\s_{1}(s,X_{s}^{\e,\g},Y_{s}^{\e,\g})\|^2\dif s,
\de
and 
\ce
2M_1\mE\left| K^{1,\e,\g} \right|_{0}^{T}&\leq&|x_0|^2+2(M_2+M_3)T+(2M _2+1)\int_0^T\mE{\left| X^{\e,\g}_s\right|^2}\dif s\\
&&+\mE\int_0^T|b_{1}(s,X_{s}^{\e,\g},Y_{s}^{\e,\g})|^2\dif s+\mE\int_0^T\|\s_{1}(s,X_{s}^{\e,\g},Y_{s}^{\e,\g})\|^2\dif s\\
&\leq& |x_0|^2+2(M_2+M_3)T+(2M _2+1)\int_0^T\mE{\left| X^{\e,\g}_s\right|^2}\dif s\\
&&+C\int_0^T(1+\mE{\left| X^{\e,\g}_s\right|^2}+\mE{\left| Y^{\e,\g}_s\right|^2})\dif s,
\de
which together with (\ref{xeb}) and (\ref{zeb}) yields (\ref{keb}). The proof is complete.
\end{proof}

By the same deduction to that in Lemma \ref{xtztc}, we obtain the following estimate.
\bl
Under the assumptions of Theorem \ref{xbarxp}, it holds that for $0<\e<1$
\be
&&\sup\limits_{t\in[0,T]}\mE|\hat{Y}_{t}^{\e,\g}|^{2}\leq C(1+|x_0|^{2}+|y_0|^{2}), \label{hatzb}\\
&&\mE\left(\sup\limits_{t\in[0,T]}|\hat{X}_{t}^{\e,\g}|^{2}\right)\leq C(1+|x_0|^{2}+|y_0|^{2}),\label{hatxb}\\
&&\mE|\hat{K}^{1,\e,\g}|_0^T\leq C(1+|x_0|^{2}+|y_0|^{2}),\label{hatkb}
\ee
where the constant $C>0$ is independent of $\e, \g$.
\el

Next, we estimate some differences. 

\bl\label{xehatxets}
Under the assumptions of Theorem \ref{xbarxp}, we have that
\be
&&\lim\limits_{l\rightarrow 0}\sup\limits_{s\in[0,T]}\mE\sup _{s \leqslant t \leqslant s+l}|X_{t}^{\e,\g}-X_{s}^{\e,\g}|^{2}=0, \label{xegts}\\
&&\lim\limits_{l\rightarrow 0}\sup\limits_{s\in[0,T]}\mE\sup _{s \leqslant t \leqslant s+l}|\hat{X}_{t}^{\e,\g}-\hat{X}_{s}^{\e,\g}|^{2}=0.\label{hatxegts}
\ee
\el
\begin{proof}
Since the proofs of (\ref{xegts}) and (\ref{hatxegts}) are similar, we only prove (\ref{xegts}). 

First of all, the It\^o formula and (\ref{b1line}) imply that for $0\leq s<t\leq T$
\be
&&|X_{t}^{\e,\g}-X_{s}^{\e,\g}|^2\no\\
&=&-2\int_s^t\<X_{r}^{\e,\g}-X_{s}^{\e,\g},\dif K_r^{1,\e,\g}\>+2\int_s^t\<X_{r}^{\e,\g}-X_{s}^{\e,\g},b_{1}(r,X_{r}^{\e,\g},Y_{r}^{\e,\g})\>\dif r\no\\
&&+2\sqrt\e\int_s^t\<X_{r}^{\e,\g}-X_{s}^{\e,\g},\s_{1}(r,X_{r}^{\e,\g},Y_{r}^{\e,\g})\dif W^1_{r}\>+\e\int_s^t\|\s_{1}(r,X_{r}^{\e,\g},Y_{r}^{\e,\g})\|^2\dif r\no\\
&\leq&-2\int_s^t\<X_{r}^{\e,\g}-X_{s}^{\e,\g},\dif K_r^{1,\e,\g}\>+\int_s^t|X_{r}^{\e,\g}-X_{s}^{\e,\g}|^2\dif r+\bar{L}_{b_{1}, \s_{1}}\int_s^t(1+|X_{r}^{\e,\g}|^2+|Y_{r}^{\e,\g}|^2)\dif r\no\\
&&+2\sqrt\e\int_s^t\<X_{r}^{\e,\g}-X_{s}^{\e,\g},\s_{1}(r,X_{r}^{\e,\g},Y_{r}^{\e,\g})\dif W^1_{r}\>.
\label{xegstito}
\ee

Next, we compute $-2\int_s^t\<X_{r}^{\e,\g}-X_{s}^{\e,\g},\dif K_r^{1,\e,\g}\>$. Note that $0\in{\rm Int}(\cD(A_1))$. Thus, there is a $\t_0>0$ such that for any $R>0$ and $\t<\t_0$,
$$
\left\{x \in B(0,R): d(x,(\overline{\cD(A_1)})^c) \geqslant \t\right\} \neq \emptyset,
$$
where $B(0,R):=\{x\in\mR^n: |x|\leq R\}$, $d(\cdot,\cdot)$ is the Euclidean distance in $\mR^n$ and $(\overline{\cD(A_1)})^c$ denotes the complement of $\overline{\cD(A_1)}$. Set
$$
g_R(\t):=\sup \left\{|z|: z \in A_1(x) \text { for all } x \in B(0,R) \text { with } d\left(x,(\overline{\cD(A_1)})^c\right) \geqslant \t\right\},
$$
and by the local boundedness of $A_1$ on ${\rm Int}(\cD(A_1))$, it holds that
$$
g_R(\t)<+\infty.
$$
Again put for any $l>0$
$$
h_R(l):=\inf \left\{\t \in\left(0, \t_0\right): g_R(\t) \leqslant l^{-1 / 2}\right\}, 
$$
and we have that
$$
g_R\left(l+h_R(l)\right) \leqslant l^{-1 / 2} \text { and } \quad \lim _{l \downarrow 0} h_R(l)=0.
$$
Take $l_R>0$ be such that $l_R+h_R\left(l_R\right)<\t_0$. For $0<l<l_R \wedge 1$, let $X_s^{\varepsilon, \g, l, R}$ be the projection of $X_s^{\e,\g}$ on $\left\{x \in B(0,R): d\left(x,(\overline{\cD(A_1)})^c\right) \geqslant l+h_R(l)\right\}$. Thus, for $Z_s^{\varepsilon,\g, l, R} \in A_1(X_s^{\varepsilon,\g, l, R}), \sup\limits_{t\in[0,T]}|X_t^{\varepsilon,\g}| \leqslant R$ and $0<t-s<l$, it holds that
\ce
&&-2 \int_s^t\left\langle X_r^{\varepsilon,\g}-X_s^{\varepsilon,\g}, \dif K_r^{1,\e,\g} \right\rangle\\
 &=&-2 \int_s^t\left\langle X_r^{\varepsilon,\g}-X_s^{\varepsilon,\g, l, R}, \dif K_r^{1,\e,\g} \right\rangle-2 \int_s^t\left\langle X_s^{\varepsilon, \g, l, R}-X_s^{\varepsilon,\g}, \dif K_r^{1,\e,\g} \right\rangle \\ 
&\leqslant& -2 \int_s^t\left\langle X_r^{\varepsilon,\g}-X_s^{\varepsilon,\g, l, R}, Z_s^{\varepsilon,\g, l, R}\right\rangle \dif r+2\left(l+h_R(l)\right)\left|K^{1,\e,\g} \right|_0^T\\
&\leqslant& 4 l^{1 / 2} R+2\left(l+h_R(l)\right)\left|K^{1,\e,\g} \right|_0^T,
\de
and furthermore by (\ref{xegstito})
\ce
&&\sup _{s \leqslant t \leqslant s+l}\left|X_t^{\varepsilon,\g}-X_s^{\varepsilon,\g}\right|^2 I_{\{\sup\limits_{t\in[0,T]}|X_t^{\varepsilon,\g}| \leqslant R\}} \\
&\leqslant&\left(4 l^{1 / 2} R+2\left(l+h_R(l)\right)\left|K^{1,\e,\g} \right|_0^T\right)+\sup _{s \leqslant t \leqslant s+l}\int_s^t|X_{r}^{\e,\g}-X_{s}^{\e,\g}|^2\dif r\\
&&+\bar{L}_{b_{1}, \s_{1}}\sup _{s \leqslant t \leqslant s+l}\int_s^t\left(1+\left|X_r^{\varepsilon,\g}\right|^2+\left|Y_r^{\varepsilon,\g}\right|^2\right) \dif r\\
&&+C \sup _{s \leqslant t \leqslant s+l}\left|\int_s^t\left\langle X_r^{\varepsilon,\g}-X_s^{\varepsilon,\g}, \sigma_1\left(r,X_r^{\varepsilon,\g}, Y_r^{\varepsilon,\g}\right)\right\rangle \dif W^1_r\right|I_{\{\sup\limits_{t\in[0,T]}|X_t^{\varepsilon,\g}| \leqslant R\}}.
\de

From the above deduction and the BDG inequality, it follows that
\ce
&&\mE\sup _{s \leqslant t \leqslant s+l}\left|X_t^{\varepsilon,\g}-X_s^{\varepsilon,\g}\right|^2 I_{\{\sup\limits_{t\in[0,T]}|X_t^{\varepsilon,\g}| \leqslant R\}}\\
&\leq&\left(4 l^{1 / 2} R+2\left(l+h_R(l)\right)\mE|K^{1,\e,\g}|_0^T\right)+C(1+|x_0|^{2}+|y_0|^{2})l\\
&&+C\mE\left(\int_s^{s+l}|X_r^{\varepsilon,\g}-X_s^{\varepsilon,\g}|^2\|\sigma_1(r,X_r^{\varepsilon,\g}, Y_r^{\varepsilon,\g})\|^2 \dif r\right)^{1/2}I_{\{\sup\limits_{t\in[0,T]}|X_t^{\varepsilon,\g}| \leqslant R\}}\\
&\leq&\left(4 l^{1 / 2} R+2\left(l+h_R(l)\right)\mE|K^{1,\e,\g}|_0^T\right)+C(1+|x_0|^{2}+|y_0|^{2})l\\
&&+\frac{1}{2}\mE\sup _{s \leqslant t \leqslant s+l}\left|X_t^{\varepsilon,\g}-X_s^{\varepsilon,\g}\right|^2 I_{\{\sup\limits_{t\in[0,T]}|X_t^{\varepsilon,\g}| \leqslant R\}}+C(1+|x_0|^{2}+|y_0|^{2})l,
\de
which yields that
\ce
&&\mE\sup _{s \leqslant t \leqslant s+l}\left|X_t^{\varepsilon,\g}-X_s^{\varepsilon,\g}\right|^2 I_{\{\sup\limits_{t\in[0,T]}|X_t^{\varepsilon,\g}| \leqslant R\}}\\
&\leq& \left(4 l^{1 / 2} R+2\left(l+h_R(l)\right)\mE|K^{1,\e,\g}|_0^T\right)+C(1+|x_0|^{2}+|y_0|^{2})l.
\de
Based on this, letting $l\rightarrow 0$ and $R\rightarrow \infty$, one could conclude (\ref{xegts}), which completes the proof.
\end{proof}

\bl
 Suppose that the assumptions of Theorem \ref{xbarxp} hold. Then there exists a constant $C>0$ independent of $\e, \g$ such that 
 \be
\sup\limits_{t\in[0,T]}\mE|Y_{t}^{\e,\g}-\hat{Y}_{t}^{\e,\g}|^{2}\leq \frac{C}{\a}
\left(\d^2+\sup\limits_{s\in[0,T]}\mE\sup _{s \leqslant r \leqslant s+\d}|X_{r}^{\e,\g}-X_{s}^{\e,\g}|^{2}\right).
\label{unzt}
\ee
\el
\begin{proof}
First of all, by (\ref{Eq1}) and (\ref{hatz}), we have that 
\ce
Y_{t}^{\e,\g}-\hat{Y}_{t}^{\e,\g}
&=&-K_t^{2,\e,\g}+\hat{K}_t^{2,\e,\g}+\frac{1}{\g}\int_{0}^{t}\(b_{2}(s,X_{s}^{\e,\g},Y_{s}^{\e,\g})
-b_{2}(s(\d),X_{s(\d)}^{\e,\g},\hat{Y}_{s}^{\e,\g})\)\dif s\\
&&+\frac{1}{\sqrt{\g}}\int_{0}^{t}\(\s_{2}(s,X_{s}^{\e,\g},Y_{s}^{\e,\g})
-\s_{2}(s(\d),X_{s(\d)}^{\e,\g},\hat{Y}_{s}^{\e,\g})\)\dif W^2_{s}.
\de
Applying the It\^{o} formula to $|Y_{t}^{\e,\g}-\hat{Y}_{t}^{\e,\g}|^{2}e^{\l t}$ for any $\l>0$ and taking the expectation, by $(\mathbf{H}^{1}_{b_{2},\s_{2}})$ and $(\mathbf{H}^{2}_{b_{2},\s_{2}})$ one could obtain that
\ce
&&\mE|Y_{t}^{\e,\g}-\hat{Y}_{t}^{\e,\g}|^{2}e^{\l t}\\
&=&\l\mE\int_{0}^{t}|Y_{s}^{\e,\g}-\hat{Y}_{s}^{\e,\g}|^{2}e^{\l s}\dif s-\mE\int_{0}^{t}2e^{\l s}\<Y_{s}^{\e,\g}-\hat{Y}_{s}^{\e,\g},\dif (K_s^{2,\e,\g}-\hat{K}_s^{2,\e,\g})\>\\
&&+\frac{1}{\g}\mE\int_{0}^{t}2e^{\l s}\<Y_{s}^{\e,\g}-\hat{Y}_{s}^{\e,\g}, b_{2}(s,X_{s}^{\e,\g},Y_{s}^{\e,\g})
-b_{2}(s(\d),X_{s(\d)}^{\e,\g},\hat{Y}_{s}^{\e,\g})\>\dif s\\
&&+\frac{1}{\g}
\mE\int_{0}^{t}e^{\l s}\|\s_{2}(s,X_{s}^{\e,\g},Y_{s}^{\e,\g})-\s_{2}(s(\d),X_{s(\d)}^{\e,\g},\hat{Y}_{s}^{\e,\g})\|^{2}\dif s\\
&\leq&\l\mE\int_{0}^{t}|Y_{s}^{\e,\g}-\hat{Y}_{s}^{\e,\g}|^{2}e^{\l s}\dif s+\frac{1}{\g}\mE\int_{0}^{t}2e^{\l s}\<Y_{s}^{\e,\g}-\hat{Y}_{s}^{\e,\g}, b_{2}(s,X_{s}^{\e,\g},Y_{s}^{\e,\g})\\
&&-b_{2}(s,X_{s}^{\e,\g},\hat{Y}_{s}^{\e,\g})\>\dif s+\frac{1}{\g}\mE\int_{0}^{t}e^{\l s}\|\s_{2}(s,X_{s}^{\e,\g},Y_{s}^{\e,\g})-\s_{2}(s,X_{s}^{\e,\g},\hat{Y}_{s}^{\e,\g})\|^{2}\dif s\\
&&+\frac{1}{\g}\mE\int_{0}^{t}2e^{\l s}\<Y_{s}^{\e,\g}-\hat{Y}_{s}^{\e,\g}, b_{2}(s,X_{s}^{\e,\g},\hat{Y}_{s}^{\e,\g})-b_{2}(s(\d),X_{s(\d)}^{\e,\g},\hat{Y}_{s}^{\e,\g})\>\dif s\\
&&+\frac{1}{\g}\mE\int_{0}^{t}e^{\l s}\|\s_{2}(s,X_{s}^{\e,\g},Y_{s}^{\e,\g})-\s_{2}(s,X_{s}^{\e,\g},\hat{Y}_{s}^{\e,\g})\|^{2}\dif s\\
&&+\frac{1}{\g}\mE\int_{0}^{t}2e^{\l s}\|\s_{2}(s,X_{s}^{\e,\g},\hat{Y}_{s}^{\e,\g})-\s_{2}(s(\d),X_{s(\d)}^{\e,\g},\hat{Y}_{s}^{\e,\g})\|^{2}\dif s\\
&\leq&(\l-\frac{\b}{\g}+\frac{2L_{b_2,\s_2}}{\g})\mE\int_{0}^{t}|Y_{s}^{\e,\g}-\hat{Y}_{s}^{\e,\g}|^{2}e^{\l s}\dif s\\
&&+\frac{C}{\g}\mE\int_{0}^{t}e^{\l s}\(|s-s(\d)|^2+|X_{s}^{\e,\g}-X_{s(\d)}^{\e,\g}|^{2}\)\dif s.
\de

Finally, it follows from the above deduction and taking $\l=\frac{\a}{\g}$ that
\ce
\mE|Y_{t}^{\e,\g}-\hat{Y}_{t}^{\e,\g}|^{2}\leq\frac{C}{\a}\left(\d^2+\sup\limits_{s\in[0,T]}\mE\sup _{s \leqslant r \leqslant s+\d}|X_{r}^{\e,\g}-X_{s}^{\e,\g}|^{2}\right).
\de
The proof is complete.
\end{proof}

Finally, we estimate $\mE\sup\limits_{t\in[0,T]}|X_t^{\e,\g}-\hat{X}_t^{\e,\g}|^2$. 

\bl
Assume that the assumptions of Theorem \ref{xbarxp} hold. Then we have that 
\be
\mE\sup\limits_{t\in[0,T]}|X_t^{\e,\g}-\hat{X}_t^{\e,\g}|^2\leq C\left(\d^2+\sup\limits_{s\in[0,T]}\mE\sup _{s \leqslant r \leqslant s+\d}|X_{r}^{\e,\g}-X_{s}^{\e,\g}|^{2}\right),
\label{xehatxe}
\ee
where the constant $C>0$ is independent of $\e, \g$.
\el
\begin{proof}
Note that
\ce
X_{t}^{\e,\g}-\hat{X}_t^{\e,\g}=-K_t^{1,\e,\g}+\hat{K}_t^{1,\e,\g}+\int_0^t \left(b_{1}(s,X_{s}^{\e,\g},Y_{s}^{\e,\g})-b_1(s(\d),X^{\e,\g}_{s(\d)},\hat{Y}_s^{\e,\g})\right)\dif s.
\de
Thus, by the integration by parts, it holds that
\ce
&&|X_{t}^{\e,\g}-\hat{X}_t^{\e,\g}|^2\\
&=&-2\int_0^t<X_{s}^{\e,\g}-\hat{X}_s^{\e,\g},\dif (K_s^{1,\e,\g}-\hat{K}_s^{1,\e,\g})\>\\
&&+2\int_0^t<X_{s}^{\e,\g}-\hat{X}_s^{\e,\g},b_{1}(s,X_{s}^{\e,\g},Y_{s}^{\e,\g})-b_1(s(\d),X^{\e,\g}_{s(\d)},\hat{Y}_s^{\e,\g})\>\dif s\\
&\leq&\int_0^t|X_{s}^{\e,\g}-\hat{X}_s^{\e,\g}|^2\dif s+L_{b_1,\s_1}\int_0^t(|s-s(\d)|^2+|X^{\e,\g}_{s}-X^{\e,\g}_{s(\d)}|^2+|Y_s^{\e,\g}-\hat{Y}_s^{\e,\g}|^2)\dif s,
\de
and
\ce
\mE\sup\limits_{s\in[0,t]}|X_s^{\e,\g}-\hat{X}_s^{\e,\g}|^2\leq\int_0^t\mE\sup\limits_{s\in[0,r]}|X_{s}^{\e,\g}-\hat{X}_s^{\e,\g}|^2\dif r+TC\left(\d^2+\sup\limits_{s\in[0,T]}\mE\sup _{s \leqslant r \leqslant s+\d}|X_{r}^{\e,\g}-X_{s}^{\e,\g}|^{2}\right).
\de
where we use (\ref{unzt}). The Gronwall inequality yields the required result. The proof is completed.
\end{proof}

\subsubsection{Some estimates for the averaging equation (\ref{Eq3})}\label{aveequ}

\bl \label{averc}
Under the assumptions of Theorem \ref{xbarxp}, Eq.(\ref{Eq3}) has a unique solution $(\bar{X}_{\cdot},\bar{K}_{\cdot})$. Moreover, it holds that
\be
&&\sup\limits_{t\in[0,T]}|\bar{X}_{t}|^{2}\leq C(1+|x_0|^{2}), \label{barxb}\\
&&\lim\limits_{l\rightarrow 0}\sup\limits_{s\in[0,T]}\sup\limits_{s\leq t\leq s+l}|\bar{X}_{t}-\bar{X}_{s}|^2=0. \label{barxts}
\ee
\el
\begin{proof}
First of all, from (\ref{barb1lip}), it follows that Eq.(\ref{Eq3}) has a unique solution $(\bar{X}_{\cdot},\bar{K}_{\cdot})$. Then, since the proofs of the required estimates are similar to that for $X^{\e,\g}$ in Lemma \ref{xtztc} and Lemma \ref{xehatxets}, we omit them. The proof is complete.
\end{proof}

\bl\label{2orde}
Suppose that the assumptions of Theorem \ref{xbarxp} hold. Then there exists a constant $C>0$ independent of $\e, \g$ such that 
\ce
\mE\(\sup_{0\leq t\leq T}|\hat{X}_{t}^{\e,\g}-\bar{X_{t}}|^{2}\)\leq C\e+\Sigma(\g),
\de
where 
\ce
\Sigma(\g)&:=&C\left(\sup\limits_{s\in[0,T]}\mE\sup _{s \leqslant r \leqslant s+\d}|\hat{X}_{r}^{\e,\g}-\hat{X}_{s}^{\e,\g}|^{2}+\sup\limits_{s\in[0,T]}\sup\limits_{s\leq r\leq s+\d}|\bar{X}_{r}-\bar{X}_{s}|^2\right)^{1/2}\\
&&+C\((\frac{\g}{\d})^{1/2}+\d^{1/2}\)+C\d^2+C\left(\sup\limits_{s\in[0,T]}\mE\sup _{s \leqslant r \leqslant s+\d}|X_{r}^{\e,\g}-X_{s}^{\e,\g}|^{2}\right)\\
&&+C\sup\limits_{s\in[0,T]}\sup\limits_{s\leq r\leq s+\d}|\bar{X}_{r}-\bar{X}_{s}|^2.
\de
\el
\begin{proof}
Note that
\ce
\hat{X}_{t}^{\e,\g}-\bar{X}_{t}
=-\hat{K}_{t}^{1,\e,\g}+\bar{K}_{t}+\int_{0}^{t}\(b_{1}(s(\d),X_{s(\d)}^{\e,\g},\hat{Y}_{s}^{\e,\g})-\bar{b}_{1}(s,\bar{X}_{s})\)\dif s
+\sqrt{\e}\int_{0}^{t}\s_{1}(s,X_{s}^{\e,\g},Y_{s}^{\e,\g})\dif W^1_{s}.
\de
Thus, by the It\^o formula, it holds that
\ce
&&|\hat{X}_{t}^{\e,\g}-\bar{X}_{t}|^2\\
&=&-2\int_0^t\<\hat{X}_{s}^{\e,\g}-\bar{X}_{s},\dif(\hat{K}_{s}^{1,\e,\g}-\bar{K}_{s})\>+2\int_0^t\<\hat{X}_{s}^{\e,\g}-\bar{X}_{s},b_{1}(s(\d),X_{s(\d)}^{\e,\g},\hat{Y}_{s}^{\e,\g})-\bar{b}_{1}(s,\bar{X}_{s})\>\dif s\\
&&+2\sqrt{\e}\int_0^t\<\hat{X}_{s}^{\e,\g}-\bar{X}_{s},\s_{1}(s,X_{s}^{\e,\g},Y_{s}^{\e,\g})\dif W^1_{s}\>+\e\int_{0}^{t}\|\s_{1}(s,X_{s}^{\e,\g},Y_{s}^{\e,\g})\|^2\dif s\\
&\leq&2\int_0^t\<\hat{X}_{s}^{\e,\g}-\bar{X}_{s},b_{1}(s(\d),X_{s(\d)}^{\e,\g},\hat{Y}_{s}^{\e,\g})-\bar{b}_{1}(s,\bar{X}_{s})\>\dif s\\
 &&+2\sqrt{\e}\int_0^t\<\hat{X}_{s}^{\e,\g}-\bar{X}_{s},\s_{1}(s,X_{s}^{\e,\g},Y_{s}^{\e,\g})\dif W^1_{s}\>+\e \bar L_{b_1,\s_1}\int_{0}^{t}(1+|X_{s}^{\e,\g}|^2+|Y_{s}^{\e,\g}|^2)\dif s.
\de
Moreover, based on the BDG inequality, we get that
\be
&&\mE\(\sup_{0\leq t\leq T}|\hat{X}_{t}^{\e,\g}-\bar{X}_{t}|^{2}\)\no\\
&\leq&2\mE\sup_{0\leq t\leq T}\left|\int_0^t\<\hat{X}_{s}^{\e,\g}-\bar{X}_{s},b_{1}(s(\d),X_{s(\d)}^{\e,\g},\hat{Y}_{s}^{\e,\g})-\bar{b}_{1}(s,\bar{X}_{s})\>\dif s\right|\no\\
&&+2\sqrt{\e}\mE\sup_{0\leq t\leq T}\left|\int_0^t\<\hat{X}_{s}^{\e,\g}-\bar{X}_{s},\s_{1}(s,X_{s}^{\e,\g},Y_{s}^{\e,\g})\dif W^1_{s}\>\right|+C\e\no\\
&\leq&2\mE\sup_{0\leq t\leq T}\left|\int_0^t\<\hat{X}_{s}^{\e,\g}-\bar{X}_{s},b_{1}(s(\d),X_{s(\d)}^{\e,\g},\hat{Y}_{s}^{\e,\g})-\bar{b}_{1}(s,\bar{X}_{s})\>\dif s\right|\no\\
&&+2\sqrt{\e}C\mE\left(\int_0^T|\hat{X}_{s}^{\e,\g}-\bar{X}_{s}|^2\|\s_{1}(s,X_{s}^{\e,\g},Y_{s}^{\e,\g})\|^2\dif s\right)^{1/2}+C\e\no\\
&\leq&2\mE\sup_{0\leq t\leq T}\left|\int_0^t\<\hat{X}_{s}^{\e,\g}-\bar{X}_{s},b_{1}(s(\d),X_{s(\d)}^{\e,\g},\hat{Y}_{s}^{\e,\g})-\bar{b}_{1}(s,\bar{X}_{s})\>\dif s\right|\no\\
&&+\frac{1}{2}\mE\(\sup_{0\leq t\leq T}|\hat{X}_{t}^{\e,\g}-\bar{X}_{t}|^{2}\)+C\e.
\label{hatxbarx}
\ee

Next, set
$$
I:=\mE\sup_{0\leq t\leq T}\left|\int_0^t\<\hat{X}_{s}^{\e,\g}-\bar{X}_{s},b_{1}(s(\d),X_{s(\d)}^{\e,\g},\hat{Y}_{s}^{\e,\g})-\bar{b}_{1}(s,\bar{X}_{s})\>\dif s\right|,
$$
and we are devoted to estimating $I$. Note that
\be
I&\leq&\mE\sup_{0\leq t\leq T}\left|\int_0^t\<\hat{X}_{s}^{\e,\g}-\hat{X}_{s(\d)}^{\e,\g}-\bar{X}_{s}+\bar{X}_{s(\d)},b_{1}(s(\d),X_{s(\d)}^{\e,\g},\hat{Y}_{s}^{\e,\g})-\bar{b}_{1}(s,\bar{X}_{s})\>\dif s\right|\no\\
&&+\mE\sup_{0\leq t\leq T}\left|\int_0^t\<\hat{X}_{s(\d)}^{\e,\g}-\bar{X}_{s(\d)},b_{1}(s(\d),X_{s(\d)}^{\e,\g},\hat{Y}_{s}^{\e,\g})-\bar{b}_{1}(s(\d),X_{s(\d)}^{\e,\g})\>\dif s\right|\no\\
&&+\mE\sup_{0\leq t\leq T}\left|\int_0^t\<\hat{X}_{s(\d)}^{\e,\g}-\bar{X}_{s(\d)},\bar{b}_{1}(s(\d),X_{s(\d)}^{\e,\g})-\bar{b}_{1}(s(\d),\bar{X}_{s(\d)})\>\dif s\right|\no\\
&&+\mE\sup_{0\leq t\leq T}\left|\int_0^t\<\hat{X}_{s(\d)}^{\e,\g}-\bar{X}_{s(\d)},\bar{b}_{1}(s(\d),\bar{X}_{s(\d)})-\bar{b}_{1}(s,\bar{X}_{s})\>\dif s\right|\no\\
&=:&I_{1}+I_{2}+I_{3}+I_{4}.
\label{i1234}
\ee
Thus, we estimate $I_{1}, I_{2}, I_{3}, I_{4}$, respectively.

For $I_{1}$, by the H\"older inequality and the linear growth for $b_1, \bar{b}_{1}$, it holds that
\be
I_{1}&\leq& C\left(\int_0^T\mE(|\hat{X}_{s}^{\e,\g}-\hat{X}_{s(\d)}^{\e,\g}|^2+|\bar{X}_{s}-\bar{X}_{s(\d)}|^2)\dif s\right)^{1/2}\no\\
&&\times\left(\int_0^T\mE(1+|X_{s(\d)}^{\e,\g}|^2+|\hat{Y}_{s}^{\e,\g}|^2+|\bar{X}_{s}|^2)\dif s\right)^{1/2}\no\\
&\leq& C\left(\sup\limits_{s\in[0,T]}\mE\sup _{s \leqslant r \leqslant s+\d}|\hat{X}_{r}^{\e,\g}-\hat{X}_{s}^{\e,\g}|^{2}+\sup\limits_{s\in[0,T]}\sup\limits_{s\leq r\leq s+\d}|\bar{X}_{r}-\bar{X}_{s}|^2\right)^{1/2}.
\label{i1}
\ee

For $I_{2}$, by the deduction in {\bf Step 2.}, we know that 
\be
I_{2}\leq C\((\frac{\g}{\d})^{1/2}+\d^{1/2}\).
\label{i2}
\ee

For $I_{3}$, the  Lipschitz continuity of $\bar{b}_{1}$ and (\ref{xehatxe}) imply that
\be
I_{3}&\leq& C\mE\int_0^T|\hat{X}_{s(\d)}^{\e,\g}-\bar{X}_{s(\d)}|(|X_{s(\d)}^{\e,\g}-\hat{X}_{s(\d)}^{\e,\g}|+|\hat{X}_{s(\d)}^{\e,\g}-\bar{X}_{s(\d)}|)\dif s\no\\
&\leq&C\int_0^T\mE|\hat{X}_{s(\d)}^{\e,\g}-\bar{X}_{s(\d)}|^2\dif s+C\int_0^T\mE|X_{s(\d)}^{\e,\g}-\hat{X}_{s(\d)}^{\e,\g}|^2\dif s\no\\
&\leq&C\int_0^T\mE|\hat{X}_{s(\d)}^{\e,\g}-\bar{X}_{s(\d)}|^2\dif s+C\left(\d^2+\sup\limits_{s\in[0,T]}\mE\sup _{s \leqslant r \leqslant s+\d}|X_{r}^{\e,\g}-X_{s}^{\e,\g}|^{2}\right).
\label{i3}
\ee

For $I_{4}$, by the  Lipschitz continuity of $\bar{b}_{1}$ and (\ref{barxts}), it holds that
\be
I_{4}&\leq& C\int_0^T\mE|\hat{X}_{s(\d)}^{\e,\g}-\bar{X}_{s(\d)}|^2\dif s+C\int_0^T\(\d^2+|\bar{X}_{s}-\bar{X}_{s(\d)}|^2\)\dif s\no\\
&\leq&C\int_0^T\mE|\hat{X}_{s(\d)}^{\e,\g}-\bar{X}_{s(\d)}|^2\dif s+C\d^2T+CT\sup\limits_{s\in[0,T]}\sup\limits_{s\leq r\leq s+\d}|\bar{X}_{r}-\bar{X}_{s}|^2.
\label{i4}
\ee

Combining (\ref{i1})-(\ref{i4}) with (\ref{i1234}), we have that
\ce
I\leq\Sigma(\g)+C\int_0^T\mE|\hat{X}_{s(\d)}^{\e,\g}-\bar{X}_{s(\d)}|^2\dif s,
\de
which together with (\ref{hatxbarx}) implies that
\ce
\mE\(\sup_{0\leq t\leq T}|\hat{X}_{t}^{\e,\g}-\bar{X}_{t}|^{2}\)\leq C\e+\Sigma(\g)+C\int_0^T\mE\(\sup_{0\leq r\leq s}|\hat{X}_{r}^{\e,\g}-\bar{X}_{r}|^{2}\)\dif s.
\de
Finally, the Gronwall inequality yields the required estimate.

{\bf Step 2.} We prove (\ref{i2}).

For $I_2$, it holds that
\ce
I_{2}&=&\Bigg(\mE\sup_{0\leq t\leq T}
 \Big|\int_{0}^{[\frac{t}{\d}]\d}\<\hat{X}_{s(\d)}^{\e,\g}-\bar{X}_{s(\d)},b_{1}(s(\d),X_{s(\d)}^{\e,\g},\hat{Y}_{s}^{\e,\g})
 -\bar{b}_{1}(s(\d),X_{s(\d)}^{\e,\g})\>\dif s\no\\
 &&\quad\quad\quad\quad+\int_{[\frac{t}{\d}]\d}^{t}\<\hat{X}_{s(\d)}^{\e,\g}-\bar{X}_{s(\d)},b_{1}(s(\d),X_{s(\d)}^{\e,\g},\hat{Y}_{s}^{\e,\g})
 -\bar{b}_{1}(s(\d),X_{s(\d)}^{\e,\g})\>\dif s\Big|\Bigg)\no\\
&\leq&\Bigg(\mE\sup_{0\leq t\leq T}
\Big|\int_{0}^{[\frac{t}{\d}]\d}\<\hat{X}_{s(\d)}^{\e,\g}-\bar{X}_{s(\d)},b_{1}(s(\d),X_{s(\d)}^{\e,\g},\hat{Y}_{s}^{\e,\g})
 -\bar{b}_{1}(s(\d),X_{s(\d)}^{\e,\g})\>\dif s\Big|\Bigg)\no\\
&&+\Bigg(\mE\sup_{0\leq t\leq T}
\Big|\int_{[\frac{t}{\d}]\d}^{t}\<\hat{X}_{s(\d)}^{\e,\g}-\bar{X}_{s(\d)},b_{1}(s(\d),X_{s(\d)}^{\e,\g},\hat{Y}_{s}^{\e,\g})
 -\bar{b}_{1}(s(\d),X_{s(\d)}^{\e,\g})\>\dif s\Big|\Bigg)\no\\
 &=:&I_{21}+I_{22}.
\de

Next, we estimate $I_{21}$. Note that
\be
I_{21}
&=&\Bigg(\mE\sup_{0\leq t\leq T}
\Big|\int_{0}^{[\frac{t}{\d}]\d}\<\hat{X}_{s(\d)}^{\e,\g}-\bar{X}_{s(\d)},b_{1}(s(\d),X_{s(\d)}^{\e,\g},\hat{Y}_{s}^{\e,\g})
 -\bar{b}_{1}(s(\d),X_{s(\d)}^{\e,\g})\>\dif s\Big|\Bigg)\no\\
&=&\mE\Bigg(\sup_{0\leq t\leq T}
\Big|\sum\limits_{k=0}^{[\frac{t}{\d}]-1}\int_{k\d}^{(k+1)\d}\<\hat{X}_{s(\d)}^{\e,\g}-\bar{X}_{s(\d)},b_{1}(s(\d),X_{s(\d)}^{\e,\g},\hat{Y}_{s}^{\e,\g})
 -\bar{b}_{1}(s(\d),X_{s(\d)}^{\e,\g})\>\dif s\Big|\Bigg)\no\\
&\leq&\mE\Bigg(\sup_{0\leq t\leq T}\sum\limits_{k=0}^{[\frac{t}{\d}]-1}\left|\int_{k\d}^{(k+1)\d}\<\hat{X}_{k\d}^{\e,\g}-\bar{X}_{k\d},b_{1}(k\d,X_{k\d}^{\e,\g},\hat{Y}_{s}^{\e,\g})
 -\bar{b}_{1}(k\d,X_{k\d}^{\e,\g})\>\dif s\right|\Bigg)\no\\
&\leq&\sum\limits_{k=0}^{[\frac{T}{\d}]-1}
\mE\Bigg(\Big|\int_{k\d}^{(k+1)\d}\<\hat{X}_{k\d}^{\e,\g}-\bar{X}_{k\d},b_{1}(k\d,X_{k\d}^{\e,\g},\hat{Y}_{s}^{\e,\g})
 -\bar{b}_{1}(k\d,X_{k\d}^{\e,\g})\>\dif s\Big|\Bigg)\no\\
&\leq&[\frac{T}{\d}]\sup_{0\leq k\leq [\frac{T}{\d}]-1}
\mE\Bigg(\Big|\int_{k\d}^{(k+1)\d}\<\hat{X}_{k\d}^{\e,\g}-\bar{X}_{k\d},b_{1}(k\d,X_{k\d}^{\e,\g},\hat{Y}_{s}^{\e,\g})
 -\bar{b}_{1}(k\d,X_{k\d}^{\e,\g})\>\dif s\Big|\Bigg)\no\\
&\leq&\g(\frac{T}{\d})\sup_{0\leq k\leq [\frac{T}{\d}]-1}
\mE\Bigg(\Big|\<\hat{X}_{k\d}^{\e,\g}-\bar{X}_{k\d},\int_{0}^{\d/\g}(b_{1}(k\d,X_{k\d}^{\e,\g},\hat{Y}_{\g s+k\d}^{\e,\g})
 -\bar{b}_{1}(k\d,X_{k\d}^{\e,\g}))\dif s\>\Big|\Bigg)\no\\
&\leq&\g(\frac{T}{\d})\sup_{0\leq k\leq [\frac{T}{\d}]-1}
\mE|\hat{X}_{k\d}^{\e,\g}-\bar{X}_{k\d}|\left|\int_{0}^{\d/\g}(b_{1}(k\d,X_{k\d}^{\e,\g},\hat{Y}_{\g s+k\d}^{\e,\g})
 -\bar{b}_{1}(k\d,X_{k\d}^{\e,\g}))\dif s\right|\no\\
 &\leq&\g(\frac{T}{\d})\sup_{0\leq k\leq [\frac{T}{\d}]-1}(\mE|\hat{X}_{k\d}^{\e,\g}-\bar{X}_{k\d}|^2)^{1/2}\no\\
 &&\times\left(\mE\left|\int_{0}^{\d/\g}(b_{1}(k\d,X_{k\d}^{\e,\g},\hat{Y}_{\g s+k\d}^{\e,\g})
 -\bar{b}_{1}(k\d,X_{k\d}^{\e,\g}))\dif s\right|^2\right)^{1/2}.
\label{b41c}
\ee
Thus, we compute $\mE\left|\int_{0}^{\d/\g}(b_{1}(k\d,X_{k\d}^{\e,\g},\hat{Y}_{\g s+k\d}^{\e,\g})-\bar{b}_{1}(k\d,X_{k\d}^{\e,\g}))\dif s\right|^2$. It is easy to see that
\ce
&&\mE\left|\int_{0}^{\d/\g}(b_{1}(k\d,X_{k\d}^{\e,\g},\hat{Y}_{\g s+k\d}^{\e,\g})-\bar{b}_{1}(k\d,X_{k\d}^{\e,\g}))\dif s\right|^2\\
&=&2\mE\int_{0}^{\d/\g}\int_{r}^{\d/\g}\<b_{1}(k\d,X_{k\d}^{\e,\g},\hat{Y}_{\g s+k\d}^{\e,\g})-\bar{b}_{1}(k\d,X_{k\d}^{\e,\g}),b_{1}(k\d,X_{k\d}^{\e,\g},\hat{Y}_{\g r+k\d}^{\e,\g})-\bar{b}_{1}(k\d,X_{k\d}^{\e,\g})\>\dif s\dif r.
\de
Therefore, set for $0<r<s\leq\d/\g$
$$
\Phi(s,r):=\mE\<b_{1}(k\d,X_{k\d}^{\e,\g},\hat{Y}_{\g s+k\d}^{\e,\g})-\bar{b}_{1}(k\d,X_{k\d}^{\e,\g}),b_{1}(k\d,X_{k\d}^{\e,\g},\hat{Y}_{\g r+k\d}^{\e,\g})-\bar{b}_{1}(k\d,X_{k\d}^{\e,\g})\>,
$$
and we only investigate $\Phi(s,r)$.

For any $w>0$ and random variables $\xi, \zeta\in L^2(\Omega,\mathscr{F}_w,\mP), \xi\in\overline{\cD(A_1)}, \zeta\in\overline{\cD(A_2)}$, we consider the following equation
\ce\left\{\begin{array}{l}
\dif \check{Y}_t^{\g, w, \xi,\zeta}\in -A_2(\check{Y}_t^{\g, w, \xi,\zeta})\dif t+\frac{1}{\g}b_2(w,\xi, \check{Y}_t^{\g, w, \xi,\zeta})\dif t+\frac{1}{\sqrt{\g}}\s_2(w,\xi, \check{Y}_t^{\g, w, \xi,\zeta})\dif W^2_t,\quad t\geq w,\\
\check{Y}_w^{\g, w, \xi,\zeta}=\zeta.
\end{array}
\right.
\de
Under $(\mathbf{H}_{A_{2}})$ and $(\mathbf{H}^1_{b_{2}, \s_{2}})$, the above multivalued SDE has a unique solution $(\check{Y}_t^{\g, w, \xi,\zeta}, \check{K}_t^{2,\g, w, \xi,\zeta})$ (\cite{rwzx}). Then it holds that
$$
\hat{Y}_t^{\e,\g}=\check{Y}_t^{\g, k \delta, X_{k \delta}^{\e,\g}, \hat{Y}_{k \delta}^{\e,\g}}, \quad \hat{K}_t^{2,\e,\g}=\check{K}_t^{2,\g, k \delta, X_{k \delta}^{\e,\g}, \hat{Y}_{k \delta}^{\e,\g}}, \quad t \in[k \delta,(k+1) \delta],
$$
and furthermore
\ce
\Phi(s,r)=\mE\<b_{1}(k\d,X_{k\d}^{\e,\g},\check{Y}^{\g, k \delta, X_{k \delta}^{\e,\g}, \hat{Y}_{k \delta}^{\e,\g}}_{\g s+k\d})
 -\bar{b}_{1}(k\d,X_{k\d}^{\e,\g}),b_{1}(k\d,X_{k\d}^{\e,\g},\check{Y}^{\g, k \delta, X_{k \delta}^{\e,\g}, \hat{Y}_{k \delta}^{\e,\g}}_{\g r+k\d})
 -\bar{b}_{1}(k\d,X_{k\d}^{\e,\g})\>.
\de
Note that $X_{k \delta}^{\e,\g}, \hat{Y}_{k \delta}^{\e,\g}$ are $\sF_{k\d}$-measurable, and for any $x\in\overline{\cD(A_1)}, y\in\overline{\cD(A_2)}$, $\check{Y}_t^{\g,k\d, x,y}$ is independent of $\sF_{k\d}$. Thus, we have that
\ce
\Phi(s,r)
 &=&\mE\Bigg[\mE\Bigg[\<b_{1}(k\d,X_{k\d}^{\e,\g},\check{Y}^{\g, k \delta, X_{k \delta}^{\e,\g}, \hat{Y}_{k \delta}^{\e,\g}}_{\g s+k\d})
 -\bar{b}_{1}(k\d,X_{k\d}^{\e,\g}),\\
 &&\qquad\qquad b_{1}(k\d,X_{k\d}^{\e,\g},\check{Y}^{\g, k \delta, X_{k \delta}^{\e,\g}, \hat{Y}_{k \delta}^{\e,\g}}_{\g r+k\d})
 -\bar{b}_{1}(k\d,X_{k\d}^{\e,\g})\>\Bigg{|}\sF_{k\d}\Bigg]\Bigg]\\
&=&\mE\Bigg[\mE\Bigg[\<b_{1}(k\d,x,\check{Y}^{\g, k \delta, x, y}_{\g s+k\d})
 -\bar{b}_{1}(k\d,x),b_{1}(k\d,x,\check{Y}^{\g, k \delta, x, y}_{\g r+k\d})
 -\bar{b}_{1}(k\d,x)\>\Bigg]\Bigg{|}_{(x,y)=(X_{k\d}^{\e,\g},\hat{Y}_{k \delta}^{\e,\g})}\Bigg].
\de

Here, we investigate $\check{Y}^{\g, k \delta, x, y}_{\g s+k\d}$. On one hand, it holds that
\ce
\check{Y}^{\g, k \delta, x, y}_{\g s+k\d}&=&y-\check{K}_{\g s+k\d}^{2,\g, k\d, x,y}+\check{K}_{k\d}^{2,\g, k\d, x,y}+\frac{1}{\g} \int_{k\d}^{\g s+k\d}b_2(k\d,x, \check{Y}^{\g, k \delta, x, y}_{r})\dif r\\
&&+\frac{1}{\sqrt{\g}} \int_{k\d}^{\g s+k\d} \s_2(k\d, x, \check{Y}^{\g, k \delta, x, y}_{r})\dif W^2_r\\
&=&y-\check{K}_{\g s+k\d}^{2,\g, k\d, x,y}+\check{K}_{k\d}^{2,\g, k\d, x,y}+\frac{1}{\g} \int_{0}^{\g s}b_2(k\d, x, \check{Y}^{\g, k \delta, x, y}_{\tau+k\d})\dif \tau\\
&&+\frac{1}{\sqrt{\g}} \int_{0}^{\g s} \s_2(k\d, x, \check{Y}^{\g, k \delta, x, y}_{\tau+k\d})\dif \tilde{W}^2_\tau\\
&=&y-\check{\tilde{\check{K}}}_{s}^{2,\g, k\d, x,y}+\int_{0}^{s}b_2(k\d, x, \check{Y}^{\g, k \delta, x, y}_{\g v+k\d})\dif v+\int_{0}^{s} \s_2(k\d, x, \check{Y}^{\g, k \delta, x, y}_{\g v+k\d})\dif \check{\tilde{W}}^2_v,
\de
where $\tilde{W}^2_\cdot:=W^2_{\cdot+k\d}-W^2_{k\d}$ and $\check{\tilde{W}}^2_\cdot:=\frac{1}{\sqrt{\g}}\tilde{W}^2_{\g \cdot}$ are two $m$-dimensional standard Brownian motions, and $\check{\tilde{\check{K}}}_{s}^{2,\g, k\d, x,y}:=\check{K}_{\g s+k\d}^{2,\g, k\d, x,y}-\check{K}_{k\d}^{2,\g, k\d, x,y}$. On the other hand, note that the frozen equation (\ref{Eq2}) is written as
\ce
Y_{s}^{k\d,x, y}=y-K^{2,k\d,x, y}_s+\int_0^s b_{2}(k\d, x,Y_{r}^{k\d,x, y})\dif r+\int_0^s\s_{2}(k\d, x,Y_{r}^{k\d,x, y})\dif W^2_{r}.
\de
Thus, for $s\in[0,\d/\g]$, $\check{Y}^{\g, k \delta, x, y}_{\g s+k\d}$ and $Y_{s}^{k\d,x, y}$ have the same distribution, which implies that
\ce
&&\mE\Bigg[\<b_{1}(k\d,x,\check{Y}^{\g, k \delta, x, y}_{\g s+k\d})-\bar{b}_{1}(k\d,x),b_{1}(k\d,x,\check{Y}^{\g, k \delta, x, y}_{\g r+k\d})-\bar{b}_{1}(k\d,x)\>\Bigg]\\
&=&\mE\<b_{1}(k\d,x,Y_{s}^{k\d,x, y})-\bar{b}_{1}(k\d,x),b_{1}(k\d,x,Y_{r}^{k\d,x, y})-\bar{b}_{1}(k\d,x)\>\\
&=&\mE\left[\mE\left[\Bigg<b_{1}(k\d,x,Y_{s}^{k\d,x, y})-\bar{b}_{1}(k\d,x), b_{1}(k\d,x,Y_{r}^{k\d,x, y})
-\bar{b}_{1}(k\d,x)\Bigg>\Bigg{|}\sF_r^{W^2}\right]\right]\\
&=&\mE\left[\Bigg<\mE\left[b_{1}(k\d,x,Y_{s}^{k\d,x, y})\Bigg{|}\sF_r^{W^2}\right]-\bar{b}_{1}(k\d,x),b_{1}(k\d,x,Y_{r}^{k\d,x, y})
-\bar{b}_{1}(k\d,x)\Bigg>\right]\\
&\leq&\left(\mE\left|\mE\left[b_{1}(k\d,x,Y_{s}^{k\d,x, y})\Bigg{|}\sF_r^{W^2}\right]-\bar{b}_{1}(k\d,x)\right|^2\right)^{1/2}\left(\mE|b_{1}(k\d,x,Y_{r}^{k\d,x, y})
-\bar{b}_{1}(k\d,x)|^2\right)^{1/2},
\de
where $\mathscr{F}_{r}^{W^2}=\sigma\{W^2_{\tau},0\leq \tau \leq r\}\vee \cN$, and $\cN$ denotes  the collection of all $\mP$-zero sets. Moreover, based on (\ref{meu2}), (\ref{memu2}), (\ref{inu2}), we obtain that
\ce
&&\left(\mE\left|\mE\left[b_{1}(k\d,x,Y_{s}^{k\d,x,y})\Bigg{|}\sF_r^{W^2}\right]-\bar{b}_{1}(k\d,x)\right|^2\right)^{1/2}\\
&=&\left(\mE\left|\mE\left[b_{1}(k\d,x,Y_{s-r}^{k\d,x,\hat{y}})\right]|_{\hat{y}=Y_{r}^{k\d,x,y}}-\bar{b}_{1}(k\d,x)\right|^2\right)^{1/2}\\
&\leq&\(Ce^{-\a(s-r)}(1+|x|^{2}+\mE|Y_{r}^{k\d,x,y}|^{2})\)^{\frac{1}{2}}\\
&\leq&\(Ce^{-\a(s-r)}(1+|x|^{2}+|y|^{2})\)^{\frac{1}{2}}\\
&\leq& Ce^{-\a(s-r)/2}(1+|x|+|y|),
\de
and
\ce
&&\left(\mE|b_{1}(k\d,x,Y_{r}^{k\d,x,y})-\bar{b}_{1}(k\d,x)|^2\right)^{1/2}\\
&=&\left(\mE|b_{1}(k\d,x,Y_{r}^{k\d,x,y})-\int_{\overline{\cD(A_2)}}b_{1}(k\d,x,z)\nu^{k\d,x}(\dif z)|^2\right)^{1/2}\\
&\leq&\left(\mE\int_{\overline{\cD(A_2)}}|b_{1}(k\d,x,Y_{r}^{k\d,x,y})-b_{1}(k\d,x,z)|^2\nu^{k\d,x}(\dif z)\right)^{1/2}\\
&\leq&C\left(\int_{\overline{\cD(A_2)}}\mE|Y_{r}^{k\d,x,y}-z|^2\nu^{k\d,x}(\dif z)\right)^{1/2}\\
&\leq&C\left(|y|^{2}e^{-\frac{\a}{2} r}+C(1+|x|^{2})\right)^{1/2}\\
&\leq&C(1+|x|+|y|).
\de

Combining the above deduction, by (\ref{xeb}), (\ref{hatzb}) one can have that
\ce
\Phi(s,r)\leq Ce^{-\a(s-r)/2}.
\de
Inserting the above inequality in (\ref{b41c}), we get that
\be
I_{21}\leq C\g(\frac{T}{\d})\sup_{0\leq k\leq [\frac{T}{\d}]-1}\left(\int_{0}^{\frac{\d}{\g}}\int_{r}^{\frac{\d}{\g}}Ce^{-\a(s-r)/2}\dif s\dif r\right)^{1/2}\leq C(\frac{\g}{\d})^{1/2}.
\label{b4de1}
\ee

Finally, we estimate $I_{22}$. By (\ref{b1line}), (\ref{xeb}), (\ref{hatzb}), (\ref{hatxb}), (\ref{barxb})  and the H\"older inequality, one could get that
\ce
I_{22}
&\leq& \Bigg(\mE\sup_{0\leq t\leq T}
\int_{[\frac{t}{\d}]\d}^{t}|\hat{X}_{s(\d)}^{\e,\g}-\bar{X}_{s(\d)}||b_{1}(s(\d),X_{s(\d)}^{\e,\g},\hat{Y}_{s}^{\e,\g})
 -\bar{b}_{1}(s(\d),X_{s(\d)}^{\e,\g})|\dif s\Bigg)\\
 &\leq& \d^{1/2}\Bigg(\mE\sup_{0\leq s\leq T}|\hat{X}_{s}^{\e,\g}-\bar{X}_{s}|^2\Bigg)^{1/2}\Bigg(\mE\sup_{0\leq t\leq T}\int_{[\frac{t}{\d}]\d}^{t}|b_{1}(s(\d),X_{s(\d)}^{\e,\g},\hat{Y}_{s}^{\e,\g})-\bar{b}_{1}(s(\d),X_{s(\d)}^{\e,\g})|^2\dif s\Bigg)^{1/2}\\
 &\leq& \d^{1/2}\Bigg(\mE\sup_{0\leq s\leq T}|\hat{X}_{s}^{\e,\g}-\bar{X}_{s}|^2\Bigg)^{1/2}\Bigg(\mE\int_{0}^{T}|b_{1}(s(\d),X_{s(\d)}^{\e,\g},\hat{Y}_{s}^{\e,\g})-\bar{b}_{1}(s(\d),X_{s(\d)}^{\e,\g})|^2\dif s\Bigg)^{1/2}\\
 &\leq& C\d^{1/2}\Bigg(\mE\sup_{0\leq s\leq T}|\hat{X}_{s}^{\e,\g}|^2+\sup_{0\leq s\leq T}|\bar{X}_{s}|^2\Bigg)^{1/2}\Bigg(\int_{0}^{T}(1+\mE|X_{s(\d)}^{\e,\g}|^2+\mE|\hat{Y}_{s}^{\e,\g}|^2)\dif s\Bigg)^{1/2}\\
&\leq& C\d^{1/2},
\de
which together with (\ref{b4de1}) implies (\ref{i2}). The proof is complete.
\end{proof}

At present, we are ready to prove Theorem \ref{xbarxp}.

{\bf Proof of Theorem \ref{xbarxp}.} By (\ref{xehatxe}) and Lemma \ref{2orde}, we get that 
\ce
\mE\(\sup_{0\leq t\leq T}|X_{t}^{\e,\g}-\bar{X_{t}}|^{2}\)&\leq& 2\mE\(\sup_{0\leq t\leq T}|X_{t}^{\e,\g}-\hat{X}_{t}^{\e,\g}|^{2}\)+2\mE\(\sup_{0\leq t\leq T}|\hat{X}_{t}^{\e,\g}-\bar{X_{t}}|^{2}\)\\
&\leq&C\left(\sup\limits_{s\in[0,T]}\mE\sup _{s \leqslant r \leqslant s+\d}|X_{r}^{\e,\g}-X_{s}^{\e,\g}|^{2}\right)+C\e+\Sigma(\g).
\de
Let us take $\d=\g^{\iota}$ for $0<\iota<1$. Thus, as $\e$ tends to $0$, $\g\rightarrow 0$, $\d\rightarrow 0, \g/\d\rightarrow 0$ and (\ref{xegts}), (\ref{hatxegts}), (\ref{barxts}) imply the required limit, which completes the proof.

\subsection{Proof of Corollary \ref{convrate}}

In this subsection, we prove Corollary \ref{convrate}. We start with a key estimate.

\bl\label{convrateproo}
Under the assumptions of Corollary \ref{convrate}, we have that for $\tau>0$ small enough 
\be
&&\sup\limits_{s\in[0,T]}\mE\sup _{s \leqslant t \leqslant s+\tau}|X_{t}^{\e,\g}-X_{s}^{\e,\g}|^{2}\leq C\tau, \label{xegtscr}\\
&&\sup\limits_{s\in[0,T]}\mE\sup _{s \leqslant t \leqslant s+\tau}|\hat{X}_{t}^{\e,\g}-\hat{X}_{s}^{\e,\g}|^{2}\leq C\tau,\label{hatxegtscr}\\
&&\sup\limits_{s\in[0,T]}\sup\limits_{s\leq t\leq s+\tau}|\bar{X}_{t}-\bar{X}_{s}|^2\leq C\tau, \label{barxtscr}
\ee
where the constant $C>0$ is independent of $\e, \d, \tau$.
\el
\begin{proof}
Since the proofs of (\ref{xegtscr}), (\ref{hatxegtscr}) and (\ref{barxtscr}) are similar, we only prove (\ref{xegtscr}).

First of all, note that for $0\leq s\leq t\leq s+\tau\leq T$, 
\ce
X_{t}^{\e,\g}-X_{s}^{\e,\g}=-K_t^{1,\e,\g}+K_s^{1,\e,\g}+\int_{s}^{t}b_{1}(r,X_{r}^{\e,\g},Y_{r}^{\e,\g})\dif r+\sqrt\e\int_{s}^{t}\s_{1}(r,X_{r}^{\e,\g},Y_{r}^{\e,\g})\dif W^1_{r}.
\de
Thus, the It\^o formula, (\ref{b1line}) and \cite[Lemma 2.4]{arrst} imply that
\ce
|X_{t}^{\e,\g}-X_{s}^{\e,\g}|^2&=&-2\int_s^t\<X_{r}^{\e,\g}-X_{s}^{\e,\g},\dif K_r^{1,\e,\g}\>\no\\
&&+2\int_s^t\<X_{r}^{\e,\g}-X_{s}^{\e,\g},b_{1}(r,X_{r}^{\e,\g},Y_{r}^{\e,\g})\>\dif r\no\\
&&+2\sqrt\e\int_s^t\<X_{r}^{\e,\g}-X_{s}^{\e,\g},\s_{1}(r,X_{r}^{\e,\g},Y_{r}^{\e,\g})\dif W^1_{r}\>\no\\
&&+\e\int_s^t\|\s_{1}(r,X_{r}^{\e,\g},Y_{r}^{\e,\g})\|^2\dif r\no\\
&\leq&\int_s^t|X_{r}^{\e,\g}-X_{s}^{\e,\g}|^2\dif r+C\int_s^t(1+|X_{r}^{\e,\g}|^2+|Y_{r}^{\e,\g}|^2)\dif r\\
&&+2\sqrt\e\left|\int_s^t\<X_{r}^{\e,\g}-X_{s}^{\e,\g},\s_{1}(r,X_{r}^{\e,\g},Y_{r}^{\e,\g})\dif W^1_{r}\>\right|.
\de
By the BDG inequality, it holds that
\ce
\mE\sup\limits_{t\in[s,s+\tau]}|X_{t}^{\e,\g}-X_{s}^{\e,\g}|^2&\leq& \int_s^{s+\tau}\mE|X_{r}^{\e,\g}-X_{s}^{\e,\g}|^2\dif r+C\int_s^{s+\tau}(1+\mE|X_{r}^{\e,\g}|^2+\mE|Y_{r}^{\e,\g}|^2)\dif r\\
&&+C\mE\(\int_s^{s+\tau}|X_{r}^{\e,\g}-X_{s}^{\e,\g}|^2\|\s_{1}(r,X_{r}^{\e,\g},Y_{r}^{\e,\g})\|^2\dif r\)^{1/2}\\
&\leq& \int_s^{s+\tau}\mE\sup\limits_{u\in[s,r]}|X_{u}^{\e,\g}-X_{s}^{\e,\g}|^2\dif r+C\tau+\frac{1}{2}\mE\sup\limits_{t\in[s,s+\tau]}|X_{t}^{\e,\g}-X_{s}^{\e,\g}|^2\\
&&+C\mE\int_s^{s+\tau}\|\s_{1}(r,X_{r}^{\e,\g},Y_{r}^{\e,\g})\|^2\dif r\\
&\leq& \int_s^{s+\tau}\mE\sup\limits_{u\in[s,r]}|X_{u}^{\e,\g}-X_{s}^{\e,\g}|^2\dif r+C\tau+\frac{1}{2}\mE\sup\limits_{t\in[s,s+\tau]}|X_{t}^{\e,\g}-X_{s}^{\e,\g}|^2,
\de
which together with the Gronwall inequality yields (\ref{xegtscr}). The proof is complete.
\end{proof} 

{\bf Proof of Corollary \ref{convrate}.}

Replacing (\ref{xegts}), (\ref{hatxegts}) and (\ref{barxts}) by (\ref{xegtscr}), (\ref{hatxegtscr}) and (\ref{barxtscr}), we can prove Corollary \ref{convrate} along the same lines as Theorem \ref{xbarxp}.

\section{Proof of Theorem \ref{ldpmmsde} and Corollary \ref{charate}}\label{prooseco}

In this section, we prove Theorem \ref{ldpmmsde} and Corollary \ref{charate}.

\subsection{Proof of Theorem \ref{ldpmmsde}}

In this subsection, we prove Theorem \ref{ldpmmsde}. 

By Theorem \ref{well}, we know that the system (\ref{Eq1ldp2}) has a unique strong solution $(X_{\cdot}^{\e,\g},K_{\cdot}^{1,\e,\g},\\Y_{\cdot}^{\e,\g},K_{\cdot}^{2,\e,\g})$. Thus, there exists a functional $\Psi^{\epsilon}: C([0,T];\mathbb{R}^{d_1+d_2})\mapsto C([0,T],\overline{\cD(A_1)})$ such that 
\ce
X^{\epsilon,\g}=\Psi^{\epsilon}(\sqrt{\e}W), \quad W:=(W^1,W^2).
\de
In order to prove the Laplace principle for $X^{\epsilon,\g}$, we will verify Condition \ref{cond} with $\mathbb{S}=C([0,T],\overline{\cD(A_1)})$.

First of all, we consider the following controlled processes:
\be\left\{\begin{array}{l}
\dif X_{t}^{\e,\g,u}\in -A_1(X_{t}^{\e,\g,u})\dif t+b_{1}(X_{t}^{\e,\g,u},Y_{t}^{\e,\g,u})\dif t+\sigma_1(X^{\epsilon,\g,u}_{t})\pi_1u(t)\dif t\\
\qquad\qquad+\sqrt{\e}\s_{1}(X_{t}^{\e,\g,u})\dif W^1_{t},\\
X_{0}^{\e,\g,u}=x_0\in\overline{\cD(A_1)},\quad  0\leq t\leq T,\\
\dif Y_{t}^{\e,\g,u}\in -A_2(Y_{t}^{\e,\g,u})\dif t+\frac{1}{\g}b_{2}(X_{t}^{\e,\g,u},Y_{t}^{\e,\g,u})\dif t+\frac{1}{\sqrt{\g \e}}\sigma_2(X^{\epsilon,\g,u}_{t},Y^{\epsilon,\g,u}_{t})\pi_2u(t)\dif t\\
\qquad\qquad +\frac{\sqrt{\e}}{\sqrt{\g\e}}\s_{2}(X_{t}^{\e,\g,u},Y_{t}^{\e,\g,u})\dif W^2_{t},  \\
Y_{0}^{\e,\g,u}=y_0\in\overline{\cD(A_2)},\quad  0\leq t\leq T,\quad u\in\mathbf{A}_2^{N},
\end{array}
\right.
\label{contproc}
\ee
where $\pi_1: \mR^{d_1+d_2}\mapsto \mR^{d_1}, \pi_2: \mR^{d_1+d_2}\mapsto \mR^{d_2}$ are two projection operators. Thus, by the Girsanov theorem, the system $(\ref{contproc})$ have a unique strong solution denoted by $(X^{\epsilon,\g,u}, K^{1,\epsilon,\g,u},Y^{\epsilon,\g,u}, K^{2,\epsilon,\g,u})$. Moreover, $X^{\epsilon,\g,u}=\Psi^{\epsilon}(\sqrt{\epsilon}W+\int_{0}^{\cdot}u(s)\dif s)$. 

Consider the following multivalued differential equation:
\be\left\{\begin{array}{l}
\dif\bar{X}^{u}_{t}\in -A_1(\bar{X}^{u}_{t})\dif t+\bar{b}_{1}(\bar{X}^u_{t})\dif t+\s_{1}(\bar{X}^u_{t})\pi_1u(t)\dif t, \quad u\in\mathbf{A}_2^{N},\\
\bar{X}^{u}_{0}=x_0\in\overline{\cD(A_1)}.
\end{array}
\right.
\label{deteequa}
\ee
By (\ref{barb1lip}) and $(\mathbf{H}^1_{b_{1}, \s_{1}})$, it holds that Eq.(\ref{deteequa}) has a unique solution $(\bar{X}^{u}, \bar{K}^{u})$. Define the measurable map $\Psi^{0}: C([0,T]; \mR^{d_1+d_2})\mapsto \mS$ by $\Psi^{0}(\int_{0}^{\cdot}u(s)\dif s)=\bar{X}^{u}$, and we verify Condition \ref{cond} through $\Psi^{\epsilon}, \Psi^{0}$. 

In the following, we divide the concrete deduction into two parts. In the first part (Subsubsection \ref{esti}), we present some key estimates for $X^{\e,\g,u}, Y^{\e,\g,u}, \bar{X}^{u}$. Then we justify Condition \ref{cond} in the second part (Subsubsection \ref{just}).

\subsubsection{Some estimates for $X^{\e,\g,u}, Y^{\e,\g,u}, \bar{X}^{u}$}\label{esti}

\bl \label{xutzutc}
Under the assumptions of Theorem \ref{ldpmmsde}, for $\{u_{\epsilon}, \e\in(0,1)\}\subset\mathbf{A}_{2}^{N}$, there exists a constant $C>0$ independent of $\e, \g$ such that 
\be
&&\mE\left(\sup\limits_{t\in[0,T]}|X_{t}^{\e,\g,u_\e}|^{2}\right)\leq C(1+|x_0|^{2}+|y_0|^{2}), \label{xeub}\\
&&\int_0^T\mE|Y_{r}^{\e,\g,u_\e}|^{2}\dif r\leq C(1+|x_0|^{2}+|y_0|^{2}), \label{yeub}\\
&&\mE|K^{1,\e,\g,u_\e}|_0^T\leq C(1+|x_0|^{2}+|y_0|^{2}). \label{keub}
\ee
\el
\begin{proof}
First of all, we estimate $X_{t}^{\e,\g,u_\e}$. Note that $X_{t}^{\e,\g,u_\e}$ satisfies the following equation:
\ce
X_{t}^{\e,\g,u_\e}&=&x_0-K_t^{1,\e,\g,u_\e}+\int_0^t b_{1}(X_{s}^{\e,\g,u_\e},Y_{s}^{\e,\g,u_\e})\dif s+\int_{0}^{t}\sigma_1(X^{\epsilon,\g,u_\e}_{s})\pi_1u_\e(s)\dif s\\
&&+\sqrt{\e}\int_0^t\s_{1}(X_{s}^{\e,\g,u_\e})\dif W^1_{s}.
\de
The It\^o formula yields that
\be
|X_{t}^{\e,\g,u_\e}|^2&=&|x_0|^2-2\int_0^t\<X_{s}^{\e,\g,u_\e}, \dif K_s^{1,\e,\g,u_\e}\>+2\int_0^t\<X_{s}^{\e,\g,u_\e}, b_{1}(X_{s}^{\e,\g,u_\e},Y_{s}^{\e,\g,u_\e})\>\dif s\no\\
&&+2\int_0^t\<X_{s}^{\e,\g,u_\e},\sigma_1(X^{\epsilon,\g,u_\e}_{s})\pi_1u_\e(s)\>\dif s+2\sqrt{\e}\int_0^t\<X_{s}^{\e,\g,u_\e},\s_{1}(X_{s}^{\e,\g,u_\e})\dif W^1_{s}\>\no\\
&&+\e\int_0^t\|\s_{1}(X_{s}^{\e,\g,u_\e})\|^2\dif s.
\label{itoxegu}
\ee
Then the BDG inequality, Lemma \ref{equi} and $(\ref{b1line})$ imply that for any $v'\in A_1(0)$
\ce
&&\mE\left(\sup\limits_{s\in[0,t]}|X_{s}^{\e,\g,u_\e}|^{2}\right)\\
&\leq&(|x_0|^2+2|v'|T)+(2|v'|+1)\mE\int_0^t|X_{r}^{\e,\g,u_\e}|^2\dif r+\bar{L}_{b_{1}, \s_{1}}\mE\int_0^t(1+|X_{r}^{\e,\g,u_\e}|^2+|Y_{r}^{\e,\g,u_\e}|^2)\dif r\no\\
&&+2\mE\sup\limits_{s\in[0,t]}\left|\int_0^s\<X_{r}^{\e,\g,u_\e},\sigma_1(X^{\epsilon,\g,u_\e}_{r})\pi_1u_\e(r)\>\dif r\right|+2\mE\sup\limits_{s\in[0,t]}\left|\int_0^s\<X_{r}^{\e,\g,u_\e},\s_{1}(X_{r}^{\e,\g,u_\e})\dif W^1_{r}\>\right|\no\\ 
&\leq&C(|x_0|^{2}+1)+C\mE\int_{0}^{t}|X_{r}^{\e,\g,u_\e}|^2\dif r+C\mE\int_{0}^{t}|Y_{r}^{\e,\g,u_\e}|^{2}\dif r+\frac{1}{4}\mE\left(\sup\limits_{s\in[0,t]}|X_{s}^{\e,\g,u_\e}|^{2}\right)\no\\
&&+C\mE\left(\int_0^t\|\sigma_1(X^{\epsilon,\g,u_\e}_{r})\|^2\dif r\right)\left(\int_0^t|u_\e(r)|^2\dif r\right)+C\mE\left(\int_0^t|X_{r}^{\e,\g,u_\e}|^2\|\s_{1}(X_{r}^{\e,\g,u_\e})\|^2\dif r\right)^{1/2}\no\\
&\leq&C(|x_0|^{2}+1)+C\mE\int_{0}^{t}|X_{r}^{\e,\g,u_\e}|^2\dif r+C\mE\int_{0}^{t}|Y_{r}^{\e,\g,u_\e}|^{2}\dif r\no\\
&&+\frac{1}{2}\mE\left(\sup\limits_{s\in[0,t]}|X_{s}^{\e,\g,u_\e}|^{2}\right)+C\mE\int_0^t\|\s_{1}(X_{r}^{\e,\g,u_\e})\|^2\dif r,
\de
and furthermore
\be
\mE\left(\sup\limits_{s\in[0,t]}|X_{s}^{\e,\g,u_\e}|^{2}\right)\leq C(|x_0|^{2}+1)+C\int_{0}^{t}\mE|X_{r}^{\e,\g,u_\e}|^{2}\dif r+C\int_{0}^{t}\mE|Y_{r}^{\e,\g,u_\e}|^{2}\dif r.
\label{exqcu}
\ee

For $Y_{t}^{\e,\g,u_\e}$, fix $v\in A_2(0)$. Applying the It\^{o} formula to $|Y_{t}^{\e,\g,u_\e}|^{2}e^{\l t}$ for any $\l>0$ and taking the expectation, one could obtain that
\ce
\mE|Y_{t}^{\e,\g,u_\e}|^{2}e^{\l t}
&=& |y_0|^{2}+\l\mE\int_0^t|Y_{s}^{\e,\g,u_\e}|^{2}e^{\l s}\dif s-2\mE\int_0^te^{\l s}\<Y_{s}^{\e,\g,u_\e},\dif K_s^{2,\e,\g,u_\e}\>\\
&&+\frac{2}{\g}\mE\int_{0}^{t}e^{\l s}\<Y_{s}^{\e,\g,u_\e}, b_{2}(X_{s}^{\e,\g,u_\e},Y_{s}^{\e,\g,u_\e})\>\dif s\\
&&+\frac{1}{\g}\mE\int_{0}^{t}e^{\l s}\|\s_{2}(X_{s}^{\e,\g,u_\e},Y_{s}^{\e,\g,u_\e})\|^2\dif s\\
&&+\frac{2}{\sqrt{\g \e}}\mE\int_{0}^{t}e^{\l s}\<Y_{s}^{\e,\g,u_\e},\sigma_2(X^{\epsilon,\g,u_\e}_{s},Y^{\epsilon,\g,u_\e}_{s})\pi_2u_\e(s)\>\dif s\\
&\leq& |y_0|^{2}+\l\mE\int_0^t|Y_{s}^{\e,\g,u_\e}|^{2}e^{\l s}\dif s+2\mE\int_0^te^{\l s}|v||Y_{s}^{\e,\g,u_\e}|\dif s\\
&&+\frac{1}{\g}\mE\int_{0}^{t}e^{\l s}\(-\a|Y_{s}^{\e,\g,u_\e}|^{2}+C(1+|X_{s}^{\e,\g,u_\e}|^{2})\)\dif s\\
&&+\frac{2}{\sqrt{\g \e}}\mE\int_{0}^{t}e^{\l s}|Y_{s}^{\e,\g,u_\e}|\|\sigma_2(X^{\epsilon,\g,u_\e}_{s},Y^{\epsilon,\g,u_\e}_{s})\||u_\e(s)|\dif s,
\de
where we use (\ref{bemu}) in the first inequality. Note that for the third and last terms of the right side for the above inequality
\ce
2|v||Y_{s}^{\e,\g,u_\e}|\leq \frac{\a}{3\g}|Y_{s}^{\e,\g,u_\e}|^2+C|v|^2,
\de
and
\ce
&&\frac{2}{\sqrt{\g \e}}|Y_{s}^{\e,\g,u_\e}|\|\sigma_2(X^{\epsilon,\g,u_\e}_{s},Y^{\epsilon,\g,u_\e}_{s})\||u_\e(s)|\\
&\leq&\frac{\a}{3\g}|Y_{s}^{\e,\g,u_\e}|^{2}+\frac{C}{\e}\|\sigma_2(X^{\epsilon,\g,u_\e}_{s},Y^{\epsilon,\g,u_\e}_{s})\|^2|u_\e(s)|^2\\
&\leq&\frac{\a}{3\g}|Y_{s}^{\e,\g,u_\e}|^{2}+\frac{C}{\e}|u_\e(s)|^2,
\de
where $(\mathbf{H}^3_{\s_{2}})$ is used. Thus, we have that
\ce
\mE|Y_{t}^{\e,\g,u_\e}|^{2}e^{\l t}&\leq& |y_0|^{2}+\(\l+\frac{2\a}{3\g}-\frac{\a}{\g}\)\mE\int_0^t|Y_{s}^{\e,\g,u_\e}|^{2}e^{\l s}\dif s+C\int_0^te^{\l s}|v|^2\dif s\\
&&+\frac{C}{\g}\mE\int_{0}^{t}e^{\l s}(1+|X_{s}^{\e,\g,u_\e}|^{2})\dif s+\frac{C}{\e}\mE\int_{0}^{t}e^{\l s}|u_\e(s)|^2\dif s\\
&\leq& |y_0|^{2}+\(\l+\frac{2\a}{3\g}-\frac{\a}{\g}\)\mE\int_0^t|Y_{s}^{\e,\g,u_\e}|^{2}e^{\l s}\dif s+C|v|^2\frac{e^{\l t}-1}{\l}\\
&&+C\frac{e^{\l t}-1}{\g\l}\left(\mE\left(\sup\limits_{s\in[0,t]}|X_{s}^{\e,\g,u_\e}|^{2}\right)+1\right)+\frac{C}{\e}\mE\int_{0}^{t}e^{\l s}|u_\e(s)|^2\dif s.
\de
From this and taking $\l=\frac{\a}{3\g}$, it follows that
\ce
\mE|Y_{t}^{\e,\g,u_\e}|^{2}\leq C(|y_0|^{2}+1)+C\mE\left(\sup\limits_{s\in[0,t]}|X_{s}^{\e,\g,u_\e}|^{2}\right)+\frac{C}{\e}\mE\int_{0}^{t}e^{-\frac{\a}{3\g}(t-s)}|u_\e(s)|^2\dif s,
\de
and furthermore
\be
\int_0^t\mE|Y_{r}^{\e,\g,u_\e}|^{2}\dif r&\leq& CT(|y_0|^{2}+1)+C\int_0^t\mE\left(\sup\limits_{s\in[0,r]}|X_{s}^{\e,\g,u_\e}|^{2}\right)\dif r\no\\
&&+\frac{C}{\e}\mE\int_0^t\int_{0}^{r}e^{-\frac{\a}{3\g}(r-s)}|u_\e(s)|^2\dif s\dif r\no\\
&\leq& CT(|y_0|^{2}+1)+C\int_0^t\mE\left(\sup\limits_{s\in[0,r]}|X_{s}^{\e,\g,u_\e}|^{2}\right)\dif r\no\\
&&+C\left(\frac{\g}{\e}\right)\mE\int_0^t|u_\e(s)|^2\dif s\no\\
&\leq& CT(|y_0|^{2}+1)+C\int_0^t\mE\left(\sup\limits_{s\in[0,r]}|X_{s}^{\e,\g,u_\e}|^{2}\right)\dif r.
\label{zeues}
\ee

Inserting (\ref{zeues}) in (\ref{exqcu}), by the Gronwall inequality one can get (\ref{xeub}) and (\ref{yeub}).

Finally, for $K^{1,\e,\g,u_\e}$, by (\ref{itoxegu}) and Lemma \ref{inteineq}, it holds that
\ce
|X_{T}^{\e,\g,u_\e}|^2&=&|x_0|^2-2\int_0^T\<X_{s}^{\e,\g,u_\e}, \dif K_s^{1,\e,\g,u_\e}\>+2\int_0^T\<X_{s}^{\e,\g,u_\e}, b_{1}(X_{s}^{\e,\g,u_\e},Y_{s}^{\e,\g,u_\e})\>\dif s\no\\
&&+2\int_0^T\<X_{s}^{\e,\g,u_\e},\sigma_1(X^{\epsilon,\g,u_\e}_{s})\pi_1u_\e(s)\>\dif s+2\sqrt{\e}\int_0^T\<X_{s}^{\e,\g,u_\e},\s_{1}(X_{s}^{\e,\g,u_\e})\dif W^1_{s}\>\\
&&+\e\int_0^T\|\s_{1}(X_{s}^{\e,\g,u_\e})\|^2\dif s\no\\
&\leq&|x_0|^2-2M_1\left| K^{1,\e,\g,u_\e} \right|_{0}^{T}+2M _2\int_0^T{\left| X^{\e,\g,u_\e}_s\right|}\dif s+2M_3T+\int_0^T|X_{s}^{\e,\g,u_\e}|^2\dif s\\
&&+\int_0^T|b_{1}(X_{s}^{\e,\g,u_\e},Y_{s}^{\e,\g,u_\e})|^2\dif s+\left(\sup\limits_{s\in[0,T]}|X_{s}^{\e,\g,u_\e}|^2\right)\\
&&+\int_0^T\|\sigma_1(X^{\epsilon,\g,u_\e}_{s})\|^2\dif s\int_0^T|u_\e(s)|^2\dif s+2\sqrt{\e}\int_0^T\<X_{s}^{\e,\g,u_\e},\s_{1}(X_{s}^{\e,\g,u_\e})\dif W^1_{s}\>\\
&&+\e\int_0^T\|\s_{1}(X_{s}^{\e,\g,u_\e})\|^2\dif s,
\de
and 
\ce
2M_1\mE\left| K^{1,\e,\g,u_\e} \right|_{0}^{T}&\leq& |x_0|^2+2(M_2+M_3)T+(2M _2+1)\int_0^T\mE{\left| X^{\e,\g,u_\e}_s\right|^2}\dif s\\
&&+\mE\left(\sup\limits_{s\in[0,T]}|X_{s}^{\e,\g,u_\e}|^2\right)+C\int_0^T(1+\mE{\left| X^{\e,\g,u_\e}_s\right|^2}+\mE{\left| Y^{\e,\g,u_\e}_s\right|^2})\dif s,
\de
which together with (\ref{xeub}) and (\ref{yeub}) yields (\ref{keub}). The proof is complete.
\end{proof}

\bl
Under the assumptions of Theorem \ref{ldpmmsde}, for $\{u_{\epsilon}, \e\in(0,1)\}\subset\mathbf{A}_{2}^{N}$, we have that
\be
\lim\limits_{l\rightarrow 0}\sup\limits_{s\in[0,T]}\mE\sup _{s \leqslant t \leqslant s+l}|X_{t}^{\e,\g,u_\e}-X_{s}^{\e,\g,u_\e}|^{2}=0. 
\label{xegutse}
\ee
\el
\begin{proof}
Note that for $0\leq s<t\leq T$,
\ce
X_{t}^{\e,\g,u_\e}-X_{s}^{\e,\g,u_\e}&=&-K_t^{1,\e,\g,u_\e}+K_s^{1,\e,\g,u_\e}+\int_{s}^{t}b_{1}(X_{r}^{\e,\g,u_\e},Y_{r}^{\e,\g,u_\e})\dif r\\
&&+\int_{s}^{t}\sigma_1(X^{\epsilon,\g,u_\e}_{r})\pi_1u_\e(r)\dif r+\int_{s}^{t}\s_{1}(X_{r}^{\e,\g,u_\e})\dif W^1_{r}.
\de
Thus, the It\^o formula and (\ref{b1line}) imply that
\be
&&|X_{t}^{\e,\g,u_\e}-X_{s}^{\e,\g,u_\e}|^2\no\\
&=&-2\int_s^t\<X_{r}^{\e,\g,u_\e}-X_{s}^{\e,\g,u_\e},\dif K_r^{1,\e,\g,u_\e}\>+2\int_s^t\<X_{r}^{\e,\g,u_\e}-X_{s}^{\e,\g,u_\e},b_{1}(X_{r}^{\e,\g,u_\e},Y_{r}^{\e,\g,u_\e})\>\dif r\no\\
&&+2\int_s^t\<X_{r}^{\e,\g,u_\e}-X_{s}^{\e,\g,u_\e},\sigma_1(X^{\epsilon,\g,u_\e}_{r})\pi_1u_\e(r)\>\dif r+2\int_s^t\<X_{r}^{\e,\g,u_\e}-X_{s}^{\e,\g,u_\e},\s_{1}(X_{r}^{\e,\g,u_\e})\dif W^1_{r}\>\no\\
&&+\int_s^t\|\s_{1}(X_{r}^{\e,\g,u_\e})\|^2\dif r\no\\
&\leq&-2\int_s^t\<X_{r}^{\e,\g,u_\e}-X_{s}^{\e,\g,u_\e},\dif K_r^{1,\e,\g,u_\e}\>+\int_s^t|X_{r}^{\e,\g,u_\e}-X_{s}^{\e,\g,u_\e}|^2\dif r\no\\
&&+\bar{L}_{b_{1}, \s_{1}}\int_s^t(1+|X_{r}^{\e,\g,u_\e}|^2+|Y_{r}^{\e,\g,u_\e}|^2)\dif r+2\int_s^t\<X_{r}^{\e,\g,u_\e}-X_{s}^{\e,\g,u_\e},\sigma_1(X^{\epsilon,\g,u_\e}_{r})\pi_1u_\e(r)\>\dif r\no\\
&&+2\int_s^t\<X_{r}^{\e,\g,u_\e}-X_{s}^{\e,\g,u_\e},\s_{1}(X_{r}^{\e,\g,u_\e})\dif W^1_{r}\>.
\label{xegustito}
\ee

Next, we compute $-2\int_s^t\<X_{r}^{\e,\g,u_\e}-X_{s}^{\e,\g,u_\e},\dif K_r^{1,\e,\g,u_\e}\>$. Let $X_s^{\varepsilon, \g, u_\e, l, R}$ be the projection of $X_s^{\e,\g,u_\e}$ on $\{x \in B(0,R): d(x,(\overline{\cD(A_1)})^c) \geqslant l+h_R(l)\}$ (See the deduction in Lemma \ref{xehatxets} for $l, h_R$). Thus, for $Z_s^{\varepsilon,\g,u_\e, l, R} \in A_1(X_s^{\varepsilon,\g,u_\e, l, R}), \sup\limits_{t\in[0,T]}|X_t^{\varepsilon,\g,u_\e}| \leqslant R$ and $0<t-s<l$, it holds that
\ce
&&-2 \int_s^t\left\langle X_r^{\e,\g,u_\e}-X_s^{\e,\g,u_\e}, \dif K_r^{1,\e,\g,u_\e} \right\rangle\\
 &=&-2 \int_s^t\left\langle X_r^{\varepsilon,\g,u_\e}-X_s^{\varepsilon,\g,u_\e, l, R}, \dif K_r^{1,\e,\g,u_\e} \right\rangle-2 \int_s^t\left\langle X_s^{\varepsilon, \g,u_\e, l, R}-X_s^{\varepsilon,\g,u_\e}, \dif K_r^{1,\e,\g,u_\e} \right\rangle \\ 
&\leqslant& -2 \int_s^t\left\langle X_r^{\varepsilon,\g,u_\e}-X_s^{\varepsilon,\g,u_\e, l, R}, Z_s^{\varepsilon,\g,u_\e, l, R}\right\rangle \dif r+2\left(l+h_R(l)\right)\left|K^{1,\e,\g,u_\e} \right|_0^T\\
&\leqslant& 4 l^{1 / 2} R+2\left(l+h_R(l)\right)\left|K^{1,\e,\g,u_\e} \right|_0^T,
\de
and furthermore by (\ref{xegustito})
\ce
&&\sup _{s \leqslant t \leqslant s+l}\left|X_t^{\varepsilon,\g,u_\e}-X_s^{\varepsilon,\g,u_\e}\right|^2 I_{\{\sup\limits_{t\in[0,T]}|X_t^{\varepsilon,\g,u_\e}| \leqslant R\}} \\
&\leqslant&\left(4 l^{1 / 2} R+2\left(l+h_R(l)\right)\left|K^{1,\e,\g,u_\e} \right|_0^T\right)+\sup _{s \leqslant t \leqslant s+l}\int_s^t|X_{r}^{\e,\g,u_\e}-X_{s}^{\e,\g,u_\e}|^2\dif r\\
&&+\bar{L}_{b_{1}, \s_{1}}\sup _{s \leqslant t \leqslant s+l}\int_s^t\left(1+\left|X_r^{\varepsilon,\g,u_\e}\right|^2+\left|Y_r^{\varepsilon,\g,u_\e}\right|^2\right) \dif r\\
&&+2\sup _{s \leqslant t \leqslant s+l}\left|\int_s^t\<X_{r}^{\e,\g,u_\e}-X_{s}^{\e,\g,u_\e},\sigma_1(X^{\epsilon,\g,u_\e}_{r})\pi_1u_\e(r)\>\dif r\right|I_{\{\sup\limits_{t\in[0,T]}|X_t^{\varepsilon,\g,u_\e}| \leqslant R\}}\no\\
&&+C \sup _{s \leqslant t \leqslant s+l}\left|\int_s^t\left\langle X_r^{\varepsilon,\g,u_\e}-X_s^{\varepsilon,\g,u_\e}, \sigma_1\left(X_r^{\varepsilon,\g,u_\e}\right)\right\rangle \dif W^1_r\right|I_{\{\sup\limits_{t\in[0,T]}|X_t^{\varepsilon,\g,u_\e}| \leqslant R\}}.
\de

From the above deduction and the BDG inequality, it follows that
\ce
&&\mE\sup _{s \leqslant t \leqslant s+l}\left|X_t^{\varepsilon,\g,u_\e}-X_s^{\varepsilon,\g,u_\e}\right|^2 I_{\{\sup\limits_{t\in[0,T]}|X_t^{\varepsilon,\g,u_\e}| \leqslant R\}}\\
&\leq&\left(4 l^{1 / 2} R+2\left(l+h_R(l)\right)\mE|K^{1,\e,\g,u_\e}|_0^T\right)+C(1+|x_0|^{2}+|y_0|^{2})l\\
&&+2\mE\int_s^{s+l}|X_{r}^{\e,\g,u_\e}-X_{s}^{\e,\g,u_\e}|\|\sigma_1(X^{\epsilon,\g,u_\e}_{r})\||u_\e(r)|\dif rI_{\{\sup\limits_{t\in[0,T]}|X_t^{\varepsilon,\g,u_\e}| \leqslant R\}}\no\\
&&+C\mE\left(\int_s^{s+l}|X_r^{\varepsilon,\g,u_\e}-X_s^{\varepsilon,\g,u_\e}|^2\|\sigma_1(X_r^{\varepsilon,\g,u_\e})\|^2 \dif r\right)^{1/2}I_{\{\sup\limits_{t\in[0,T]}|X_t^{\varepsilon,\g,u_\e}| \leqslant R\}}\\
&\leq&\left(4 l^{1 / 2} R+2\left(l+h_R(l)\right)\mE|K^{1,\e,\g,u_\e}|_0^T\right)+C(1+|x_0|^{2}+|y_0|^{2})l\\
&&+\frac{1}{2}\mE\sup _{s \leqslant t \leqslant s+l}\left|X_t^{\varepsilon,\g,u_\e}-X_s^{\varepsilon,\g,u_\e}\right|^2 I_{\{\sup\limits_{t\in[0,T]}|X_t^{\varepsilon,\g,u_\e}| \leqslant R\}}+C(1+|x_0|^{2}+|y_0|^{2})l,
\de
which yields that
\ce
&&\mE\sup _{s \leqslant t \leqslant s+l}\left|X_t^{\varepsilon,\g,u_\e}-X_s^{\varepsilon,\g,u_\e}\right|^2 I_{\{\sup\limits_{t\in[0,T]}|X_t^{\varepsilon,\g,u_\e}| \leqslant R\}}\\
&\leq& \left(4 l^{1 / 2} R+2\left(l+h_R(l)\right)\mE|K^{1,\e,\g,u_\e}|_0^T\right)+C(1+|x_0|^{2}+|y_0|^{2})l.
\de
Based on this, letting $l\rightarrow 0$ and $R\rightarrow \infty$, one could conclude (\ref{xegutse}).
\end{proof}

\bl
Under the assumptions of Theorem \ref{ldpmmsde}, for any $u\in\mathbf{A}_{2}^{N}$, it holds that
\be
&&\sup\limits_{t\in[0,T]}|\bar{X}^u_{t}|^{2}\leq C(1+|x_0|^{2}), a.s.\label{barxub}\\
&&|\bar{K}^{u}|_0^T\leq C(1+|x_0|^{2}), a.s. \label{barkub}\\
&&\lim\limits_{l\rightarrow 0}\sup\limits_{s\in[0,T]}\sup\limits_{s\leq t\leq s+l}|\bar{X}^u_{t}-\bar{X}^u_{s}|^2=0, a.s.,\label{barxuts}
\ee
where the constant $C>0$ is nonrandom.
\el

Since the proofs of (\ref{barxub}), (\ref{barkub}) and (\ref{barxuts}) are similar to that of (\ref{xeub}), (\ref{keub}) and (\ref{xegutse}), respectively, we omit them.

Next, set for $u_\e, u\in\mathbf{A}_{2}^{N}$
\ce
g_\e(t):=\int_0^t\sigma_1(\bar{X}^{u}_{r})\pi_1(u_{\epsilon}(r)-u(r))\dif r,
\de
and we have the following result.
\bl\label{ge}
suppose that the assumptions of Theorem \ref{ldpmmsde} hold, and $u_{\epsilon}$ converges to $u$ almost surely as $\epsilon\rightarrow 0$. Then $g_\e(\cdot)$ tends to $0$ almost surely in $C([0,T],\mR^n)$ as $\epsilon\rightarrow 0$. 
\el
\begin{proof}
On one hand, $\{g_\e; \e\in (0,1)\}$ is relatively compact in $C([0,T],\mR^n)$. In fact, by the Ascoli-Arzel\'a lemma, 
we only need to justify that 

$(i)$ $\sup\limits_{\e\in (0,1)}\sup\limits_{t\in[0,T]}\left|g_\e(t)\right|<\infty$,

$(ii)$ $\{[0,T]\ni t\mapsto g_\e(t); \e\in (0,1)\}$ is equi-continuous.

For $0\leq s<t\leq T$, it holds that
\ce
|g_\e(s)-g_\e(t)|&\leq&\left|\int_{s}^{t}\s_1(\bar{X}^{u}_{r})\pi_1(u_\e(r)-u(r))\dif r\right|\leq 2N^{1/2}\left(\int_{s}^{t}\|\s_1(\bar{X}^{u}_{r})\|^2\dif r\right)^{1/2}\\
&\leq&2N^{1/2}C\left(\int_{s}^{t}(1+|\bar{X}^{u}_{r}|^2)\dif r\right)^{1/2},
\de
where (\ref{b1line}) is used. Letting $s=0$, by (\ref{barxub})  we have that
\ce
|g_\e(t)|\leq 2N^{1/2}C(1+|x_0|),
\de
where $C$ is independent of $\e$. So, $(i)$ holds.

For $(ii)$, noticing that 
\ce
|g_\e(s)-g_\e(t)|\leq 2N^{1/2}C(1+|x_0|)(t-s)^{1/2},
\de
we know that $(ii)$ holds.

On the other hand, note that 
\ce
\int_{0}^{t}\|\s_1(\bar{X}^{u}_{r})\|^2\dif r\leq C\int_{0}^{t}(1+|\bar{X}^{u}_{r}|^2)\dif r<\infty. 
\de
Since $u_{\e}\rightarrow u$ almost surely in $\mathbf{A}_2^{N}$ as ${\e}\rightarrow0$, one get that for any $t\in[0,T]$
$$
\lim\limits_{{\e}\rightarrow0}g_\e(t)=0,
$$
which implies that
\ce
\lim\limits_{{\e}\rightarrow0}\sup\limits_{t\in[0,T]}\left|g_\e(t)\right|=0.
\de
The proof is complete.
\end{proof}

Finally we introduce the following auxiliary process: 
\be\left\{\begin{array}{l}
\dif\hat{Y}_{t}^{\e,\g,u_\e}\in -A_2(\hat{Y}_{t}^{\e,\g,u_\e})\dif t+\frac{1}{\g}b_{2}(X_{t(\d)}^{\e,\g,u_\e},\hat{Y}_{t}^{\e,\g,u_\e})\dif t+\frac{1}{\sqrt{\g}}\s_{2}(X_{t(\d)}^{\e,\g,u_\e},\hat{Y}_{t}^{\e,\g,u_\e})\dif W^2_{t}, \\
\hat{Y}_{0}^{\e,\g,u_\e}=y_0.
\end{array}
\right.
\label{hatzu}
\ee
Under $(\mathbf{H}_{A_{2}})$ and $(\mathbf{H}^1_{b_{2}, \s_{2}})$, the above multivalued SDE has a unique solution $(\hat{Y}_{t}^{\e,\g,u_\e}, \hat{K}_{t}^{\e,\g,u_\e})$ (\cite{rwzx}). 

\bl
Under the assumptions of Theorem \ref{ldpmmsde}, for $\{u_{\epsilon}, \e\in(0,1)\}\subset\mathbf{A}_{2}^{N}$, there exists a constant $C>0$ independent of $\e, \g$ such that 
 \be
 &&\sup\limits_{t\in[0,T]}\mE|\hat{Y}_{t}^{\e,\g,u_\e}|^{2}\leq C(1+|x_0|^{2}+|y_0|^{2}), \label{hatzub}\\
&&\int_0^T\mE|Y_{t}^{\e,\g,u_\e}-\hat{Y}_{t}^{\e,\g,u_\e}|^{2}\dif t\leq C\frac{\g}{\e}+\frac{CT}{\a}\left(\sup\limits_{s\in[0,T]}\mE\sup _{s \leqslant r \leqslant s+\d}|X_{r}^{\e,\g,u_\e}-X_{s}^{\e,\g,u_\e}|^{2}\right).
\label{unztu}
\ee
\el
\begin{proof}
Since the proof of (\ref{hatzub}) is similar to that for (\ref{zeb}), we only prove (\ref{unztu}).

First of all, by (\ref{contproc}) and (\ref{hatzu}), we have that 
\ce
Y_{t}^{\e,\g,u_\e}-\hat{Y}_{t}^{\e,\g,u_\e}
&=&-K_t^{2,\e,\g,u_\e}+\hat{K}_t^{2,\e,\g,u_\e}+\frac{1}{\g}\int_{0}^{t}\(b_{2}(X_{s}^{\e,\g,u_\e},Y_{s}^{\e,\g,u_\e})
-b_{2}(X_{s(\d)}^{\e,\g,u_\e},\hat{Y}_{s}^{\e,\g,u_\e})\)\dif s\\
&&+\frac{1}{\sqrt{\g}}\int_{0}^{t}\(\s_{2}(X_{s}^{\e,\g,u_\e},Y_{s}^{\e,\g,u_\e})
-\s_{2}(X_{s(\d)}^{\e,\g,u_\e},\hat{Y}_{s}^{\e,\g,u_\e})\)\dif W^2_{s}\\
&&+\frac{1}{\sqrt{\g \e}}\int_{0}^{t}\s_{2}(X_{s}^{\e,\g,u_\e},Y_{s}^{\e,\g,u_\e})\pi_2u_\e(s)\dif s.
\de
Applying the It\^{o} formula to $|Y_{t}^{\e,\g,u_\e}-\hat{Y}_{t}^{\e,\g,u_\e}|^{2}e^{\l t}$ for any $\l>0$ and taking the expectation, by Lemma \ref{equi} one could obtain that
\ce
&&\mE|Y_{t}^{\e,\g,u_\e}-\hat{Y}_{t}^{\e,\g,u_\e}|^{2}e^{\l t}\\
&=&\l\mE\int_{0}^{t}|Y_{s}^{\e,\g,u_\e}-\hat{Y}_{s}^{\e,\g,u_\e}|^{2}e^{\l s}\dif s-\mE\int_{0}^{t}2e^{\l s}\<Y_{s}^{\e,\g,u_\e}-\hat{Y}_{s}^{\e,\g,u_\e},\dif (K_s^{2,\e,\g,u_\e}-\hat{K}_s^{2,\e,\g,u_\e})\>\\
&&+\frac{1}{\g}\mE\int_{0}^{t}2e^{\l s}\<Y_{s}^{\e,\g,u_\e}-\hat{Y}_{s}^{\e,\g,u_\e}, b_{2}(X_{s}^{\e,\g,u_\e},Y_{s}^{\e,\g,u_\e})
-b_{2}(X_{s(\d)}^{\e,\g,u_\e},\hat{Y}_{s}^{\e,\g,u_\e})\>\dif s\\
&&+\frac{1}{\sqrt{\g \e}}\mE\int_{0}^{t}2e^{\l s}\<Y_{s}^{\e,\g,u_\e}-\hat{Y}_{s}^{\e,\g,u_\e}, \s_{2}(X_{s}^{\e,\g,u_\e},Y_{s}^{\e,\g,u_\e})\pi_2u_\e(s)\>\dif s\\
&&+\frac{1}{\g}\mE\int_{0}^{t}e^{\l s}\|\s_{2}(X_{s}^{\e,\g,u_\e},Y_{s}^{\e,\g,u_\e})
-\s_{2}(X_{s(\d)}^{\e,\g,u_\e},\hat{Y}_{s}^{\e,\g,u_\e})\|^{2}\dif s\\
&\leq&\l\mE\int_{0}^{t}|Y_{s}^{\e,\g,u_\e}-\hat{Y}_{s}^{\e,\g,u_\e}|^{2}e^{\l s}\dif s+\frac{1}{\g}\mE\int_{0}^{t}2e^{\l s}\<Y_{s}^{\e,\g,u_\e}-\hat{Y}_{s}^{\e,\g,u_\e}, b_{2}(X_{s}^{\e,\g,u_\e},Y_{s}^{\e,\g,u_\e})\\
&&-b_{2}(X_{s}^{\e,\g,u_\e},\hat{Y}_{s}^{\e,\g,u_\e})\>\dif s+\frac{1}{\g}\mE\int_{0}^{t}e^{\l s}\|\s_{2}(X_{s}^{\e,\g,u_\e},Y_{s}^{\e,\g,u_\e})-\s_{2}(X_{s}^{\e,\g,u_\e},\hat{Y}_{s}^{\e,\g,u_\e})\|^{2}\dif s\\
&&+\frac{1}{\sqrt{\g \e}}\mE\int_{0}^{t}2e^{\l s}\<Y_{s}^{\e,\g,u_\e}-\hat{Y}_{s}^{\e,\g,u_\e}, \s_{2}(X_{s}^{\e,\g,u_\e},Y_{s}^{\e,\g,u_\e})\pi_2u_\e(s)\>\dif s\\
&&+\frac{1}{\g}\mE\int_{0}^{t}2e^{\l s}\<Y_{s}^{\e,\g,u_\e}-\hat{Y}_{s}^{\e,\g,u_\e}, b_{2}(X_{s}^{\e,\g,u_\e},\hat{Y}_{s}^{\e,\g,u_\e})-b_{2}(X_{s(\d)}^{\e,\g,u_\e},\hat{Y}_{s}^{\e,\g,u_\e})\>\dif s\\
&&+\frac{1}{\g}\mE\int_{0}^{t}e^{\l s}\|\s_{2}(X_{s}^{\e,\g,u_\e},Y_{s}^{\e,\g,u_\e})-\s_{2}(X_{s}^{\e,\g,u_\e},\hat{Y}_{s}^{\e,\g,u_\e})\|^{2}\dif s\\
&&+\frac{1}{\g}\mE\int_{0}^{t}2e^{\l s}\|\s_{2}(X_{s}^{\e,\g,u_\e},\hat{Y}_{s}^{\e,\g,u_\e})-\s_{2}(X_{s(\d)}^{\e,\g,u_\e},\hat{Y}_{s}^{\e,\g,u_\e})\|^{2}\dif s.
\de
By $(\mathbf{H}^{2}_{b_{2},\s_{2}})$, it holds that
\ce
&&\frac{1}{\g}\mE\int_{0}^{t}2e^{\l s}\<Y_{s}^{\e,\g,u_\e}-\hat{Y}_{s}^{\e,\g,u_\e}, b_{2}(X_{s}^{\e,\g,u_\e},Y_{s}^{\e,\g,u_\e})-b_{2}(X_{s}^{\e,\g,u_\e},\hat{Y}_{s}^{\e,\g,u_\e})\>\dif s\\
&&+\frac{1}{\g}\mE\int_{0}^{t}e^{\l s}\|\s_{2}(X_{s}^{\e,\g,u_\e},Y_{s}^{\e,\g,u_\e})-\s_{2}(X_{s}^{\e,\g,u_\e},\hat{Y}_{s}^{\e,\g,u_\e})\|^{2}\dif s\\
&\leq&-\frac{\b}{\g}\mE\int_{0}^{t}e^{\l s}|Y_{s}^{\e,\g,u_\e}-\hat{Y}_{s}^{\e,\g,u_\e}|^{2}\dif s.
\de
$(\mathbf{H}^{3}_{\s_{2}})$ implies that
\ce
&&\frac{1}{\sqrt{\g \e}}\mE\int_{0}^{t}2e^{\l s}\<Y_{s}^{\e,\g,u_\e}-\hat{Y}_{s}^{\e,\g,u_\e}, \s_{2}(X_{s}^{\e,\g,u_\e},Y_{s}^{\e,\g,u_\e})\pi_2u_\e(s)\>\dif s\\
&\leq&\frac{\a}{2\g}\mE\int_{0}^{t}e^{\l s}|Y_{s}^{\e,\g,u_\e}-\hat{Y}_{s}^{\e,\g,u_\e}|^{2}\dif s+\frac{C}{\e}\mE\int_{0}^{t}e^{\l s}|u_\e(s)|^2\dif s.
\de
And by $(\mathbf{H}^{1}_{b_{2},\s_{2}})$ we know that
\ce
&&\frac{1}{\g}\mE\int_{0}^{t}2e^{\l s}\<Y_{s}^{\e,\g,u_\e}-\hat{Y}_{s}^{\e,\g,u_\e}, b_{2}(X_{s}^{\e,\g,u_\e},\hat{Y}_{s}^{\e,\g,u_\e})-b_{2}(X_{s(\d)}^{\e,\g,u_\e},\hat{Y}_{s}^{\e,\g,u_\e})\>\dif s\\
&&+\frac{1}{\g}\mE\int_{0}^{t}e^{\l s}\|\s_{2}(X_{s}^{\e,\g,u_\e},Y_{s}^{\e,\g,u_\e})-\s_{2}(X_{s}^{\e,\g,u_\e},\hat{Y}_{s}^{\e,\g,u_\e})\|^{2}\dif s\\
&&+\frac{1}{\g}\mE\int_{0}^{t}2e^{\l s}\|\s_{2}(X_{s}^{\e,\g,u_\e},\hat{Y}_{s}^{\e,\g,u_\e})-\s_{2}(X_{s(\d)}^{\e,\g,u_\e},\hat{Y}_{s}^{\e,\g,u_\e})\|^{2}\dif s\\
&\leq&\frac{2L_{b_2,\s_2}}{\g}\mE\int_{0}^{t}e^{\l s}|Y_{s}^{\e,\g,u_\e}-\hat{Y}_{s}^{\e,\g,u_\e}|^2\dif s+\frac{C}{\g}\mE\int_{0}^{t}e^{\l s}|X_{s}^{\e,\g,u_\e}-X_{s(\d)}^{\e,\g,u_\e}|^{2}\dif s.
\de

Collecting the above deduction, we conclude that
\ce
\mE|Y_{t}^{\e,\g,u_\e}-\hat{Y}_{t}^{\e,\g,u_\e}|^{2}e^{\l t}
&\leq&(\l-\frac{\b}{\g}+\frac{\a}{2\g}+\frac{2L_{b_2,\s_2}}{\g})\mE\int_{0}^{t}e^{\l s}|Y_{s}^{\e,\g,u_\e}-\hat{Y}_{s}^{\e,\g,u_\e}|^{2}\dif s\\
&&+\frac{C}{\e}\mE\int_{0}^{t}e^{\l s}|u_\e(s)|^2\dif s+\frac{C}{\g}\mE\int_{0}^{t}e^{\l s}|X_{s}^{\e,\g,u_\e}-X_{s(\d)}^{\e,\g,u_\e}|^{2}\dif s\\
&\leq&(\l-\frac{\b}{\g}+\frac{\a}{2\g}+\frac{2L_{b_2,\s_2}}{\g})\mE\int_{0}^{t}e^{\l s}|Y_{s}^{\e,\g,u_\e}-\hat{Y}_{s}^{\e,\g,u_\e}|^{2}\dif s\\
&&+\frac{C}{\e}\mE\int_{0}^{t}e^{\l s}|u_\e(s)|^2\dif s+\frac{C}{\g}\left(\sup\limits_{s\in[0,T]}\mE\sup _{s \leqslant r \leqslant s+\d}|X_{r}^{\e,\g,u_\e}-X_{s}^{\e,\g,u_\e}|^{2}\right)\frac{e^{\l t}-1}{\l},
\de
which together with $\l=\frac{\a}{2\g}$ yields that
\ce
\mE|Y_{t}^{\e,\g,u_\e}-\hat{Y}_{t}^{\e,\g,u_\e}|^{2}\leq \frac{C}{\e}\mE\int_{0}^{t}e^{-\frac{\a}{2\g}(t-s)}|u_\e(s)|^2\dif s+\frac{C}{\a}\left(\sup\limits_{s\in[0,T]}\mE\sup _{s \leqslant r \leqslant s+\d}|X_{r}^{\e,\g,u_\e}-X_{s}^{\e,\g,u_\e}|^{2}\right),
\de
and
\ce
\int_0^T\mE|Y_{t}^{\e,\g,u_\e}-\hat{Y}_{t}^{\e,\g,u_\e}|^{2}\dif t&\leq&\frac{C}{\e}\mE\int_0^T\int_{0}^{t}e^{-\frac{\a}{2\g}(t-s)}|u_\e(s)|^2\dif s\dif t\\
&&+\frac{CT}{\a}\left(\sup\limits_{s\in[0,T]}\mE\sup _{s \leqslant r \leqslant s+\d}|X_{r}^{\e,\g,u_\e}-X_{s}^{\e,\g,u_\e}|^{2}\right)\\
&\leq&C\left(\frac{\g}{\e}\right)+\frac{CT}{\a}\left(\sup\limits_{s\in[0,T]}\mE\sup _{s \leqslant r \leqslant s+\d}|X_{r}^{\e,\g,u_\e}-X_{s}^{\e,\g,u_\e}|^{2}\right).
\de
The proof is complete.
\end{proof}

\subsubsection{Verification for Condition \ref{cond}}\label{just}

\bl\label{auxilemm2}
Suppose that the assumptions of Theorem \ref{ldpmmsde} hold. Assume that $u,\ u_{\epsilon}\in\mathbf{A}_{2}^{N},\ \epsilon\in(0,1)$, and $u_{\epsilon}$ converges to $u$ almost surely as $\epsilon\rightarrow 0$. Then $\Psi^{\epsilon}(\sqrt{\epsilon}W+\int_{0}^{\cdot}u_{\epsilon}(s)\dif s)\rightarrow \Psi^{0}(\int_{0}^{\cdot}u(s)\dif s)$ in probability.
\el
\begin{proof}
We divide the proof into four steps.

{\bf Step 1.} We estimate the difference of $\Psi^{\epsilon}(\sqrt{\epsilon}W+\int_{0}^{\cdot}u_{\epsilon}(s)\dif s), \Psi^{0}(\int_{0}^{\cdot}u(s)\dif s)$ in $\mS$.

Note that
$$
X^{\epsilon,\g,u_{\epsilon}}=\Psi^{\epsilon}(\sqrt{\epsilon}W+\int_{0}^{\cdot}u_{\epsilon}(s)\dif s), \quad \bar{X}^{u}=\Psi^{0}(\int_{0}^{\cdot}u(s)\dif s).
$$
Thus, set $Z^{\epsilon,u_{\epsilon}}(t)=X^{\epsilon,\g,u_{\epsilon}}_{t}-\bar{X}^{u}_{t}$, and it holds that
\ce
Z^{\epsilon,u_{\epsilon}}(t)
&=&-(K^{1,\epsilon,\g,u_{\epsilon}}_{t}-\bar{K}^{u}_{t})+\int_{0}^{t}\left[b_1(X^{\epsilon,\g,u_{\epsilon}}_{s},Y^{\epsilon,\g,u_{\epsilon}}_{s})-\bar{b}_1(\bar{X}^{u}_{s})\right]\dif s\no\\
&&+\int_{0}^{t}\left[\s_1(X^{\epsilon,\g,u_{\epsilon}}_{s})\pi_1u_{\epsilon}(s)-\sigma_1(\bar{X}^{u}_{s})\pi_1u(s)\right]\dif s+\sqrt{\epsilon}\int_{0}^{t}\s_1(X^{\epsilon,\g,u_{\epsilon}}_{s})\dif W^1_s.
\de
By It\^o's formula and Lemma \ref{equi}, we get that
\be
&&|Z^{\epsilon,u_{\epsilon}}(t)|^{2}\no\\
&=&-2\int_{0}^{t}\langle Z^{\epsilon,u_{\epsilon}}(s),\dif (K^{1,\epsilon,\g,u_{\epsilon}}_{s}-\bar{K}^{u}_{s})   \rangle +2\int_{0}^{t}\langle   Z^{\epsilon,u_{\epsilon}}(s), b_1(X^{\epsilon,\g,u_{\epsilon}}_{s},Y^{\epsilon,\g,u_{\epsilon}}_{s})-\bar{b}_1(\bar{X}^{u}_{s})\rangle  \dif s   \no\\
&&+2\int_{0}^{t}\langle   Z^{\epsilon,u_{\epsilon}}(s), \s_1(X^{\epsilon,\g,u_{\epsilon}}_{s})\pi_1u_{\epsilon}(s)-\sigma_1(\bar{X}^{u}_{s})\pi_1u(s) \rangle  \dif s  \no\\
&&+2\sqrt{\epsilon} \int_{0}^{t}\langle  Z^{\epsilon,u_{\epsilon}}(s),  \s_1(X^{\epsilon,\g,u_{\epsilon}}_{s})\dif W^1_s\rangle +\epsilon\int_{0}^{t}\|\s_1(X^{\epsilon,\g,u_{\epsilon}}_{s})\|^{2} \dif s\no\\
&\leq&2\int_{0}^{t}\langle   Z^{\epsilon,u_{\epsilon}}(s), b_1(X^{\epsilon,\g,u_{\epsilon}}_{s},Y^{\epsilon,\g,u_{\epsilon}}_{s})-\bar{b}_1(\bar{X}^{u}_{s}) \rangle  \dif s   \no\\
&&+2\int_{0}^{t}\langle   Z^{\epsilon,u_{\epsilon}}(s), \s_1(X^{\epsilon,\g,u_{\epsilon}}_{s})\pi_1u_{\epsilon}(s)-\sigma_1(\bar{X}^{u}_{s})\pi_1u(s) \rangle  \dif s  \no\\
&&+2\sqrt{\epsilon} \int_{0}^{t}\langle  Z^{\epsilon,u_{\epsilon}}(s),  \s_1(X^{\epsilon,\g,u_{\epsilon}}_{s}) \dif W^1_s\rangle +\epsilon\int_{0}^{t}\|\s_1(X^{\epsilon,\g,u_{\epsilon}}_{s})\|^{2} \dif s\no\\
&=:&J_{1}(t)+J_{2}(t)+J_{3}(t)+J_{4}(t).
\label{j1j2j3j4}
\ee

For $J_{1}(t)$, note that
\ce
J_{1}(t)&=&2\int_{0}^{t}\langle   Z^{\epsilon,u_{\epsilon}}(s), b_1(X^{\epsilon,\g,u_{\epsilon}}_{s},Y^{\epsilon,\g,u_{\epsilon}}_{s})-b_1(X^{\epsilon,\g,u_{\epsilon}}_{s(\d)},\hat{Y}^{\epsilon,\g,u_{\epsilon}}_{s}) \rangle  \dif s\\
&&+2\int_{0}^{t}\langle   Z^{\epsilon,u_{\epsilon}}(s),-\bar{b}_1(X^{\epsilon,\g,u_{\epsilon}}_{s})+\bar{b}_1(X^{\epsilon,\g,u_{\epsilon}}_{s(\d)})\rangle \dif s\\
&&+2\int_{0}^{t}\langle   Z^{\epsilon,u_{\epsilon}}(s),\bar{b}_1(X^{\epsilon,\g,u_{\epsilon}}_{s})-\bar{b}_1(\bar{X}^{u}_{s}) \rangle  \dif s\\
&&+2\int_{0}^{t}\langle   Z^{\epsilon,u_{\epsilon}}(s)-Z^{\epsilon,u_{\epsilon}}(s(\d)),b_1(X^{\epsilon,\g,u_{\epsilon}}_{s(\d)},\hat{Y}^{\epsilon,\g,u_{\epsilon}}_{s})-\bar{b}_1(X^{\epsilon,\g,u_{\epsilon}}_{s(\d)})\rangle \dif s\\
&&+2\int_{0}^{t}\langle   Z^{\epsilon,u_{\epsilon}}(s(\d)),b_1(X^{\epsilon,\g,u_{\epsilon}}_{s(\d)},\hat{Y}^{\epsilon,\g,u_\e}_{s})-\bar{b}_1(X^{\epsilon,\g,u_{\epsilon}}_{s(\d)})\rangle \dif s\\
&=:&J_{11}(t)+J_{12}(t)+J_{13}(t)+J_{14}(t)+J_{15}(t).
\de
So, by the H\"older inequality and the Lipschitz continuity of $b_1, \bar{b}_1$, we get that
\be
&&\mE\left(\sup\limits_{t\in[0,T]}|J_{11}(t)|\right)+\mE\left(\sup\limits_{t\in[0,T]}|J_{12}(t)|\right)+\mE\left(\sup\limits_{t\in[0,T]}|J_{13}(t)|\right)\no\\
&\leq& C\int_{0}^{T}\mE|Z^{\epsilon,u_{\epsilon}}(s)|^2\dif s+C\int_{0}^{T}\mE|X^{\epsilon,\g,u_{\epsilon}}_{s}-X^{\epsilon,\g,u_{\epsilon}}_{s(\d)}|^2\dif s\no\\
&&+C\int_{0}^{T}\mE|Y^{\epsilon,\g,u_{\epsilon}}_{s}-\hat{Y}^{\epsilon,\g,u_{\epsilon}}_{s}|^2\dif s.
\label{j111213}
\ee
And the H\"older inequality and the linear growth of $b_1, \bar{b}_1$ imply that
\be
\mE\left(\sup\limits_{t\in[0,T]}|J_{14}(t)|\right)
&\leq& C\left(\int_0^T(\mE|X^{\epsilon,\g,u_{\epsilon}}_{s}-X^{\epsilon,\g,u_{\epsilon}}_{s(\d)}|^2+\mE|\bar{X}^{u}_{s}-\bar{X}^{u}_{s(\d)}|^2)\dif s\right)^{1/2}\no\\
&&\times\left(\int_0^T(1+\mE|X^{\epsilon,\g,u_{\epsilon}}_{s(\d)}|^2+\mE|\hat{Y}^{\epsilon,\g,u_{\epsilon}}_{s}|^2)\dif s\right)^{1/2}.
\label{j14}
\ee

For $J_{2}$, we rewrite it as
\ce
J_{2}(t)&=&2\int_{0}^{t}\langle   Z^{\epsilon,u_{\epsilon}}(s), \s_1(X^{\epsilon,\g,u_{\epsilon}}_{s})\pi_1u_{\epsilon}(s)-\sigma_1(\bar{X}^{u}_{s})\pi_1u_\e(s) \rangle  \dif s\no\\
&&+2\int_{0}^{t}\langle   Z^{\epsilon,u_{\epsilon}}(s), \sigma_1(\bar{X}^{u}_{s})\pi_1(u_{\epsilon}(s)-u(s)) \rangle  \dif s.\no\\
&=:&J_{21}(t)+J_{22}(t).
\de
For $J_{21}(t)$, by $(\mathbf{H}^1_{b_{1}, \s_{1}})$ and the H\"older inequality, it holds that
\be
\mE\sup\limits_{t\in[0,T]}|J_{21}(t)|&\leq& 2\sqrt{L_{b_1,\s_1}}\mE\int_{0}^{T}|Z^{\epsilon,u_{\epsilon}}(s)||Z^{\epsilon,u_{\epsilon}}(s)||u_{\epsilon}(s)|\dif s\no\\
&\leq&\frac{1}{4}\mE\left[\sup\limits_{t\in[0,T]}|Z^{\epsilon,u_{\epsilon}}(t)|^{2}\right]+C\mE\left(\int_{0}^{T}|Z^{\epsilon,u_{\epsilon}}(s)||u_{\epsilon}(s)|\dif s\right)^2\no\\
&\leq& \frac{1}{4}\mE\left[\sup\limits_{t\in[0,T]}|Z^{\epsilon,u_{\epsilon}}(t)|^{2}\right]+C\mE\left(\int_{0}^{T}|Z^{\epsilon,u_{\epsilon}}(s)|^2\dif s\right)\left(\int_{0}^{T}|u_{\epsilon}(s)|^{2}\dif s\right)\no\\
&\leq& \frac{1}{4}\mE\left[\sup\limits_{t\in[0,T]}|Z^{\epsilon,u_{\epsilon}}(t)|^{2}\right]+C\mE\left[\int_{0}^{T}|Z^{\epsilon,u_{\epsilon}}(s)|^{2}\dif s\right].
\label{j21}
\ee

For $J_{3}(t)$, from the Burkholder-Davis-Gundy inequality and the linear growth of $\s_1$, it follows that
\be
\mE\left(\sup\limits_{t\in[0,T]}|J_{3}(t)|\right)&\leq& 2\sqrt{\epsilon}C\mE\left(\int_{0}^{T}|Z^{\epsilon,u_{\epsilon}}(s)|^2\|\s_1(X^{\epsilon,\g,u_{\epsilon}}_{s}) \|^2 \dif s\right)^{1/2}\no\\
&\leq&\frac{1}{4}\mE\left[\sup\limits_{t\in[0,T]}|Z^{\epsilon,u_{\epsilon}}(t)|^{2}\right]+\e C\mE\int_{0}^{T} \|\s_1(X^{\epsilon,\g,u_{\epsilon}}_{s}) \|^2 \dif s\no\\
&\leq&\frac{1}{4}\mE\left[\sup\limits_{t\in[0,T]}|Z^{\epsilon,u_{\epsilon}}(t)|^{2}\right]+\e C\int_{0}^{T}(1+\mE|X^{\epsilon,\g,u_{\epsilon}}_{s}|^2)\dif s\no\\
&\leq&\frac{1}{4}\mE\left[\sup\limits_{t\in[0,T]}|Z^{\epsilon,u_{\epsilon}}(t)|^{2}\right]+\e C.
\label{j3}
\ee
For $J_{4}(t)$, by the linear growth of $\s_1$, we know
\be
\mE\left(\sup\limits_{t\in[0,T]}|J_{4}(t)|\right)\leq \e C\int_{0}^{T}(1+\mE|X^{\epsilon,\g,u_{\epsilon}}_{s}|^2)\dif s \leq\epsilon C.
\label{j4}
\ee

Combining (\ref{j111213})-(\ref{j4}) with (\ref{j1j2j3j4}), we can get
\ce
&&\mE\left(\sup\limits_{t\in[0,T]}|Z^{\epsilon,u_{\epsilon}}(t)|^{2}\right)\\
&\leq& C\int_{0}^{T}\mE\sup\limits_{r\in[0,s]}|Z^{\epsilon,u_{\epsilon}}(r)|^2\dif s+C\int_{0}^{T}\mE|X^{\epsilon,\g,u_{\epsilon}}_{s}-X^{\epsilon,\g,u_{\epsilon}}_{s(\d)}|^2\dif s+C\int_{0}^{T}\mE|Y^{\epsilon,\g,u_{\epsilon}}_{s}-\hat{Y}^{\epsilon,\g,u_{\epsilon}}_{s}|^2\dif s\\
&&+C\left(\int_0^T(\mE|X^{\epsilon,\g,u_{\epsilon}}_{s}-X^{\epsilon,\g,u_{\epsilon}}_{s(\d)}|^2+\mE|\bar{X}^{u}_{s}-\bar{X}^{u}_{s(\d)}|^2)\dif s\right)^{1/2}\\
&&\qquad\qquad \times \left(\int_0^T(1+\mE|X^{\epsilon,\g,u_{\epsilon}}_{s(\d)}|^2+\mE|\hat{Y}^{\epsilon,\g,u_{\epsilon}}_{s}|^2)\dif s\right)^{1/2}\\
&&+4\mE\left(\sup\limits_{t\in[0,T]}\left|J_{15}(t)\right|\right)+4\mE\left(\sup\limits_{t\in[0,T]}\left|J_{22}(t)\right|\right)+C\e,
\de
which together with the Gronwall inequality implies that
\be
\mE\left(\sup\limits_{t\in[0,T]}|Z^{\epsilon,u_{\epsilon}}(t)|^{2}\right)\leq C\bigg[\Gamma(\e)+\mE\left(\sup\limits_{t\in[0,T]}\left|J_{15}(t)\right|\right)+\mE\left(\sup\limits_{t\in[0,T]}\left|J_{22}(t)\right|\right)+C\e\bigg],
\label{zeue}
\ee
where
\ce
\Gamma(\e)&:=&C\int_{0}^{T}\mE|X^{\epsilon,\g,u_{\epsilon}}_{s}-X^{\epsilon,\g,u_{\epsilon}}_{s(\d)}|^2\dif s+C\int_{0}^{T}\mE|Y^{\epsilon,\g,u_{\epsilon}}_{s}-\hat{Y}^{\epsilon,\g,u_{\epsilon}}_{s}|^2\dif s\\
&&+C\left(\int_0^T(\mE|X^{\epsilon,\g,u_{\epsilon}}_{s}-X^{\epsilon,\g,u_{\epsilon}}_{s(\d)}|^2+\mE|\bar{X}^{u}_{s}-\bar{X}^{u}_{s(\d)}|^2)\dif s\right)^{1/2}\\
&&\qquad\qquad \times \left(\int_0^T(1+\mE|X^{\epsilon,\g,u_{\epsilon}}_{s(\d)}|^2+\mE|\hat{Y}^{\epsilon,\g,u_{\epsilon}}_{s}|^2)\dif s\right)^{1/2}.
\de

{\bf Step 2.} We estimate $\mE\left(\sup\limits_{t\in[0,T]}\left|J_{15}(t)\right|\right)$.

For $J_{15}(t)$, it holds that
\ce
\mE\left(\sup\limits_{t\in[0,T]}\left|J_{15}(t)\right|\right)&=&2\Bigg(\mE\sup_{0\leq t\leq T}
 \Big|\int_{0}^{[\frac{t}{\d}]\d}\<Z^{\epsilon,u_{\epsilon}}(s(\d)),b_{1}(X_{s(\d)}^{\e,\g,u_\e},\hat{Y}_{s}^{\e,\g,u_\e})
 -\bar{b}_{1}(X_{s(\d)}^{\e,\g,u_\e})\>\dif s\no\\
 &&\quad\quad\quad\quad+\int_{[\frac{t}{\d}]\d}^{t}\<Z^{\epsilon,u_{\epsilon}}(s(\d)),b_{1}(X_{s(\d)}^{\e,\g,u_\e},\hat{Y}_{s}^{\e,\g,u_\e})
 -\bar{b}_{1}(X_{s(\d)}^{\e,\g,u_\e})\>\dif s\Big|\Bigg)\no\\
&\leq&2\Bigg(\mE\sup_{0\leq t\leq T}
\Big|\int_{0}^{[\frac{t}{\d}]\d}\<Z^{\epsilon,u_{\epsilon}}(s(\d)),b_{1}(X_{s(\d)}^{\e,\g,u_\e},\hat{Y}_{s}^{\e,\g,u_\e})
 -\bar{b}_{1}(X_{s(\d)}^{\e,\g,u_\e})\>\dif s\Big|\Bigg)\no\\
&&+2\Bigg(\mE\sup_{0\leq t\leq T}
\Big|\int_{[\frac{t}{\d}]\d}^{t}\<Z^{\epsilon,u_{\epsilon}}(s(\d)),b_{1}(X_{s(\d)}^{\e,\g,u_\e},\hat{Y}_{s}^{\e,\g,u_\e})
 -\bar{b}_{1}(X_{s(\d)}^{\e,\g,u_\e})\>\dif s\Big|\Bigg)\no\\
 &=:&J_{151}+J_{152}.
\de

Next, we estimate $J_{151}$. By the similar deduction to that of $I_{21}$ in Lemma \ref{2orde}, we get that
\be
J_{151}\leq C(\frac{\g}{\d})^{1/2}.
\label{b4de}
\ee

Finally, we estimate $J_{152}$. By (\ref{b1line}), (\ref{xeub}), (\ref{hatzub}) and the H\"older inequality, one could get that
\ce
J_{152}
&\leq& 2\Bigg(\mE\sup_{0\leq t\leq T}
\int_{[\frac{t}{\d}]\d}^{t}|Z^{\epsilon,u_{\epsilon}}(s(\d))||b_{1}(X_{s(\d)}^{\e,\g,u_\e},\hat{Y}_{s}^{\e,\g,u_\e})
 -\bar{b}_{1}(X_{s(\d)}^{\e,\g,u_\e})|\dif s\Bigg)\\
 &\leq&2\d^{1/2}\left(\mE\sup_{0\leq s\leq T}|Z^{\epsilon,u_{\epsilon}}(s)|^2\right)^{1/2}\left(\mE\sup_{0\leq t\leq T}\int_{[\frac{t}{\d}]\d}^{t}|b_{1}(X_{s(\d)}^{\e,\g,u_\e},\hat{Y}_{s}^{\e,\g,u_\e})
 -\bar{b}_{1}(X_{s(\d)}^{\e,\g,u_\e})|^2\dif s\right)^{1/2}\\ 
 &\leq&2\d^{1/2}\left(\mE\sup_{0\leq s\leq T}|Z^{\epsilon,u_{\epsilon}}(s)|^2\right)^{1/2}\left(\mE\int_{0}^{T}|b_{1}(X_{s(\d)}^{\e,\g,u_\e},\hat{Y}_{s}^{\e,\g,u_\e})
 -\bar{b}_{1}(X_{s(\d)}^{\e,\g,u_\e})|^2\dif s\right)^{1/2}\\
 &\leq&C\d^{1/2}\left(\mE\sup_{0\leq s\leq T}|X^{\epsilon,\g,u_{\epsilon}}_{s}|^2+\mE\sup_{0\leq s\leq T}|\bar{X}^{u}_{s}|^2\right)^{1/2}\left(\int_{0}^{T}(1+\mE|X_{s(\d)}^{\e,\g,u_\e}|^2+\mE|\hat{Y}_{s}^{\e,\g,u_\e}|^2)\dif s\right)^{1/2}\\
&\leq& C\d^{1/2},
\de
which together with (\ref{b4de}) implies
\be
\mE\left(\sup\limits_{t\in[0,T]}\left|J_{15}(t)\right|\right)\leq C\((\frac{\g}{\d})^{1/2}+\d^{1/2}\).
\label{j15}
\ee

{\bf Step 3.} We compute $\mE\left(\sup\limits_{t\in[0,T]}\left|J_{22}(t)\right|\right)$.

For $J_{22}(t)$, applying the It\^o formula to $\<Z^{\epsilon,u_{\epsilon}}(t),g_\e(t)\>$, we have that
\ce
\frac{1}{2}J_{22}(t)&=&\<Z^{\epsilon,u_{\epsilon}}(t),g_\e(t)\>+\int_0^t\<g_\e(s),\dif (K^{1,\epsilon,\g,u_{\epsilon}}_{s}-\bar{K}^{u}_{s})\>\\
&&-\int_0^t\<g_\e(s),b_1(X^{\epsilon,\g,u_{\epsilon}}_{s},Y^{\epsilon,\g,u_{\epsilon}}_{s})-\bar{b}_1(\bar{X}^{u}_{s})\>\dif s\\
&&-\int_0^t\<g_\e(s),\s_1(X^{\epsilon,\g,u_{\epsilon}}_{s})\pi_1u_{\epsilon}(s)-\sigma_1(\bar{X}^{u}_{s})\pi_1u(s)\>\dif s\\
&&-\sqrt{\epsilon}\int_0^t\<g_\e(s),\s_1(X^{\epsilon,\g,u_{\epsilon}}_{s})\dif W^1_s\>\\
&=:&J_{221}(t)+J_{222}(t)+J_{223}(t)+J_{224}(t)+J_{225}(t).
\de

For $J_{221}(t)$, note that 
\ce
\sup\limits_{t\in[0,T]}|J_{221}(t)|\leq\sup\limits_{t\in[0,T]}|Z^{\epsilon,u_{\epsilon}}(t)|\sup\limits_{t\in[0,T]}|g_\e(t)|.
\de
Thus, by Lemma \ref{ge} it holds that
\ce
\sup\limits_{t\in[0,T]}|J_{221}(t)|\overset{a.s.}{\rightarrow} 0, \quad \e\rightarrow 0.
\de
Besides, the fact that $\sup\limits_{t\in[0,T]}|g_\e(t)|\leq C(1+|x_0|)$ implies that
\ce
\mE\sup\limits_{t\in[0,T]}|Z^{\epsilon,u_{\epsilon}}(t)|\sup\limits_{t\in[0,T]}|g_\e(t)|&\leq& C(1+|x_0|)\mE\sup\limits_{t\in[0,T]}|Z^{\epsilon,u_{\epsilon}}(t)|\\
&\leq& C(1+|x_0|)\left(\mE\sup\limits_{t\in[0,T]}|Z^{\epsilon,u_{\epsilon}}(t)|^2\right)^{1/2}\\
&\leq& C(1+|x_0|)\left(\mE\sup_{0\leq t\leq T}|X^{\epsilon,\g,u_{\epsilon}}_{t}|^2+\mE\sup_{0\leq t\leq T}|\bar{X}^{u}_{t}|^2\right)^{1/2}\\
&\leq&C(1+|x_0|^2+|y_0|^2),
\de
where the last inequality is based on (\ref{xeub}), (\ref{barxub}). Then the dominated convergence theorem yields that
\ce
\lim\limits_{\e\rightarrow 0}\mE\sup\limits_{t\in[0,T]}|J_{221}(t)|=0.
\de

For $J_{222}(t)$, noticing that
\ce
\sup\limits_{t\in[0,T]}|J_{222}(t)|\leq\sup\limits_{t\in[0,T]}|g_\e(t)|(|K^{1,\epsilon,\g,u_{\epsilon}}|_0^T+|\bar{K}^{u}|_0^T),
\de
by Lemma \ref{ge} we obtain that $\sup\limits_{t\in[0,T]}|J_{222}(t)|\overset{a.s.}{\rightarrow} 0$. Then (\ref{keub}), (\ref{barkub}) and the dominated convergence theorem imply that
\be
\lim\limits_{\e\rightarrow 0}\mE\sup\limits_{t\in[0,T]}|J_{222}(t)|=0.
\label{j222}
\ee

By the same deduction to that for (\ref{j222}), one can get that
\ce
\lim\limits_{\e\rightarrow 0}\mE\sup\limits_{t\in[0,T]}|J_{223}(t)|=0, \quad \lim\limits_{\e\rightarrow 0}\mE\sup\limits_{t\in[0,T]}|J_{224}(t)|=0.
\de

For $J_{225}(t)$, the BDG inequality yields that
\ce
\mE\sup\limits_{t\in[0,T]}|J_{225}(t)|&\leq& C\sqrt{\epsilon}\mE\left(\int_0^T|g_\e(s)|^2\|\s_1(X^{\epsilon,\g,u_{\epsilon}}_{s})\|^2\dif s\right)^{1/2}\\
&\leq& C\sqrt{\epsilon}\left(\mE\sup\limits_{t\in[0,T]}|g_\e(t)|^2+\mE\int_0^T\|\s_1(X^{\epsilon,\g,u_{\epsilon}}_{s})\|^2\dif s\right)\\
&\leq& C\sqrt{\epsilon}(1+|x_0|^2+|y_0|^2).
\de
Therefore, it holds that
\ce
\lim\limits_{\e\rightarrow 0}\mE\sup\limits_{t\in[0,T]}|J_{225}(t)|=0.
\de

Combining the above deduction, we get that 
\be
\lim\limits_{\e\rightarrow 0}\mE\sup\limits_{t\in[0,T]}\left|J_{22}(t)\right|=0.
\label{j22}
\ee

{\bf Step 4.} We prove that $\Psi^{\epsilon}(\sqrt{\epsilon}W+\int_{0}^{\cdot}u_{\epsilon}(s)\dif s)\rightarrow \Psi^{0}(\int_{0}^{\cdot}u(s)\dif s)$ in probability.

Let $\d=\g^{\iota}$ for $0<\iota<1$ and $\g\rightarrow 0, \d\rightarrow 0, \g/\d\rightarrow 0$ as $\e\rightarrow 0$. Thus, (\ref{xegutse}), (\ref{barxuts}), (\ref{unztu}), (\ref{zeue}), (\ref{j15}), (\ref{j22}) imply that $X^{\epsilon,\g,u_{\epsilon}}-\bar{X}^{u}\rightarrow 0$ in the mean square and $\Psi^{\epsilon}(\sqrt{\epsilon}W+\int_{0}^{\cdot}u_{\epsilon}(s)\dif s)\rightarrow \Psi^{0}(\int_{0}^{\cdot}u(s)\dif s)$ in probability. The proof is complete.
\end{proof}

By the similar or even simpler deduction to that of Lemma \ref{auxilemm2}, we obtain the following result.

\bl\label{auxilemm3}
Suppose that the assumptions of Theorem \ref{ldpmmsde} hold, and $h_{\e}\rightarrow h$ in $\mathbf{D}_2^{N}$ as ${\e}\rightarrow0$. Then $\Psi^{0}(\int_{0}^{\cdot}h_{\e}(s)\dif s)$ converges to $\Psi^{0}(\int_{0}^{\cdot}h(s)\dif s)$.
\el

Now, it is the position to prove Theorem \ref{ldpmmsde}.

{\bf Proof of Theorem \ref{ldpmmsde}.}

By Theorem \ref{ldpbase}, to establish LDP, it is sufficient to verify the two conditions in Condition \ref{cond}.  
In Lemma \ref{auxilemm3}, we have already proved Condition \ref{cond} $(i)$. Thus we only need to verify Condition \ref{cond} $(ii)$. 

For $\epsilon\in(0,1)$ and $\{u_{\epsilon}, \e>0\}\subset\mathbf{A}_{2}^{N}$, $u\in\mathbf{A}_{2}^{N}$, let $u_{\epsilon}$ converge to $u$ in distribution. By the Skorohod theorem, there exists a
probability space $(\tilde{\Omega}, \tilde{\sF}, \tilde{\mP})$, and ${\bf A}_2^N$-valued random variables $\{\tilde{u}_{\epsilon}\}$, $\tilde{u}$ and a $d_1+d_2$-dimensional Brownian motion $\tilde{W}$ such that

(i) $(\tilde{u}_{\epsilon},\tilde{W})\overset{d}{=}(u_{\epsilon},W)$ and $\tilde{u}\overset{d}{=}u$;

(ii) $\tilde{u}_{\epsilon}$ converges to $\tilde{u}$ almost surely.

In the following, we construct multivalued SDEs:
\be\left\{\begin{array}{l}
\dif X_{t}^{\e,\g,\tilde{u}_\e}\in -A_1(X_{t}^{\e,\g,\tilde{u}_\e})\dif t+b_{1}(X_{t}^{\e,\g,\tilde{u}_\e},Y_{t}^{\e,\g,\tilde{u}_\e})\dif t+\sigma_1(X^{\epsilon,\g,\tilde{u}_\e}_{t})\pi_1 \tilde{u}_\e(t)\dif t\\
\qquad\qquad+\sqrt{\e}\s_{1}(X_{t}^{\e,\g,\tilde{u}_\e})\dif \tilde{W}^1_{t},\\
X_{0}^{\e,\g,\tilde{u}_\e}=x_0\in\overline{\cD(A_1)},\quad  0\leq t\leq T,\\
\dif Y_{t}^{\e,\g,\tilde{u}_\e}\in -A_2(Y_{t}^{\e,\g,\tilde{u}_\e})\dif t+\frac{1}{\g}b_{2}(X_{t}^{\e,\g,\tilde{u}_\e},Y_{t}^{\e,\g,\tilde{u}_\e})\dif t+\frac{1}{\sqrt{\g \e}}\sigma_2(X^{\epsilon,\g,\tilde{u}_\e}_{t},Y^{\epsilon,\g,\tilde{u}_\e}_{t})\pi_2 \tilde{u}_\e(t)\dif t\\
\qquad\qquad +\frac{1}{\sqrt{\g}}\s_{2}(X_{t}^{\e,\g,\tilde{u}_\e},Y_{t}^{\e,\g,\tilde{u}_\e})\dif \tilde{W}^2_{t},  \quad \tilde{u}_\e\in\mathbf{A}_2^{N},\\
Y_{0}^{\e,\g,\tilde{u}_\e}=y_0\in\overline{\cD(A_2)},\quad  0\leq t\leq T,
\end{array}
\right.
\label{contproc2}
\ee
\be\left\{\begin{array}{l}
\dif\bar{X}^{\tilde{u}}_{t}\in -A_1(\bar{X}^{\tilde{u}}_{t})\dif t+\bar{b}_{1}(\bar{X}^{\tilde{u}}_{t})\dif t+\s_{1}(\bar{X}^{\tilde{u}}_{t})\pi_1\tilde{u}(t)\dif t, \quad \tilde{u}\in\mathbf{A}_2^{N},\\
\bar{X}^{\tilde{u}}_{0}=x_0.
\end{array}
\right.
\label{deteequa2}
\ee
For the system (\ref{contproc2}) and Eq.(\ref{deteequa2}), each has a unique strong solution denoted by $(X^{\epsilon,\g,\tilde{u}_\e}, \\K^{1,\epsilon,\g,\tilde{u}_\e},Y^{\epsilon,\g,\tilde{u}_\e}, K^{2,\epsilon,\g,\tilde{u}_\e})$ and $(\bar{X}^{\tilde{u}},\bar{K}^{\tilde{u}})$ respectively. Moreover, it holds that
$$
X^{\epsilon,\g,\tilde{u}_{\epsilon}}=\Psi^{\epsilon}\left(\sqrt{\epsilon}\tilde{W}+\int_{0}^{\cdot}\tilde{u}_{\epsilon}(s)\dif s\right), \quad \bar{X}^{\tilde{u}}=\Psi^{0}\left(\int_{0}^{\cdot}\tilde{u}(s)\dif s\right).
$$
By Lemma \ref{auxilemm2}, we have that $\Psi^{\epsilon}(\sqrt{\epsilon}\tilde{W}+\int_{0}^{\cdot}\tilde{u}_{\epsilon}(s)\dif s)\rightarrow \Psi^{0}(\int_{0}^{\cdot}\tilde{u}(s)\dif s)$ in probability, which yields that
$\Psi^{\epsilon}(\sqrt{\epsilon}\tilde{W}+\int_{0}^{\cdot}\tilde{u}_{\epsilon}(s)\dif s)\rightarrow \Psi^{0}(\int_{0}^{\cdot}\tilde{u}(s)\dif s)$ in distribution. Note that
\ce
\Psi^{\epsilon}\left(\sqrt{\epsilon}\tilde{W}+\int_{0}^{\cdot}\tilde{u}_{\epsilon}(s)\dif s\right)&\overset{d}{=}&\Psi^{\epsilon}\left(\sqrt{\epsilon}W+\int_{0}^{\cdot}u_{\epsilon}(s)\dif s\right),\\
\Psi^{0}\left(\int_{0}^{\cdot}\tilde{u}(s)\dif s\right)&\overset{d}{=}&\Psi^{0}\left(\int_{0}^{\cdot}u(s)\dif s\right).
\de
So, $\Psi^{\epsilon}\left(\sqrt{\epsilon}W+\int_{0}^{\cdot}u_{\epsilon}(s)\dif s\right)\rightarrow \Psi^{0}\left(\int_{0}^{\cdot}u(s)\dif s\right)$ in distribution, which is Condition \ref{cond} $(ii)$.

\subsection{Proof of Corollary \ref{charate}}

In this subsection, we prove Corollary \ref{charate}.

First of all, by Theorem \ref{ldpmmsde}, we know that the family $\{X^{\epsilon,\g},\epsilon\in(0,1)\}$ satisfies the LDP in $\mS:=C([0,T],\overline{\mathcal{D}(\p I_D)})$ with the rate function given by
$$
I(\varsigma)=\frac{1}{2} \inf\limits_{h\in {\bf D}_{\varsigma}: \varsigma=\bar{X}^{h}}\|h\|_{\mH}^2.
$$
Besides, we notice that $\g$ goes to $0$ faster than $\e$ and homogenization occurs first. Therefore, we can think of $I$ as the rate function with which the family $\{\bar{X}^{\e},\epsilon\in(0,1)\}$ satisfies the LDP in $\mS:=C([0,T],\overline{\mathcal{D}(\p I_D)})$, where $(\bar{X}^{\e},\bar{K}^{\e})$ is the unique strong solution to the following multivalued SDE:
\be\left\{\begin{array}{l}
\dif \bar{X}^\e_{t}= \triangledown \kappa(\bar{X}^\e_{t})\dif |\bar{K}^{\e}|_0^t+\bar{b}_{1}(\bar{X}^\e_{t})\dif t+\sqrt{\e}\s_{1}(\bar{X}^\e_{t})\dif W^1_t,\\
\bar{X}^\e_{0}=x_0\in\overline{\cD(\p I_D)}.
\end{array}
\right.
\label{homosyst}
\ee

Next, let $\phi\in C([0,T], D), \psi\in C([0,T], \mR^n), V\in\sV_0$ such that
\ce
\phi(t)=\psi(t)+V(t), \quad V_t=\int_0^t\triangledown \kappa(\phi_r)\dif |V|_0^r, \quad |V|_0^t=\int_0^tI_{\{\phi_r\in\p D\}}\dif |V|_0^r.
\de
For $\phi, \psi$ as above, we define a mapping $\Gamma: C([0,T], \mR^n)\rightarrow C([0,T], D)$ by
$$
\phi=\Gamma(\psi).
$$
Then by \cite{ys}, it holds that $\Gamma$ is continuous. By the contraction principle (\cite[Theorem 3.1, P.64]{fw1}) and a LDP for diffusion processes (\cite[Theorem 3.1, P.135]{fw1}), we know that 
\ce
I(\phi)=\left\{\begin{array}{l}\frac{1}{2}\inf\limits_{\psi:\phi=\Gamma(\psi)}\int_0^T\(\dot{\psi}(t)-\bar{b}_1(\phi(t))\)^*(\s_1\s_1^*)^{-1}(\phi(t))\(\dot{\psi}(t)-\bar{b}_1(\phi(t))\)\dif t,\\
 ~\mbox{if}~ \psi\in AC([0,T],\mR^n) ~\mbox{and}~ \psi(0)=x_0,\\
\infty, \quad otherwise.
\end{array}
\right.
\de
The proof is complete.

\section{An example}\label{exam}

In this section, we explain our results by an example.

\bx
The typical example that Theorem \ref{ldpmmsde} is applicable to is a particular system of slow-fast motion, where the fast motion is a fast mean reverting process and the slow motion appears due to the interest in short time asymptotics. Concretely speaking, consider the following system on $\mR\times\mR$:
\be\left\{\begin{array}{l}
\dif X_{t}^{\e,\g}\in -\p \varphi_1(X_{t}^{\e,\g})\dif t+b(Y_{t}^{\e,\g})\dif t+\sqrt{\e}\s(Y_{t}^{\e,\g})\dif W^1_{t},\\
X_{0}^{\e,\g}=x_0\in\overline{\cD(\p \varphi_1)},\quad  0\leq t\leq T,\\
\dif Y_{t}^{\e,\g}\in -\p \varphi_2(Y_{t}^{\e,\g})\dif t+\frac{1}{\g}(m-\frac{1}{2}Y_{t}^{\e,\g})\dif t+\frac{1}{\sqrt{\g}}\dif W^2_{t},\\
Y_{0}^{\e,\g}=y_0\in\overline{\cD(\p \varphi_2)},\quad  0\leq t\leq T,
\end{array}
\right.
\label{Eqex}
\ee
where $\varphi_1, \varphi_2$ are lower semicontinuous convex functions on $\mR$, $m\in\mR$, $W^1, W^2$ are $1$-dimensional Brownian motions  and $b: \mR\mapsto\mR, \s: \mR\mapsto\mR$ are measurable functions. It is obvious that
\ce
b_1(x,y)=b(y), \quad \s_1(x,y)=\s(y), \quad b_2(x,y)=m-\frac{1}{2}y, \quad \s_2(x,y)=1.
\de
If $0\in{\rm Int}(\cD(\p \varphi_1))$, $0\in{\rm Int}(\cD(\p \varphi_2))$, and $b, \s$ are Lipschitz continuous, by Theorem \ref{well}, we know that the system (\ref{Eqex}) has a unique strong solution $(X_{\cdot}^{\e,\g},K_{\cdot}^{1,\e,\g},Y_{\cdot}^{\e,\g},K_{\cdot}^{2,\e,\g})$. Note that for $y_{i}\in\mR$, $i=1, 2$
\ce
2\<y_1-y_2,b_2(x,y_1)-b_2(x,y_2)\>=2\<y_1-y_2,m-\frac{1}{2}y_1-m+\frac{1}{2}y_2\>=-|y_1-y_2|^2,
\de
where $\b=1, L_{b_2,\sigma_2}=\frac{1}{4}, \b>2L_{b_2,\sigma_2}$. Thus, when $\s(y)=\s$ for $\s\in\mR$, by Theorem \ref{ldpmmsde} we obtain that the family $\{X^{\epsilon,\g},\epsilon\in(0,1)\}$ satisfies the LDP in $\mS:=C([0,T],\overline{\mathcal{D}(\p \varphi_1)})$ with the rate function given by
$$
I(\varsigma)=\frac{1}{2} \inf\limits_{h\in {\bf D}_{\varsigma}: \varsigma=\bar{X}^{h}}\|h\|_{\mH}^2,
$$
where $\bar{X}^{h}$ solves the following equation
\ce\left\{\begin{array}{l}
\dif \bar{X}^h_{t}\in -\p \varphi_1(\bar{X}^h_{t})\dif t+\bar{b}\dif t+\s\pi_1h(t)\dif t,\\
\bar{X}^h_{0}=x_0\in\overline{\cD(\p \varphi_1)},
\end{array}
\right.
\de
$\bar{b}=\int_{\mR}b(y)\nu(\dif y)$ and $\nu$ is the unique invariant probability measure of the following equation
\ce\left\{\begin{array}{l}
\dif Y_{t}\in -\p \varphi_2(Y_{t})\dif t+(m-\frac{1}{2}Y_{t})\dif t+\dif W^2_{t},\\
Y_{0}=y_0\in\overline{\cD(\p \varphi_2)},\quad  0\leq t\leq T.
\end{array}
\right.
\de

Besides, if $\varphi_1=\varphi_2=0$, the system (\ref{Eqex}) falls into the framework in \cite[Section 6.2]{ks1}. Moreover, our LDP result is the same to the conclusion there.
\ex

\end{document}